\documentclass[12pt]{report}    %
\usepackage{hyperref}
\usepackage[x11names]{xcolor}
\hypersetup{colorlinks = true, allcolors=SteelBlue3}

\renewcommand{\baselinestretch}{1.5}\normalsize
\def\title#1{\def\ntitle{#1}}
\def\author#1{\def\nauthor{#1}}
\def\degree#1{\def\ndegree{#1}}
\def\date#1{\def\ndate{#1}}
\newcommand{\Thesistitle}{
\newpage
\null\vspace*{0.5in} 
\begin{center} 
{\Large\bf \ntitle}

\vspace*{3.5cm}
{\large\bf \nauthor}
\vspace*{2cm}

A thesis submitted in partial fulfillment of the requirements

for the degree of 

{\normalsize \ndegree}

\vspace*{3.5cm}
{\normalsize \ndate}

\vspace*{1cm}
Hong Kong Baptist University
\thispagestyle{empty}
\end{center}}

\renewenvironment{abstract}
{\newpage
\setcounter{page}{1}
\addcontentsline{toc}{chapter}{Abstract}
\begin{center}{\Large\bf Abstract}\normalsize
\end{center}
\vspace*{0.5cm}
\renewcommand{\baselinestretch}{1} \normalsize}
{\newpage}

\newcommand{\thesis}{
\newpage
\pagenumbering{arabic}
\setcounter{chapter}{0}
\setcounter{section}{0}
\setcounter{page}{1}
}

\makeatletter
\def\@chapter[#1]#2{\ifnum \c@secnumdepth >\m@ne
                         \refstepcounter{chapter}%
                         \typeout{\@chapapp\space\thechapter.}%
                         \addcontentsline{toc}{chapter}%
                                   {\protect\numberline{Chapter \thechapter}\hspace*{1.5cm} \ \ #1\ }%
                    \else
                      \addcontentsline{toc}{chapter}{#1}%
                    \fi
                    \chaptermark{#1}%
                    \addtocontents{lof}{\protect\addvspace{10\p@}}%
                    \addtocontents{lot}{\protect\addvspace{10\p@}}%
                    \if@twocolumn
                      \@topnewpage[\@makechapterhead{#2}]%
                    \else
                      \@makechapterhead{#2}%
                      \@afterheading
                    \fi}
\makeatother

\usepackage{amsmath}%
\usepackage{graphicx}
\graphicspath{{Images/}}
\usepackage{amssymb}

\renewcommand{\vec}[1]{\ensuremath{\mathbf{#1}}}
\newcommand{\vecsym}[1]{\ensuremath{\boldsymbol{#1}}}
\def\bbl{\text{\boldmath$\{$}}
\def\bbr{\text{\boldmath$\}$}}
\newcommand{\bbrace}[1]{\bbl #1 \bbr}

\newlength{\overwdth}

\DeclareMathOperator{\Order}{O}

\newcommand{\Zb}{\vec Z_{b}}
\newcommand{\Z}{\vec Z}

\DeclareMathOperator{\Err}{Err}
\newcommand{\tbfpsi}{\tilde{\vecsym \psi}}
\newcommand{\bfpsi}{\vecsym \psi} 
\newcommand{\bfphi}{\vecsym \phi} 
\newcommand{\bfa}{\vec a}
\newcommand{\bfA}{\vec A}

\newcommand{\bfc}{\vec c}
\newcommand{\bfd}{\vec d}
\newcommand{\bfe}{\vec e}
\newcommand{\bfy}{\vec y}
\newcommand{\bfz}{\vec z} 
\newcommand{\bfx}{\vec x} 
\newcommand{\bfC}{\vec C}

\newcommand{\bfL}{\vec L}
\newcommand{\bfl}{\vec l}
\newcommand{\bfh}{\vec h}
\newcommand{\bftheta}{\vecsym \theta}
\newcommand{\bfSigma}{\vecsym \Sigma}
\newcommand{\Cs}{[0,1)^{s}}

\newtheorem{theorem}{Theorem}

\begin{document}                %
\author{HONG, Hee Sun}    %

\title{Digital Nets and Sequences for Quasi-Monte Carlo Methods}
\degree{Doctor of Philosophy}   %
\date{May 2002}            %
\Thesistitle                    %

\renewcommand{\contentsname}{\centerline{Contents}}
\newpage
\renewcommand{\baselinestretch}{2}\normalsize

\newpage
\begin{abstract}                %
Quasi-Monte Carlo methods are a way of improving 
the efficiency of Monte Carlo methods. Digital nets and sequences 
are one of the low discrepancy point sets used
in quasi-Monte Carlo methods. This thesis presents the three new results
pertaining to digital nets and sequences: implementing 
randomized digital nets, finding the distribution of the discrepancy of 
scrambled digital nets, and obtaining better quality of digital nets 
through evolutionary computation.  
Finally, applications of scrambled and non-scrambled digital nets are 
provided.
\end{abstract}

\newpage          %
\addcontentsline{toc}{chapter}{Table of Contents}
\tableofcontents

\newpage         %
\addcontentsline{toc}{chapter}{List of Figures}
\listoffigures

\newpage         %
\addcontentsline{toc}{chapter}{List of Tables}
\listoftables

\thesis           %
\newpage
\begin{chapter}{Introduction}
We consider the problem of approximating $I$, the integral of
a function $f$ over the $s$-dimensional unit cube $\Cs$:
$$
I = \int_{\Cs} f(\bfy) d\bfy
$$
Computing such a multidimensional integral is important since
it has wide applications in finance \cite{CafMor96,PapTra96}, 
physics and engineering \cite{Kei96,PapTra97} and statistics 
\cite{Fan80,Gen92,Gen93}. Numerical integration methods, such as 
tensor product Newton Cotes and Gaussian quadrature,  
become impractical as $s$ increases since the error of the such 
methods converges like $\Order(N^{-p/s})$ where $p$ is the smoothness of 
function.

Monte Carlo methods are often a popular choice to estimate 
such integrals. They use the sample mean of the integrand evaluated
on a set, $P$, with $N$ randomly selected independent points 
drawn from a uniform distribution on $[0,1)^s:$ 
$$
\hat{I}_N =\frac{1}{N}\sum_{\bfy \in P} f(\bfy).
$$
Unfortunately the convergence rate of the quadrature error for Monte Carlo 
methods, $\Order(N^{-1/2})$, is relatively low.

To overcome the lower convergence rate of Monte Carlo methods, 
quasirandom sequences have been introduced.   
The basic quasirandom sequences replace a random set $P$ by a 
carefully chosen deterministic set that is more uniformly 
distributed on $[0,1)^s$, and it often yields a convergence rate of
$O(N^{-1}{\log}^{s-1}N)$. The quadrature methods based on low 
discrepancy sets are called quasi-Monte Carlo methods
\cite{Nie92,Tez95,HelLar98,Fox99}.

This chapter introduces the basic definitions of quasi-Monte Carlo 
methods and the outline of the remaining thesis.

\begin{section}{Quasi-Monte Carlo Methods}

The most important property of quasirandom sequences is 
an equidistribution property.
The quality of such uniformity for quasirandom sequences is
measured by the discrepancy, which is a distance between the 
continuous uniform distribution on $\Cs$ and the empirical  
distribution of the $\bfy_i$ for $i = 1,\cdots,N$. 
First we define the empirical distribution of the 
sequence as
\begin{equation}
F_N(\bfy) = \frac{1}{N} \sum_{i=1}^{N}1_{[\bfy_i,\infty)}(\bfy),
\end{equation}
where $1_A$ is the indicator function of the set $A$. 
We define the uniform distribution on $\Cs$ as
\begin{equation}
F_{\text{unif}}(\bfy) = \prod_{j=1}^{s} y_j, \mbox{   } 
\bfy = (y_1,\cdots,y_s)^T \in \Cs.
\end{equation}

The discrepancy arises from the worst case error analysis of 
quasi-Monte Carlo quadratures. The Koskma-Hlawka inequality, which 
provides an important theoretical justification for the 
quasi-Monte Carlo quadratures, uses the discrepancy to provide an upper
bound quadrature error of quasi-Monte Carlo methods: 
\begin{equation}
|\hat{I}_N - I| \le D^*(P)V_{HK}(f)
\end{equation}
where $V_{HK}$ is said to have bounded variation $f$ in the sense of Hardy and 
Krause. If $V_{HK}(f) < \infty $, then $f \in BVHK$. For further 
understanding of $V_{HK}$ see Niederreiter \cite{Nie92}.
$D^*(P)$ indicates the star discrepancy which will be defined next.

The star discrepancy is the most widely studied discrepancy and it 
is defined as follows:
\begin{equation}
D^*(P) = \sup_{\bfy \in \Cs} |F_{\text{unif}}(\bfy)-F_N(\bfy)|
= ||F_{\text{unif}}-F_N||_{\infty}.
\end{equation}

The star discrepancy can be thought as a special case 
$(p = \infty)$ of $L_p$-star discrepancy which is defined by
\begin{equation}
D^*_p(P) = ||F_{\text{unif}} - F_N||_p, 1 \le p \le \infty.
\end{equation}
Another discrepancy, which will be used in this thesis,
given by \cite{Hic96b}, called a generalized $L_2$-discrepancy, 
is 
\begin{eqnarray} \label{dis}
  && \hspace{-0.7cm}D^2(P) =  \nonumber \\
  && \hspace{-0.4cm}  -1+\frac{1}{N^2}\sum_{i,j=0}^{N-1}
  \prod_{r=1}^{s}  
   \left[ -\frac{(-\gamma^2)^\alpha}{(2\alpha)!} B_{2\alpha}(\{y_{ir}-y_{jr}\})+
  \sum_{k=0}^{\alpha}\frac{\gamma^{2k}}{(k!)^2}B_i(y_{ir})B_i(y_{jr}) \right],
\end{eqnarray}
where the notation $\{y\}$ means the fractional part of a number. 
The positive integer $\alpha$ indicates the degree of smoothness of 
integrands in the underlying Hilbert space, and the parameter $\gamma$
measures the relative importance given to the uniformity of the points
in low dimensional projections of the unit cube versus high dimensional
projections of the unit cube. The reproducing kernel leading to this 
discrepancy is for a Hilbert space
of integrands whose partial derivatives up to order $\alpha$ in each
coordinate direction are square integrable \cite{HicHon98a}.  It is
known that the root mean square discrepancy for scrambled
$(t,m,s)$-nets decays as $O(N^{-3/2+\epsilon})$ as $N \to \infty$
\cite{HicHon98a,HicYue00} for $\alpha \ge 2$. 
$B_i(y)$ is the $i$th Bernoulli
polynomial (see \cite{AbrSte64}). The first few Bernoulli 
polynomials, which are used in this thesis, are $B_0(y)=1$ and
\begin{eqnarray*}
  B_1(y)&=&y-\frac{1}{2},     \\
  B_2(y)&=&y^2-y+\frac{1}{6},  \\
  B_3(y)&=&y^3-\frac{3}{2}y^2+\frac{1}{2}y, \\
  B_4(y)&=&y^4-2y^3+y^2-\frac{1}{30}.
\end{eqnarray*}
See \cite{Hic96b} for further discussion about
the $L_p$-star discrepancy and several other discrepancies.

The following two sections introduce two important families of
quasi-Monte Carlo methods which are : 
\begin{itemize}
\item integration lattices \cite{Nie92,Slo85}, and
\item digital nets and sequences \cite{Nie92}.
\end{itemize} 
\end{section}

\begin{section}{Integration Lattices}
Rank-1 lattices, also known as good lattice point (glp) sets, 
were introduced by Korobov\cite{Kor59} and have been widely
studied since then. See \cite{SloJoe94,Nie92} for further details.
The formula for the node points of a rank-1 lattice is 
simply
\begin{equation}
P = \{\{ i\bfh / N \}: i = 0, \cdots, N-1 \},
\end{equation}
where $N$ is the number of points, $\bfh$ is an $s$-dimensional generating 
vector of integers
The formula for a lattice is rather simple. However finding
a good generating vector $\bfh$ that makes the lattice have low discrepancy
for the chosen $N$ and $s$ is not trivial.
Recently the formula for lattice is extended to an infinite sequence
\cite{HicEtal00}. This is done by using the radical inverse function, 
$Q_b(i)$. For any integer $b \ge 2$, let any non-negative integer $i$ 
be represented in base $b$ as $i = \cdots i_3 i_2 i_1$, where the digits 
$i_k$ take on values between $0$ and $b-1$. The function $Q_b(i)$ flips 
the digits about the decimal point, i.e.,
\begin{equation}
Q_b(i) = 0.i_1 i_2 i_3 \cdots (\mbox{ base } b) = \sum_{k=1}^{\infty}i_kb^{-k}.
\end{equation}
The sequence $\{ Q_b(i) : i = 0,1,\cdots \}$ is called the van der Corput 
sequence. An infinite sequence of imbedded lattices is defined by replacing
$i / N$ by $Q_b(i)$, i.e.,
\begin{equation}
y_i = \{ \{Q_b(i)\bfh \}: i = 0,1,\cdots \},
\end{equation} 
where the first $N$ points of this sequence are a lattice whenever 
$N$ is a power of $b$ \cite{HicEtal00}.

\begin{section}{Digital nets and Sequences}

Digital nets and sequences \cite{Lar97,Nie92} are one widely used 
types of low discrepancy point sets. 

Let $b \ge 2$ denote a prime number.
For any non-negative integer
$i=\cdots i_{3}i_{2}i_{1} (\mbox{base }b)$, we define the 
$\infty \times 1$
vector $\bfpsi(i)$ as the vector of its digits, i.e.,
$\bfpsi(i)=(i_{1}, i_{2}, \ldots)^{T}$.  For any point $z =
0.z_{1}z_{2} \cdots (\mbox{base }b) \in [0,1)$, let $\bfphi(z)=(z_{1},
z_{2}, \ldots )^{T}$ denote the $\infty \times 1$ vector of the digits
of $z$.  Let ${\bfC}_1,\ldots,{\bfC}_s$ denote predetermined $\infty
\times \infty$ {\em generator matrices}.  The {\em digital sequence}
in base $b$ is $\{\bfy_{0}, \bfy_{1}, \bfy_{2}, \ldots \}$, where each
$\bfy_{i}=(y_{i1},\ldots,y_{is})^{T}\in \Cs$ is defined by
\begin{equation} \label{digdef}
\bfphi(y_{ij}) = \bfC_{j}\bfpsi(i), \quad j=1, \ldots,s, \ i=0,1, \ldots.
\end{equation}
Here all arithmetic operations take place in mod $b$.
Thus, the right side of \eqref{digdef} should not give a vector ending
in an infinite trail of $b-1$s.

Digital nets and sequences are the special cases of $(t,m,s)$-nets and 
$(t,s)$-sequences and further explanation is given in Chapter 4.

\end{section}

\begin{section}{\texorpdfstring{$(t,m,s)$-Nets and $(t,s)$-Sequences}{(t,m,s)-Nets and (t,s)-Sequences}}
The concept of $t$ in $(t,m,s)$-nets and $(t,s)$-sequences provides
one way to measure the quality of low discrepancy sequences.
See \cite {Nie92} for more detail explanations. 
Here we provide a brief outline of elementary
intervals of $(t,m,s)$-nets and $(t,s)$-sequences.

Let $s \ge 1$ and $b \ge 2$ be integers. An elementary interval in base 
$b$ is a subinterval of  $[0,1)^s$ of the form
$$
\prod_{r=1}^{s} \left [ \frac{a_r}{b^{k_r}},\frac{a_r+1}{b^{k_r}} \right )
$$ 
with integer $k_r \ge 0$, $0 \le a_r < b^{k_r}$. Such an elementary interval
has volume $b^{-(k_1+\cdots+k_s)}$.

Let $m \ge 0$ be an integer. A finite set of $N= b^m$ points from $\Cs$ 
is a $(t,m,s)$-net in base $b$ if every elementary interval in base $b$ with 
the volume $b^{t-m}$ contains exactly $b^{t}$ points of a set, 
where $t$ is a nonnegative integer. Smaller $t$ values imply a better equidistribution
property for the net. Obviously the best case is when $t=0$, that is
every $b$-ary box of volume $1 / N$ contains exactly one of the $N$ points 
of the set.
 
For a given $t \ge 0$, an infinite sequence of points from  $\Cs$ is a 
$(t,s)$ sequence in base $b$ if for all integers $k \ge 0$ and $m \ge t$
the finite sequence $\{ y_{kb^m+1},\cdots,y_{(k+1)b^m} \}$ is a 
$(t,m,s)$-net in base $b$.

The constructions of $(t,m,s)$-nets and $(t,s)$-sequences in base $b = 2$ 
were introduced by Sobol' \cite{Sob67}. Later Faure \cite{Fau82} 
provided $(0,m,s)$-nets and $(0,s)$-sequences for prime bases $ b \ge s$. 
The general definitions are given by Niederreiter \cite{Nie87}.
Mullen, Mahalanabis and Niederreiter \cite{Nie95} provide tables
of attainable $t$-values for nets. 

The advantage of using nets taken from $(t,s)$-sequences 
is that one can increase $N$ through a sequence of values 
$N = \lambda b^m$, $1 \le \lambda < b$, and all the points
used in $\hat{I}_{\lambda b^m}$ are also used in $\hat{I}_{(\lambda+1)b^m}$.
Owen \cite{Owe96} introduces the $(\lambda,t,m,s)$-nets to 
describe such sequences.

The initial $\lambda b^m$ points of a $(t,s)$-sequence are well
distributed. For integers $m,t,b,\lambda$ with $m \ge 0$,
$ 0 \ge t \le m$, and $1 \le \lambda < b$, a sequence ${\bfy_i}$ of 
$\lambda b^m$ points is called a $(\lambda,t,m,s)$-net in base $b$ 
if every $b$-ary box of volume $b^{t-m}$
contains $\lambda b^t$ points of the sequence and no $b$-ary box of volume
$b^{t-m-1}$ contains more than $b^t$ points of the sequence.

For functions of bounded variation in the sense of Hardy and Krause,
the numerical integration by averaging over the points of a $(t,m,s)$-net
has an error of order $N^{-1}(\log N)^{s-1}$. Niederreiter \cite{Nie92}
discusses more precise error bounds for $(t,m,s)$-nets and 
$(t,s)$-sequences.  
\end{section}

\begin{section}{Randomized Quasi-Monte Carlo Methods}
Quasi-Monte Carlo methods can obtain a better convergence rate than
Monte Carlo methods, because the points are chosen to be more uniform.
However deterministic quasi-Monte Carlo methods have disadvantages 
compared with Monte Carlo methods. 
First, quasi-Monte Carlo methods are statistically
biased since the points are chosen in certain deterministic way.
Therefore the mean of the quadrature estimate is not the desired 
integral. Second, quasi-Monte Carlo methods do not facilitate simple error
estimates while Monte Carlo methods provide straightforward probabilistic 
error estimates. 

Randomizing quasi-Monte Carlo point is one way to 
overcome such disadvantages while still preserving their higher convergence 
rates. One may think of randomized quasi-Monte Carlo as a sophisticated 
variance reduction technique. There are two types of randomization.
One is adding the same $s$-dimensional random shift to every point.
Let $P$ be the original low discrepancy point set, 
then $P_{sh} = \{ i\bfh / N + \Delta \}: i = 0, \cdots, N-1 \}$, where 
$\Delta$ is a random vector uniformly distributed on $\Cs$. 
This type of randomization, which was 
introduced in \cite{CraPat76}, is often used with lattice rule because shifted
lattice rules retain their lattice structure. However, shifted nets do not 
necessarily remain as nets. Another type of randomization, which was proposed
by Art Owen \cite{Owe95}, is a more sophisticated method that randomly
permutes the digits of each point. And it often applies to nets because 
scrambled nets retain their nets properties. However, scrambled lattice 
rules do not necessary remain lattice rules. The more detail
discussions on Owen's scrambling can be found in Chapter 2.
\begin{figure}[h]
\begin{center}
\includegraphics[height=18cm,width=7.4cm]{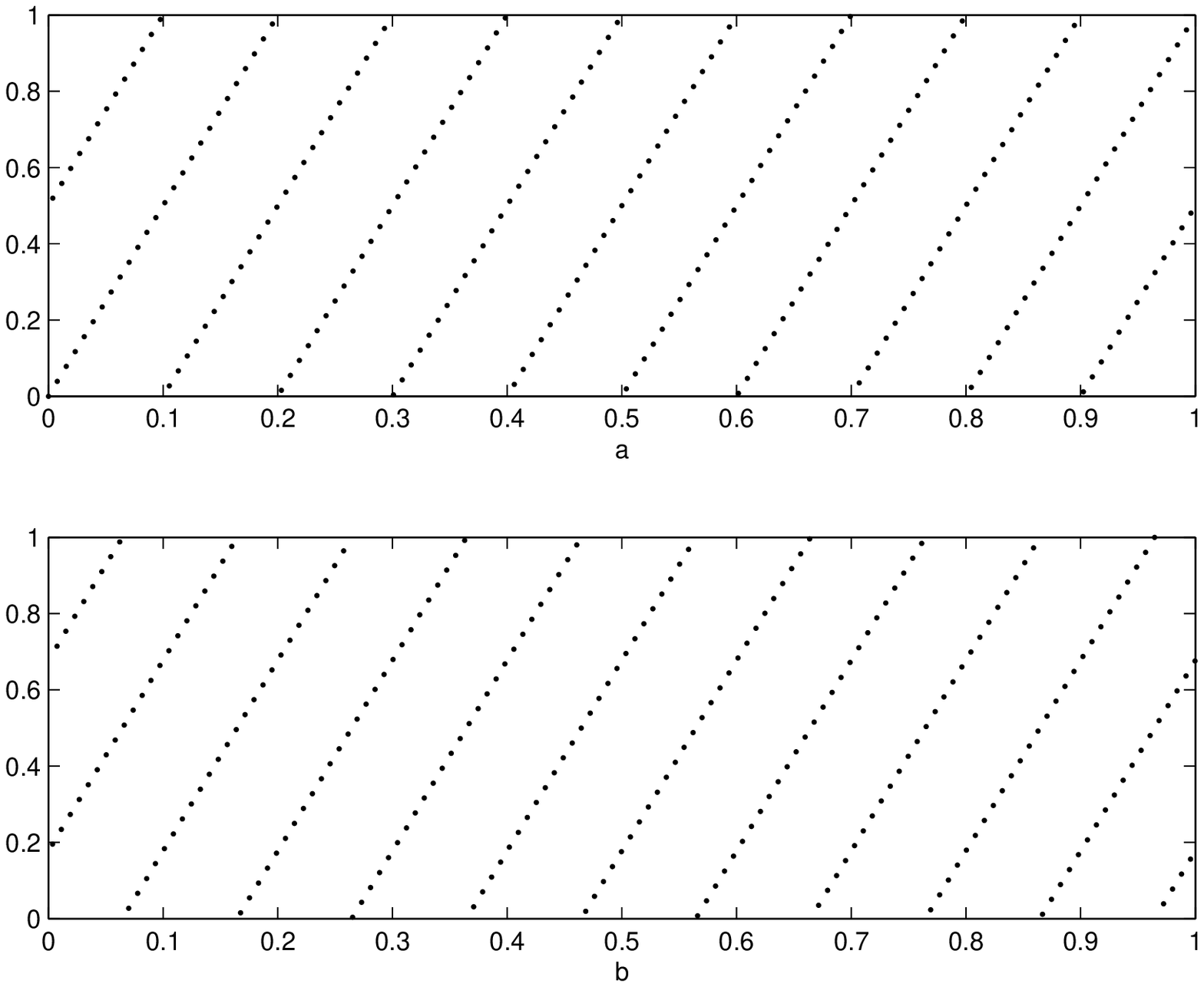}
\end{center}
\caption{\label{glp89} Dimensions 7 and 8 of the a) glp and
b) randomly shifted glp for $N= 256$}
\end{figure}
\end{section}
\begin{figure}[h]
\begin{center}
\includegraphics[height=18cm,width=7.4cm]{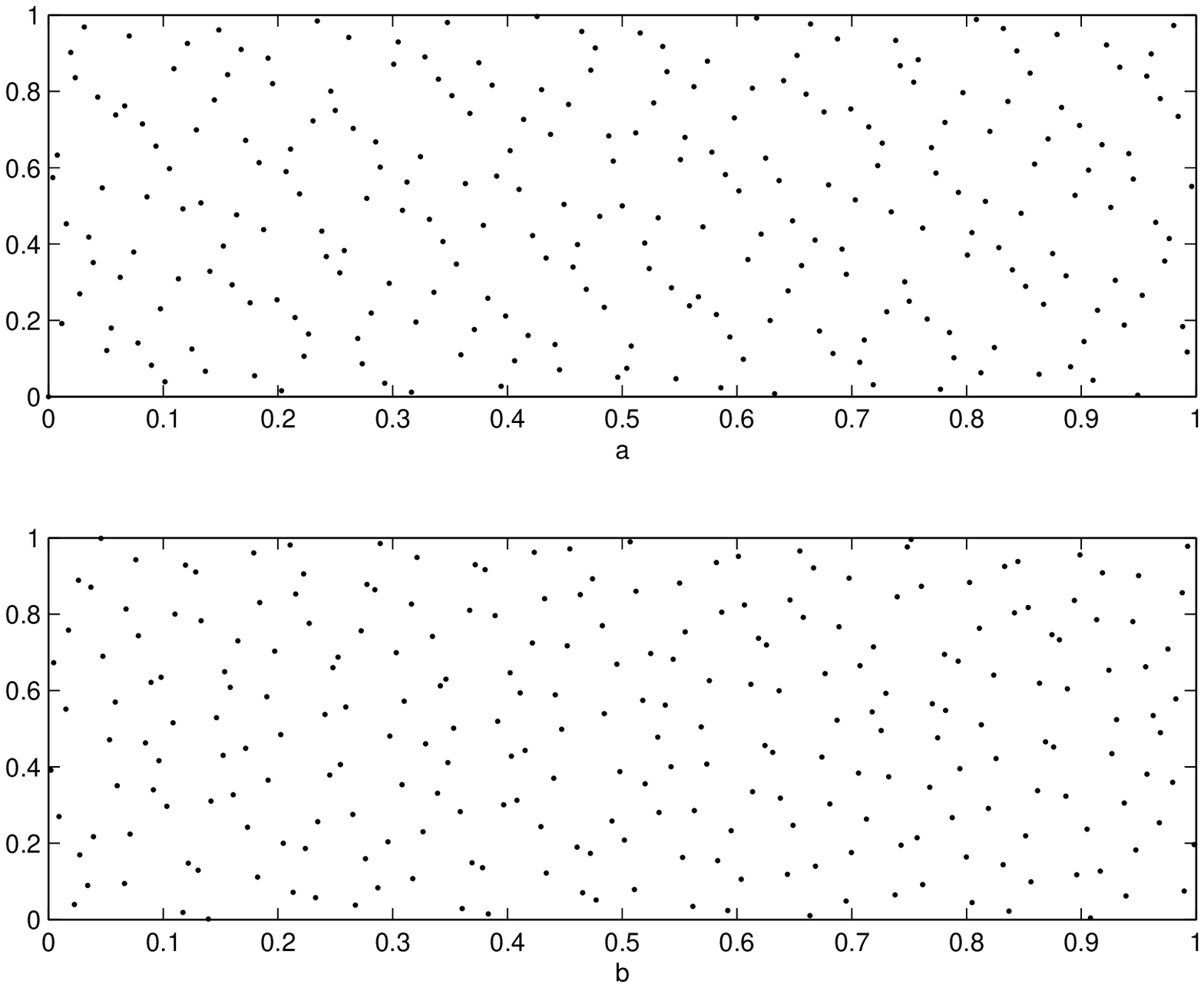}
\end{center}
\caption{\label{ss25} Dimensions 2 and 5 of the a) Sobol' and 
b) Scrambled Sobol' sequences for $N= 256$}
\end{figure}
\end{section}

\begin{section}{Outline of the Thesis}
In this chapter we have introduced the necessary background of
quasi-Monte Carlo methods. 

Chapter 2 discusses the randomization of the quasi-Monte Carlo methods 
specifically randomizing digital nets. Detailed descriptions
of Owen's randomization (or scrambling) and Faure-Tezuka's 
randomization (or tumbling) are provided. 
The actual implementation of randomization for digital $(t,s)$-sequences
and some numerical results are presented.

Chapter 3 explores the distribution of the discrepancy of scrambled 
digital $(t,m,s)$-nets which is based on recent theory. 
We fit the empirical distribution by a sum of chi squared random variables.
The distribution of the discrepancy of randomized lattice points is 
presented also.

Finding better digital nets is the main content of Chapter 4. Here we use 
optimization methods to find the generator matrices for digital 
$(t,m,s)$-nets. Evolutionary computation is introduced as a tool for 
a finding better net.

Chapter 5 presents the applications of the sequences that 
we have generated. Various problems are explored and the performance
of the new sequences is compared with other existing low discrepancy 
sequences. 

Finally the thesis ends with some concluding remarks. 

\end{section}
\end{chapter}

\newpage

\begin{chapter}{Implementing Randomized Digital \texorpdfstring{$(t,s)$}{(t,s)}-Sequences}

This chapter introduces randomized digital $(t,s)$-sequences.
Section 2.1 explains the details of Owen-scrambling and Section 2.2
introduces Faure-Tezuka-tumbling. Section 2.3 mainly deals with our method 
of scrambling that is nearly as general as Owen's scrambling and the detailed  
explanation of the implementation.
Section 2.4 explains our efforts to improve the 
efficiency of the generation time of the scrambled and non-scrambled 
sequences and discusses the time 
and the complexity of our algorithm. Section 5 presents 
various numerical results on the discrepancy for the newly generated randomized 
$(t,s)$-sequences.

\begin{section}{Owen's Scrambling}

Owen \cite{Owe95,Owe96,Owe97,Owe98b,Owe99a} proposes 
scrambled $(t,s)$-sequences as a hybrid of Monte Carlo and
quasi-Monte Carlo methods. This clever randomizations of quasi-Monte 
Carlo methods can combine the best of both by yielding higher accuracy 
with practical error estimates.
The following is a detailed description of Owen's scrambling.

Let $\{\bfx_{0}, \bfx_{1}, \bfx_{2}, \ldots \}$ denote the randomly
scrambled sequence proposed by Owen.  Let $x_{ijk}$ denote the
$k^{\text{th}}$ digit of the $j^{\text{th}}$ component of $\bfx_{i}$,
and similarly for $y_{ijk}$.  Then
\begin{eqnarray*}
x_{ij1} &=& \pi_j(y_{ij1}), \\
x_{ij2} &=& \pi_{y_{ij1}}(y_{ij2}), \\
x_{ij3} &=& \pi_{y_{ij1},y_{ij2}}(y_{ij3}), \quad \ldots, \\
x_{ijk} &=& \pi_{y_{ij1},y_{ij2},\cdots,y_{ijk-1}}(y_{ijk}), \ldots
\end{eqnarray*}
where the $\pi_{a_1 a_2 \dots}$ are random permutations of the
elements in mod $b$ chosen uniformly and mutually independently. 
Owen \cite{Owe98b} provides a geometrical description to help visualize 
his scrambling as follows:
The rule for choosing $x_{ij}$ is like cutting the unit cube into $b$
equal parts along the $\bfx_j$ axis and then reassembling these parts in 
random order to reform the cube. The rule for choosing $x_{ij2}$ is like 
cutting the unit cube into $b^2$ equal parts along the $\bfx_j$ axis, 
taking  them as $b$ groups of $b$ consecutive parts, and reassembling 
the $b$ parts within each groups in random order. The rule for $x_{ijk}$ involves cutting
the cube into $b^k$ equal parts along the $\bfx_j$ axis, 
forming $b^{k-1}$ groups of $b$ equal parts, and reassembling the $b$ 
parts within each group in random order.

Owen \cite{Owe95, Owe97} states that a randomized 
net inherits certain equidistribution properties of $(t,m,s)$-nets
by proving the following propositions:

{\bf Proposition 1} \quad
         {\it  If $\{ \bfy_i\}$ is a $(\lambda,t,m,s)$-net in base $b$
               then $\{ \bfx_i\}$ is a $(\lambda,t,m,s)$-net in base $b$ 
               with probability 1.}

{\bf Proposition 2} \quad
         {\it Let $\bfy$ be a point in $[0,1)^s$ and let $\bfx$
              be the scrambled version of $\bfy$ as described above.  
              Then $\bfx$ has the uniform distribution on $[0,1)^s$. }
\end{section}

\begin{section}{Faure-Tezuka's Tumbling}
Faure and Tezuka \cite{FauTez01} proposed another type of randomizing 
digital nets and sequences. However the effect of the 
Faure-Tezuka-randomization can be thought as re-ordering 
the original sequence, rather than permuting its digits like 
the Owen-scrambling. Thus we refer to this method as 
Faure-Tezuka-tumbling as suggested by Art Owen (personal communication).
The following describes details of Faure-Tezuka-tumbling.

For any $m,\lambda=0,1, \ldots$, let
$i=\lambda b^{m}, \ldots, (\lambda+1)b^{m}-1$.  Then $\bfpsi(i)$
takes all possible values in its top $m$ rows.  By the same token
$\bfL^{T}\bfpsi(i) + \bfe$ takes on all possible values in its top $m$
rows, but not necessarily in the same order.  Therefore, the
Faure-Tezuka-tumbled $(t,m,s)$-net, $\{\bfz_{\lambda b^{m}}, \ldots,
\bfz_{(\lambda+1)b^{m}-1}\}$, obtained by replacing
$\bfpsi(i)$ of \eqref{digdef} by $\bfL^{T}\bfpsi(i) + \bfe$
is just the same set as the original $(t,m,s)$-net, $\{\bfy_{\nu b^{m}}, \ldots,
\bfy_{(\nu+1)b^{m}-1}\}$, given by \eqref{digdef} for some $\nu$.
\end{section}

\begin{section}{Implementation}
There is a cost in scrambling $(t,s)$-sequences.  The manipulation of
each digit of each point requires more time than for an unscrambled
sequence. In addition Owen's scrambling can be rather tedious because 
of the bookkeeping required to store the many permutations.  Here we present
an alternative method that is only slightly less general than Owen's proposal 
but is an efficient method of scrambling that minimizes this cost.  
This is called an Owen-scrambling to recognize that it is done in the
spirit of Owen's original proposal. 

\begin{subsection}{The Method of Scrambling} 

Let $\bfL_{1},\cdots,\bfL_{s}$, be nonsingular lower triangular
$\infty \times \infty$ matrices and let $\bfe_{1},\cdots,\bfe_{s}$ be
an $\infty \times 1$ vector.  Assume that for all $j$ any linear
combination of columns of $\bfL_{j}\bfC_{j}$ plus $\bfe_{j}$ does not
end in an infinite trail of $b-1$s.  A particular {\em
Owen-scrambling} $\{\bfx_{i}\}$ of a digital sequence $\{\bfy_{i}\}$
is defined as
\begin{equation}\label{Owescrdef}
    \bfphi(x_{ij}) = \bfL_{j}\bfphi(y_{ij}) + \bfe_{j} =
    \bfL_{j}\bfC_{j}\bfpsi(i) + \bfe_{j}, \quad j=1, \ldots,s, \
    i=0,1, \ldots.
\end{equation}
The left multiplication of each generator matrix $\bfC_{j}$ by a
nonsingular lower triangular matrix $\bfL_{j}$ yields a new generator
matrix that satisfies the same condition \eqref{tcond} as the
original.  The addition of an $\bfe_{j}$ just acts as a digital shift
and also does not alter the $t$-value.  Therefore the scrambled
sequence is also a digital $(t,s)$-sequence.

To get a randomly scrambled sequence one chooses the elements of the
$\bfL_{j}$ and $\bfe_{j}$ randomly.  The resulting randomized sequence
has properties listed in the theorem below.  Although these properties
are not equivalent to the Owen's original randomization, they are
sufficient for the analysis of scrambled net quadrature rules given in
\cite{Owe95,Hic96b,Owe96,Owe97,Owe98b,HicHon98a,Yue98,YueMao99,HicYue00}.

\begin{theorem} 
Let $\{\bfx_{i}\}$ be an Owen-scrambling of a digital $(t,s)$-sequence
where elements of $\bfL_{1},\cdots,\bfL_{s}, \bfe_{1},\cdots,\bfe_{s}$ 
are all chosen randomly and independently.  The diagonal elements of 
$\bfL_{1},\cdots,\bfL_{s}$ are chosen uniformly on $\{1, \ldots, 
b-1\}$, and the other elements are chosen uniformly on $\Zb$.  Then 
this randomly scrambled sequence satisfies the following properties 
almost surely:
\begin{enumerate}
    \setlength{\labelsep}{2ex}
    \renewcommand{\labelenumi}{\roman{enumi}.}
    
    \item $x_{ijk}$ is uniformly distributed on $\Zb$;
    
    \item if $y_{ijh}=y_{\ell jh}$ for $h < k$ and $y_{ijk} \not=
    y_{\ell jk}$, then
   \begin{enumerate}
    \renewcommand{\labelenumii}{\alph{enumii}.}
     \item $x_{ijh}=x_{\ell jh}$ for $h<k$
     \item $(x_{ijk},x_{\ell jk})$ are uniformly distributed on $\{(x,y) 
     \in \Zb^{2} : x\ne y\}$;
     \item $x_{ijh},x_{\ell jh}$ are independent for $h>k$  
   \end{enumerate}
   \item $(x_{ij},x_{\ell j})$ is independent from $(x_{ir},x_{\ell r})$ for
   $j \not= r$
\end{enumerate}
\end{theorem}
{\bf Proof} \quad
The qualification ``almost surely'' is required to rule
out the zero probability event that for some $j$ there exists a linear
combination of columns of $\bfL_{j}\bfC_{j}$ plus $\bfe_{j}$ that ends
in an infinite trail of $b-1$s.  Assertion iii.\ follows from the fact
that the elements of $\bfL_{j}, \bfe_{j},\bfL_{r}$, and $\bfe_{r}$ are
chosen independently of each other.  Now the other two assertions are
proved.
    
Note from \eqref{Owescrdef} that $x_{ijk} =
\bfl^T_{jk}\bfphi(y_{ij})+e_{jk}$, where $\bfl^T_{jk}$ is the
$k^{\text{th}}$ row of $\bfL_{j}$, and all arithmetic operations are
assumed to take place in the finite field.  Since $e_{jk}$ is
distributed uniformly on $\Zb$, it follows that $x_{ijk}$ is also
distributed uniformly on $\Zb$.

To prove ii.\ consider $\bfphi(x_{ij})-\bfphi(x_{\ell j})=
\bfL_{j}[\bfphi(y_{ij})-\bfphi(y_{\ell j})]$.  Under the assumption of
ii.\ the first $k-1$ rows of $[\bfphi(y_{ij})-\bfphi(y_{\ell j})]$ are
zero, which implies that the first $k-1$ rows of
$\bfphi(x_{ij})-\bfphi(x_{\ell j})$ are also zero, since $\bfL_{j}$ is
lower triangular.  This immediately implies iia.  Furthermore,
$x_{ijk}-x_{\ell j k}= \bfl^T_{jk}[\bfphi(y_{ij})-\bfphi(y_{\ell j})]
= l_{jkk}(y_{ijk}-y_{\ell jk})$, where $l_{jkk}$ is the
$k^{\text{th}}$ diagonal element of $\bfL_{j}$.  Since
$y_{ijk}-y_{\ell jk} \ne 0$, and $l_{jkk}$ is a random nonzero element
of $\Zb$, it follows that $x_{ijk}-x_{\ell j k}$ is a random nonzero
element of $\Zb$.  Combined with i., this implies iib.  For $h>k$ it
follows that $x_{ijh}-x_{\ell j h}=
\bfl^T_{jh}[\bfphi(y_{ij})-\bfphi(y_{\ell j})] = \cdots +
l_{jhk}(y_{ijk}-y_{\ell jk}) + \cdots$.  The fact that $y_{ijk}-y_{\ell 
jk} \ne 0$, and $l_{jhk}$ is uniformly distributed on $\Zb$, combined 
with i.\ now imply iic.

To understand why the randomization given by \eqref{Owescrdef} is not
as rich as that originally prescribed by Owen, consider randomizing
$y_{i11}$, the first digit of $y_{i1}$ for different values of $i$. 
The $y_{i11}$ take on $b$ different values, and there are $b!$
possible permutations of these values.  However, the formula given by
\eqref{Owescrdef} is $x_{i11}=l_{111}y_{i11}+e_{11}$, where $l_{jhk}$
denotes the $h,k$ element of $\bfL_{j}$, and $e_{jk}$ denotes the
$k^{\text{th}}$ element of $\bfe_{j}$.  Since there are only $b-1$
possible values for $l_{111}$ and only $b$ possible values for
$e_{11}$, this formula cannot give at most $b(b-1)$ permutations.

The Faure-Tezuka-tumbling \cite{FauTez01} randomizes the digits of
$i$ before multiplying by the generator matrices.  Let $\bfL$, be a
nonsingular lower triangular $\infty \times \infty$ matrix, and let
$\bfe$ be an $\infty \times 1$ vector with a finite number of nonzero
elements.  A particular {\em tumbling}, $\{\bfz_{i}\}$,
of a digital sequence is defined as
\begin{equation}\label{FTscrdef}
    \bfphi(z_{ij}) = \bfC_{j}[\bfL^{T}\bfpsi(i) + \bfe], \quad j=1,
    \ldots,s, \ i=0,1, \ldots.
\end{equation}
A particular {\em scrambling-tumbling}, $\{\bfx_{i}\}$,
of a digital sequence is defined as
\begin{equation}\label{OFTscrdef}
    \bfphi(x_{ij}) = \bfL_{j}\bfC_{j}[\bfL^{T}\bfpsi(i) + \bfe] +
    \bfe_{j}, \quad j=1, \ldots,s, \ i=0,1, \ldots.
\end{equation}
in order to obtain random Faure-Tezuka-tumbling the elements of $\bfL$ and $\bfe$ 
are chosen randomly.

The scrambling method described here applies only to {\em digital} 
$(t,s)$-sequences.  Since all known general constructions 
of $(t,m,s)$-nets and $(t,s)$-sequences are digital constructions 
\cite{Sob67,Fau82,Nie88,Lar98,Lar97,NieXin98a}, this restriction is not 
serious. The algorithm we implement here builds on the
algorithms of \cite{BraFox88,BraFoxNie92} and includes the recent digital
Niederreiter-Xing sequences by \cite{NieXin96,NieXin98a}.
\end{subsection}

\begin{subsection}{Generators for the Randomized Nets} 

Generators for scrambled Sobol' \cite{Sob67}, Faure \cite{Fau82}, and
Niederreiter \cite{Nie88} points have been programmed following the
Fortran codes of \cite{BraFox88,BraFoxNie92}.  Below the differences
and improvements that we have made are highlighted and explained.  The
code for the Niederreiter points has been extended to higher dimension
by Giray \"Okten, and we have employed this extension.  A generator
for Niederreiter-Xing points has been coded based on the generator
matrices provided by \cite{Pir02}.

The different scrambled or non-scrambled sequences given by
\eqref{digdef}, \eqref{Owescrdef}, and \eqref{FTscrdef} are just
special cases of \eqref{OFTscrdef}.  Thus, \eqref{OFTscrdef} is the
basic algorithm to be implemented, with different choices given to the
user as to how to choose the matrices $\bfL_{1}, \ldots, \bfL_{s},
\bfL$ and the vectors $\bfe_{1}, \ldots, \bfe_{s}, \bfe$.

In practice one cannot store all the entries of an $\infty \times
\infty$ matrix.  Assuming that $b^{m_{max}}$ is the maximum amount 
of points that is required in
one run, and that $K$ is the number of digits to be
scrambled.  Then, one can restrict $\bfL_{1}, \ldots, \bfL_{s}$ to
$K\times m_{max}$ matrices, $\bfC_{1}, \ldots, \bfC_{s}, \bfL$ to 
$m_{max} \times m_{max}$ matrices, $\bfe_{1}, \ldots, \bfe_{s}$ to $K\times 1$ vectors, and
$\bfe$ to a $m_{max} \times 1$ vector.  One might think that it is best to
scramble all available digits, but Section 2.4 provides
evidence that this is not necessary.

The value of $K_{\max}$, the maximum possible value of $K$, is 31,
which corresponds to the number of bits in the largest integer
available.  When $b=2$, the algorithms for generating digital
sequences may be implemented more efficiently by storing an array of
bits as an integer.  For the scrambled Sobol' generator the value of
$s$, the maximum dimension is $40$, the same as in
\cite{BraFox88}.  For the scrambled Faure sequence, $s=500$, and for
the scrambled Niederreiter generator $s=318$.  For the
Niederreiter-Xing points $S=16$ based on the generator matrices
available so far.  The program may be modified to allow for higher
dimensions as the generator matrices become available.

To save time some matrices and vectors are premultiplied so that
\eqref{OFTscrdef} becomes
\begin{equation}\label{OFTscrimp}
    \bfphi(x_{ij}) = \tilde{\bfC}_{j} \bfpsi(i) + \tilde{\bfe}_{j},
    \quad j=1, \ldots,s, \ i=0,1, \ldots.
\end{equation}
where $\tilde{\bfC}_{j}=\bfL_{j}\bfC_{j}\bfL^{T}$, and 
$\tilde{\bfe}_{j}=\bfL_{j}\bfC_{j}\bfe+\bfe_{j}$.

Another time saving feature is to use a gray code \cite{Lic98}. 
Recursively define the $\infty \times 1$ vector function $\tbfpsi(i)$
as follows:
\begin{eqnarray*}
    \tbfpsi(0)&=&(0,0,\ldots)^{T}\\
    \tbfpsi(i+1)&=&\tbfpsi(i) - \delta_{k, \hat k_{i+1}}, \quad \\
    \text{where } \hat k_{i}&=& \min\{k : \lfloor i b^{-k} \rfloor \ne i 
    b^{-k} \},
\end{eqnarray*}
$\delta_{kl}$ is the Kronecker delta function, and where again all 
arithmetic operations take place in the finite field.  For example, 
if $b=3$, then
\begin{eqnarray*}
\tbfpsi(1) &=& (2,0,0\ldots)^{T}, \\
\tbfpsi(2) &=& (1,0,0\ldots)^{T}, \\
\quad \tbfpsi(3) &=& (1,2,0\ldots)^{T}, \\ \vdots
\end{eqnarray*}
Note that only one digit of $\tbfpsi$ changes as the argument is
increased by one.  Replacing $\bfpsi(i)$ by $\tbfpsi(i)$ in
\eqref{OFTscrimp} still results in a scrambled digital sequence; it 
just has the effect of re-ordering the points.  The efficiency 
advantage is that the digits of the $i+1^{\text{st}}$ point can be 
obtained from those of the  $i^{\text{th}}$ point by the iteration:
\begin{equation} \label{uppoint}
\bfphi(x_{i+1,j}) = \bfphi(x_{ij}) - \tilde{\bfC}_{j \hat k_{i+1}},
\end{equation}
where $\tilde{\bfC}_{j l}$ is the $l^{\text{th}}$ column of 
$\tilde{\bfC}_{j}$.

The implementations of the scrambled Sobol' and Niederreiter sequences
closely followed the structures of the original program and our
changes were rather minor.  However, the structure of the Faure
sequence generator has been altered more substantially to improve
efficiency.  The scrambled Sobol', Niederreiter, and Niederreiter-Xing
generators all have base $b=2$, which allows the use of logical
operations instead of additions and multiplications.

\begin{figure}[h]
\begin{center}
\includegraphics[height=20cm,width=5.2cm]{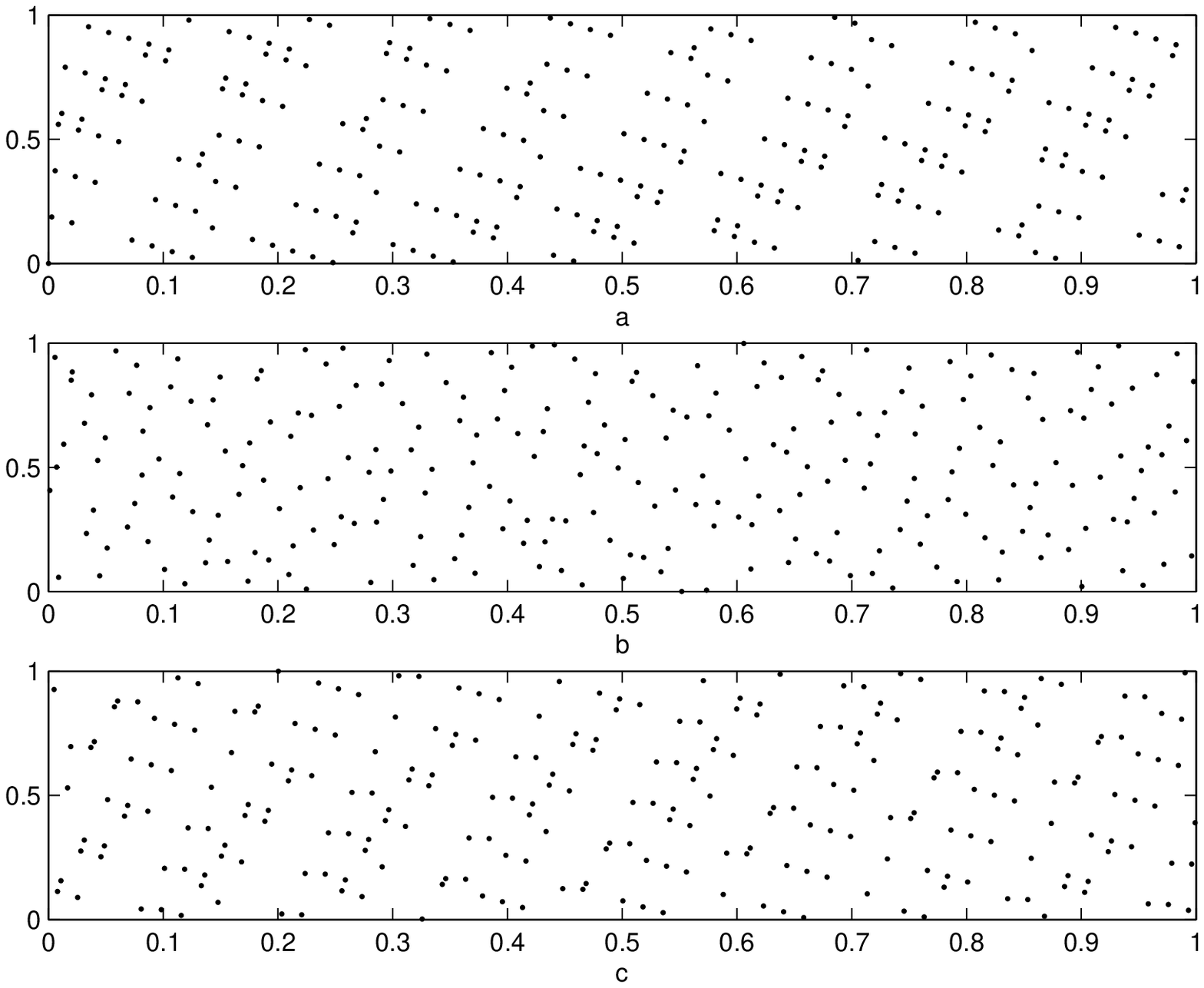}
\end{center}
\caption{\label{nx} Dimensions 6 and 7 of the a) Faure, 
b) Scrambled Faure, and c)Tumbled Faure sequences for $N= 256$}
\end{figure}
\end{subsection}

\begin{subsection}{Time and Space complexity}

This section considers the number of operations and the amount of memory
required to generate the first $N$ points of a scrambled
$(t,s$)-sequence with $b^{m-1} < N \le b^{m}$.  The space required to
store the scrambled generator matrices is $\Order(smK)$, where $K \ge
m$ is the number of digits scrambled.  The time required to generate
and pre-multiply the matrices is $\Order(sm^{2}K)$.  To generate one
new point from the scrambled sequence requires just $\Order(sK)$
operations according to \eqref{uppoint}, so the total number of
operations to generate all points is $\Order(sNK)$.  Thus, the
preparations steps take proportionally negligible time as $m \to
\infty$.  Note that for unscrambled points the corresponding time and
space complexities can be obtained by replacing $K$ by $m$.

The following table shows time comparison of sequence generating time
of unscrambled and scrambled digital point sets using $K=m+15$.  The
times have been normalized so that the time to generate the Sobol' set
equals 1.
\begin{table}
    \caption{Time required to compute low discrepancy points}
\begin{center}
    \begin{tabular}{rccccccc} 
  Algorithm & SOBOL & SSOBOL & NIED & SNIED & Old FAURE & FAURE & 
  SFAURE  \\\hline
  Time & 1 & 2 & 2 & 2 & 10 & 3 & 6 
\end{tabular}
\end{center}
\end{table}
From the table it shows that generating time for scrambled Sobol',
Niederreiter and Faure sequences are about 2 times more than that of
the original sequences based on the same number of extra digit
scrambling.  This is due the fact that one must manipulate $K$ digits
of each number, not just $m$.  The new Faure generator takes about
$1/3$ the time of the original one due to the use of the gray code.
\end{subsection}
\end{section}

\begin{section}{Numerical Results}
The discrepancy used in here is the squared discrepancy \eqref{dis}
and it has been scaled by a dimension-dependent constant so that 
the root mean square scaled discrepancy of a simple random sample 
is $N^{-1/2}$.

Figures \ref{ss2sd} shows the histogram of the distribution of the 
square discrepancy of the particular Owen-scrambled Sobol' sequence 
for the different numbers of scrambled digits, $K$. $K$ has been chosen to 
be $10$, $20$, and $30$ 
for $s = 2$ and $m = 10$ with 200 independent replications.
The range of the $x$-axis is are $10^{-15},\ldots,10^{-6}$ in this figure. 
From the figure, there is only little distinction in the distribution 
of the square discrepancy for $K= 20$ and $30$. However a larger 
number of scrambled digits produces a wider range of discrepancies. 
From Figure \ref{ssd2sd} the root mean square discrepancies 
for $K = 20$ and $30$ are nearly the 
same. However the root mean square discrepancy for $K = 10$
looses the benefit from scrambling as $m$ increases. 

Figures \ref{ssd2sd} and \ref{ssd5sd} plot the discrepancy
of the particular Owen-scrambled Sobol' sequence for different choices of
$K$ and the choices are same as before. 
The choices of dimension are $s=2$ and $5$. These figures show 
the root mean square discrepancy of 100 different replications.
There is almost no difference in the root mean square 
discrepancy between scrambling $K = 20$ or $30$ digits. 
However for $K=10$ the discrepancy looses its superior performance 
as $N$ increases, which indicates an insufficient number of scrambled 
digits. Also notice that choosing the number of scrambled digits is 
independent from $s$.

\begin{figure}[ht]
\begin{center}
\includegraphics[height=12cm]{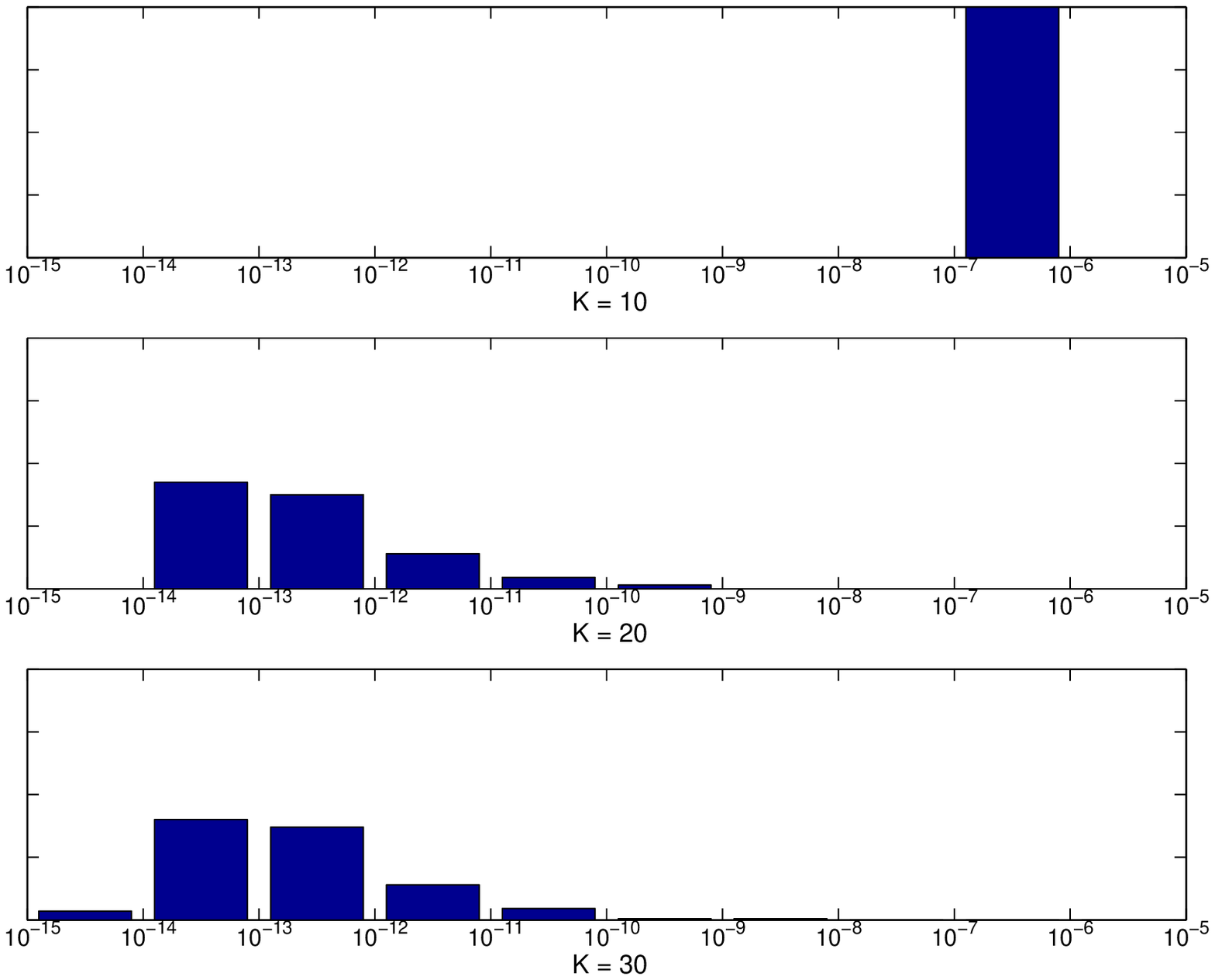}
\end{center}
\caption{\label{ss2sd} The histogram of the square
discrepancy for the Owen-scrambled Sobol' net for $K = 10$,
$20$, and $30$ for $s = 2$ and $N = 2^{10}$.}
\end{figure}

\begin{figure}[ht]
\begin{center}
\includegraphics[height=8cm]{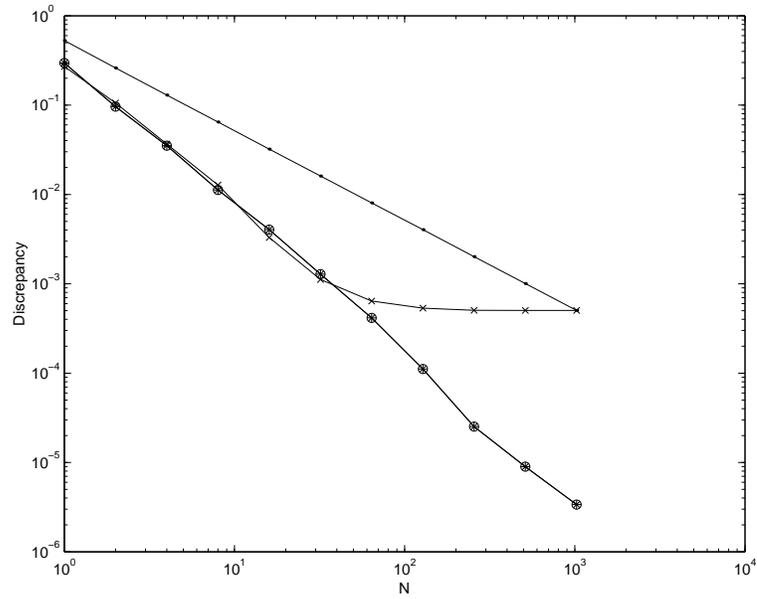}
\end{center}
\caption{\label{ssd2sd} The scaled root mean squared discrepancy with
$\alpha = 2$ as a function of $N$ for the Owen-scrambled Sobol' net for
$K = 10$ (x), $20$ (o), $30$ (*), and unscrambled (.) for 
$s = 2$.}
\end{figure}

\begin{figure}[ht]
\begin{center}
\includegraphics[height=8cm]{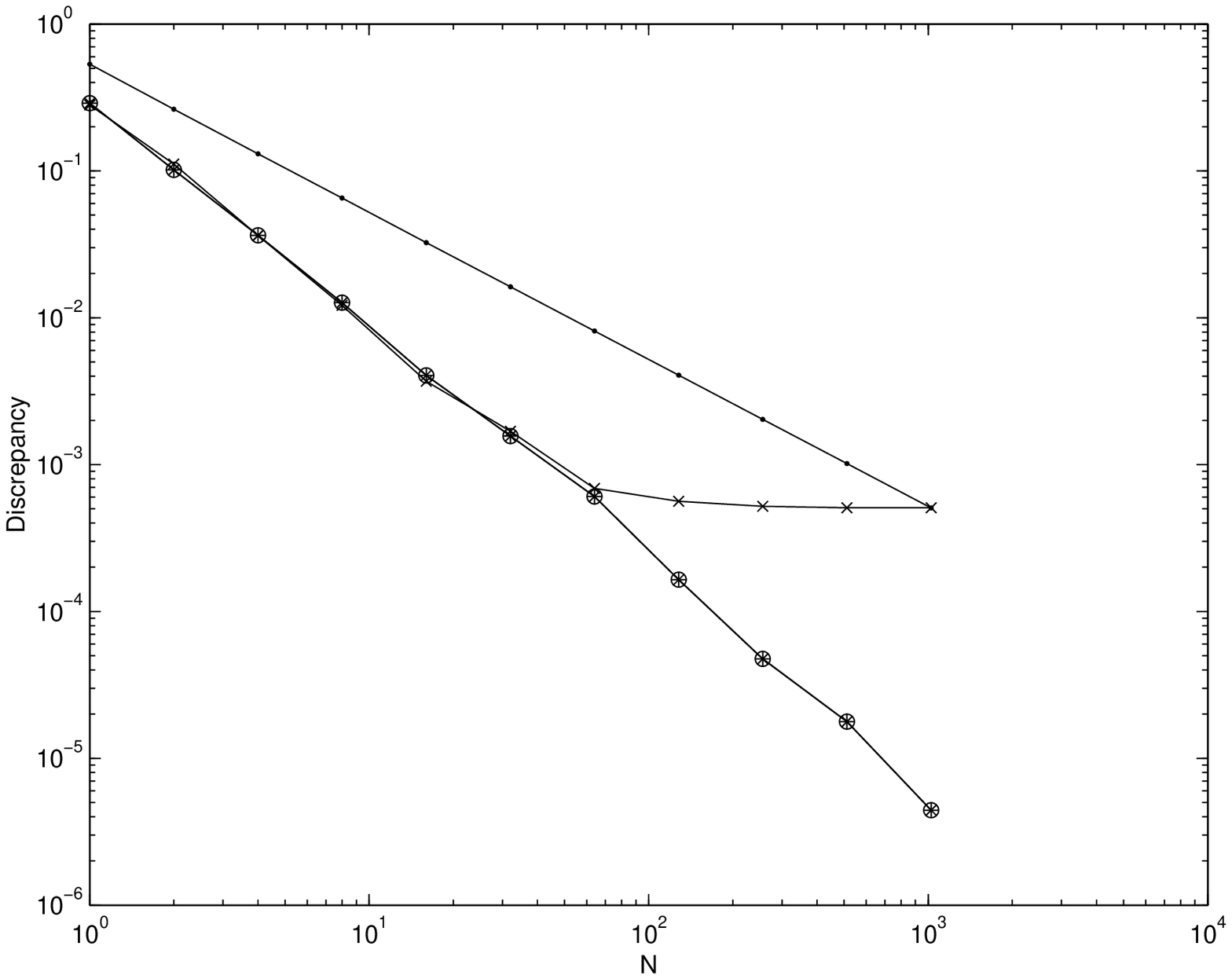}
\end{center}
\caption{\label{ssd5sd} Same as Figure \ref{ssd2sd} for $s = 5$.}
\end{figure}

Figure \ref{stf5} plots the discrepancies of the Faure sequence, 
the scrambled Faure sequence, and the tumbled Faure sequence.
The number of scrambled digits, $K$, is chosen to be 20, and 100 different
replications have been performed. The choice of dimension is $s = 5$.
Since the Faure sequence takes a base $b$ as a smallest prime number which 
is equal to or greater than $s$, $b$ assigned to be $5$.
Therefore $N$ is chosen to be $\lambda 5^{m}$, where 
$\lambda = 1,\cdots,b-1$ and $m = 1,\cdots,4$.
From Figure \ref{stf5} scrambled sequences perform the 
best among all sequences and the tumbled sequence performs better than the 
original Faure sequence. Also notice that the scrambling and tumbling 
procedure help to flatten out the humps, which are present in the original 
sequence when the value of $\lambda \neq 1$.

\begin{figure}[ht]
\begin{center}
\includegraphics[height=8cm]{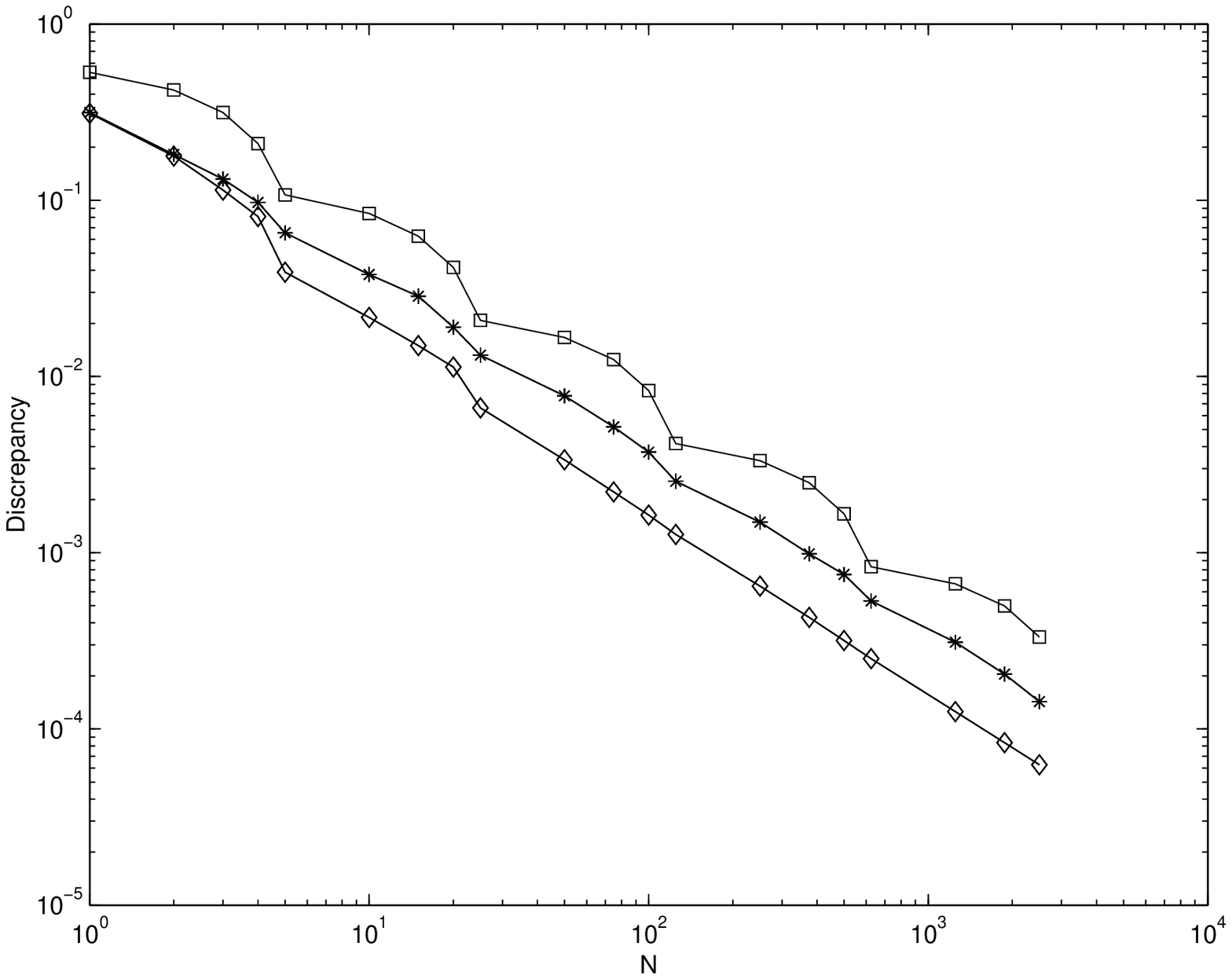}
\end{center}
\caption{\label{stf5} The scaled root mean squared discrepancy with
$\alpha = 2$ as a function of $N$ for the Faure ($\square$),
the Owen-scrambled Faure ($\diamond$), and 
the Faure-Tezuka-Tumbled Faure (*) nets for $s = 5$.}
\end{figure}
  
Figures \ref{gama1s}--\ref{gama2} show
the root mean square discrepancies of randomly scrambled
$(t,m,s)$-nets in base 2 with $N=2^{m}$ points that have been calculated
by using 100 different replications.  The choices of dimension are
$s=1$ and $10$.  The number of scrambled digits, $K$, is chosen
to be 31.

\begin{figure}[ht]
\begin{center}
\includegraphics[height=8cm]{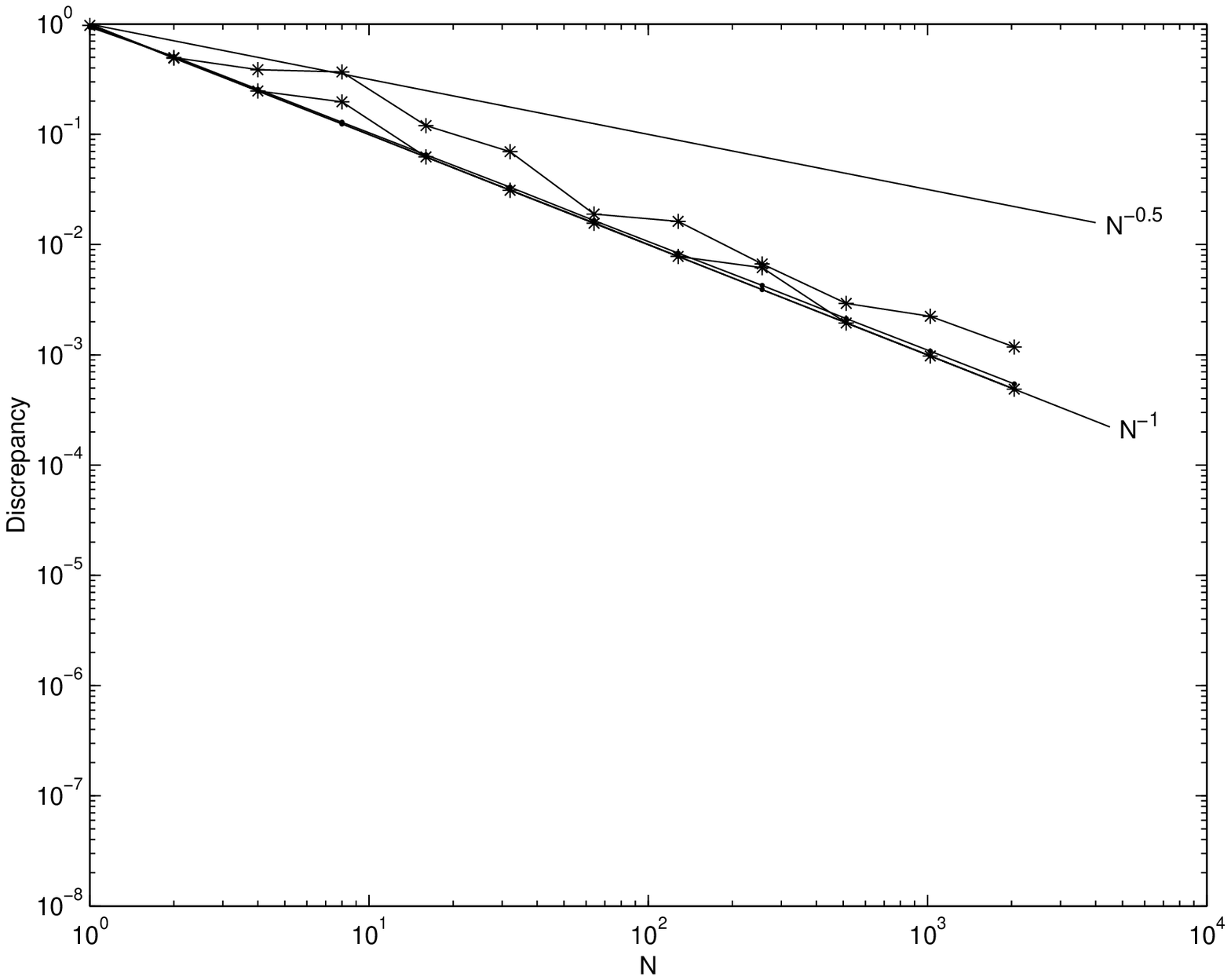}
\end{center}
\caption{\label{gama1s} The scaled root mean square discrepancy with
$\alpha = 1$ as a function of $N$ for the Owen-scrambled Sobol' ($\cdot$) 
and the Niederreiter-Xing ($*$) nets for $s=1$ and $10$.
\label{snetdisa1}}
\end{figure}
\begin{figure}[ht]
\begin{center}
\includegraphics[height=8cm]{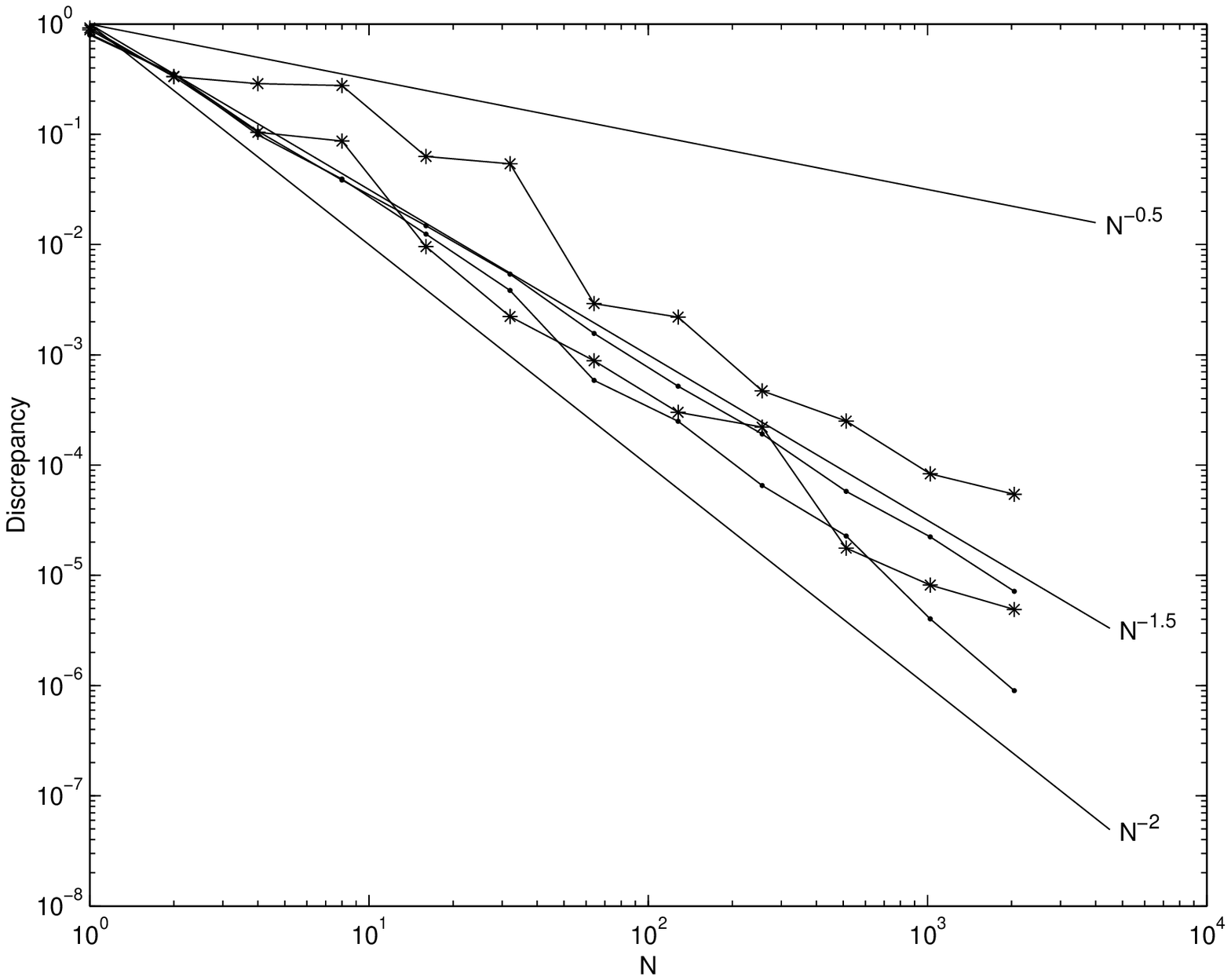}
\end{center}
\caption{\label{gama2s} The scaled root mean square discrepancy with
$\alpha = 2$ as a function of $N$ for Owen-scrambled Sobol' ($\cdot$) and
Niederreiter-Xing ($*$) nets for $s=1$ and $10$.  \label{snetdisa2}}
\end{figure}
\begin{figure}[ht]
\begin{center}
\includegraphics[height=8cm]{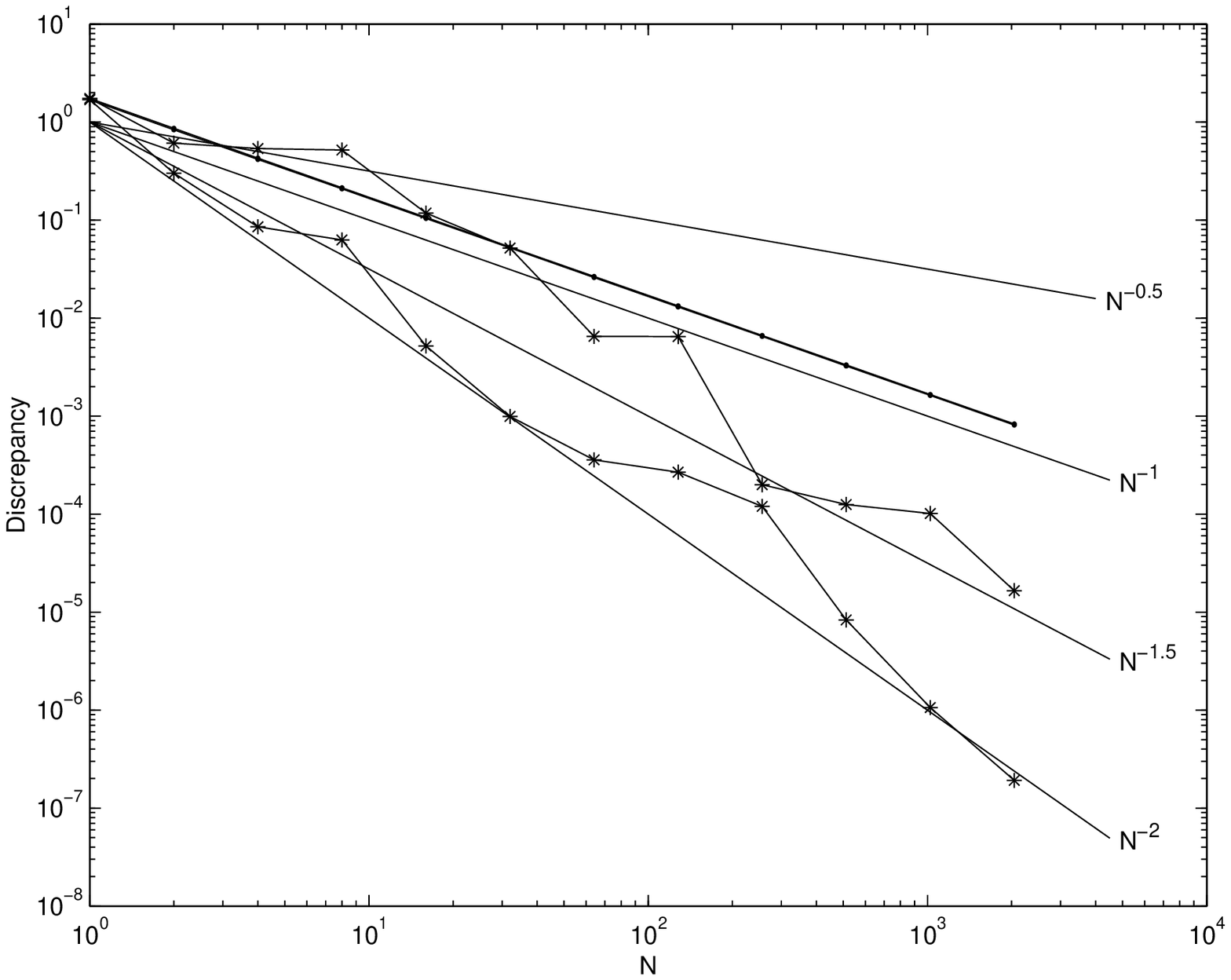}
\end{center}
\caption{\label{gama2} The scaled the root mean square discrepancy with
$\alpha = 2$ as a function of $N$ for non-scrambled Sobol' ($\cdot$) and
Niederreiter-Xing ($*$) nets for $s=1$ and $10$. \label{netdisa2}}
\end{figure}

Theory predicts that the root mean square discrepancy of nets should
decay as $\Order(N^{-1})$ for $\alpha=1$, and as $\Order(N^{-3/2})$ for 
$\alpha=2$.  Figures \ref{gama1s}--\ref{gama2} display this behavior, 
although for larger $s$ the asymptotic decay rate is only attained for 
large enough $N$.  However, even for smaller $N$ the discrepancy decays no 
worse than $\Order(N^{-1/2})$, the decay rate for a simple Monte 
Carlo sample.

The discrepancy of unscrambled nets does not attain the
$\Order(N^{-3/2})$ decay rate for $\alpha=2$.  Although for these 
experiments all possible digits were scrambled, our experiment 
suggests that scrambling about $K=m+10$ digits is enough to obtain 
the benefits of scrambling.  Scrambling more than this number does 
not improve the discrepancy but increases the time required to 
generate the scrambled points.
\end{section}
\end{chapter}

\begin{chapter}{The Distribution of the Discrepancy of Scrambled
Digital \texorpdfstring{$(t,m,s)$}{(t,m,s)}-Nets}

Recently, Loh \cite{Loh01} proved that the scrambled net quadrature obeys
a central limit theorem.  This chapter studies the empirical
distribution of the square discrepancy of scrambled nets and compares
it to what one would expect from Loh's central limit theorem.  
Here we use scrambling techniques mentioned in Chapter 2.
Organization of the chapter is as follows.  Section 3.1 provides a 
general description of the discrepancy with respect to the error of 
the quadrature and derives the theoretical asymptotic
distribution of the square discrepancy of scrambled nets.  
Section 3.2 explains the procedure for fitting the empirical distribution of 
the square discrepancy of scrambled nets to the theoretical asymptotic
distribution.  Section 3.3 discusses the experimental results.  

\begin{section}{The Distribution of the Squared Discrepancy
of Scrambled Nets}

The discrepancy measures the uniformity of the distribution of a set
of points and it can be interpreted as the maximum possible quadrature
error over a unit ball of integrands.  Since the discrepancy depends
only on the quadrature rule, it is often used to compare
different quadrature rules.

Let $K(x,y)$ be the
reproducing kernel for some Hilbert space of integrands whose domain is
$\Cs$.  To approximate the integral 
$$
\int_{\Cs} f(x) \ dx
$$
one may use the quasi-Monte Carlo quadrature rule 
$$
\frac{1}{N} \sum_{z \in P} f(z),
$$
where $P \subset \Cs$ is some well-chosen set of $N$ points.  The
error of this rule is 
$$
\Err(f)= \int_{\Cs} f(x) \ dx - \frac{1}{N} \sum_{z
\in P} f(z),
$$
and the square discrepancy is 
$$
D^2(P) = \Err_x (\Err_y(K(x,y))),
$$
where $\Err_{x}$ implies that the error functional is
applied to the argument $x$ \cite{Hic99b}.

The reproducing kernel may be written as the (usually infinite) sum
$$
K(x,y)=\sum_{\nu} \phi_{\nu}(x) \phi_{\nu}(y),
$$
where $\{\phi_{\nu}\}$ is a basis for the Hilbert space of integrands
\cite{Wah90}.  This implies that the discrepancy may be written as the
sum \cite{HicWoz00a}
$$
D^2(P) = \sum_{\nu}[\Err(\phi_{\nu},P)]^2.
$$

By Loh's central limit theorem it is known that each $\Err(\phi_{\nu},P)$
is asymptotically distributed as a normal random variable with mean
zero.  However, the $\Err(\phi_{\nu},P)$ are in general correlated.

By making an orthogonal transformation one may write
\begin{equation} \label{apd}
D^2(P)  \approx \sum_{\nu =1}^{\infty}\beta_\nu X_\nu,
\end{equation}
where the $X_{\nu}$ are independent $\chi^2(1)$ random variables, and
the $\beta_{1} \ge \beta_{2} \cdots $ are some constants.  Note that
if 
$\beta_1 = \cdots = \beta_n = \beta$ 
and 
$0 = \beta_{n+1} =\beta_{n+2} = \cdots$ 
then 
$D^2(P)/\beta$ 
is approximately distributed as $\chi^2(n)$.
\end{section}

\begin{section}{Fitting The Distribution Of The Square Discrepancy}
Since equation \eqref{apd} involves an infinite
sum it must be further approximated in order to be computationally
feasible.  This is done by approximating all but the first $n$ terms
by a constant:
\begin{equation} \label{discmod}
D^2(P) \approx \sum_{\nu=1}^{n} \beta_\nu X_\nu + c
\end{equation}
where 
$$
c = E[\sum_{\nu=n+1}^{\infty} \beta_\nu X_\nu].
$$
Here $E$ denotes the expectation.

Let $Z$ denote the random variable given by the square discrepancy of
a scrambled net, and let $W$ denote the random variable given by the
right hand side of \eqref{discmod}.  Let $G_{Z}$ and $G_{W}$ denote
the probability distribution functions of these two random variables,
respectively.  Ideally, one would find the $\beta_{\nu}$ and $c$ that
minimize some loss function involving $G_{Z}$ and $G_{W}$, but these
distributions are not known exactly.  Thus, the optimal $\beta_{\nu}$
and $c$ are found by minimizing
\begin{equation} \label{lossfun}
       \sum_{i=1}^{M}
       [\log(\hat{G}_{Z}^{-1}(p_{i}))-\log(\hat{G}_{W}^{-1}(p_{i})]^2.
\end{equation}
Here, 
$$
p_{i}=(2i-1)/(2M), \text{ for } i=1, \ldots, M,
$$ 
and
$$
\hat{G}_{Z}^{-1}(p_{i})=Z_{(i)},
$$ 
the $i^{\text{th}}$ order statistic
from a sample of $M$ scrambled net square discrepancies.  Also,
$$
\hat{G}_{W}^{-1}(p_{i})=W_{(k_{i})},
$$ 
the $k_{i}^{\text{th}}$ order
statistic from a sample of $L$ independent and identically distributed
drawings of $W$, where 
$$
k_{i}=iL/M - (L-M)/(2M),
$$ 
and $L$ is an odd
multiple of $M$.  The logarithm is used in \eqref{lossfun} because for
a fixed $N$ and $s$ the square discrepancy values can vary by a factor
of 100 or more, and it is not good for the larger values to unduly
influence the fitted distribution.

Fitting the distribution of the square discrepancy by
minimizing \eqref{lossfun} relies on several approximations.  The
central limit theorem is invoked to obtain \eqref{apd}.  The infinite
sum is replaced by a finite one in \eqref{discmod}.  The probability
distributions of each side of \eqref{discmod} are approximated by
Monte Carlo sampling, and the two probability distributions are
compared at only a finite number of points.  In spite of these
approximations the fitted distribution matches the observed
distribution of the square discrepancy rather well.
\end{section}

\begin{section}{Numerical Results}

The discrepancy considered here is \eqref{dis}. With
$\alpha = 2$ and $\gamma = 1$ the discrepancy can be rewritten as
\begin{multline} \label{discex}
     D^2(P) = -1+ \frac{1}{N^2}\sum_{\bfy,\bfy' \in P}^{N} \prod_{j=1}^{s}
     \left\{1 + \gamma \left[B_1(\bfy_{j})B_1(\bfy'_{j}) + \frac 14
     B_2(\bfy_{j})B_2(\bfy'_{j}) \right .  \right .\\
      \left .  \left . - \frac{1}{24} B_{4}(\bbrace{\bfy_{j}-\bfy'_{j}})
      \right] \right \}
\end{multline}

Figures \ref{s2} and \ref{s5} show the empirical distribution (dashed
line) of the square discrepancy of scrambled nets and its fitted value
(thin line).  The fit is based on $M=1000$ replications of the
scrambled Sobol' \cite{BraFox88,HonHic00a}, Niederreiter-Xing,
and Faure nets, $L=11000$, and $n=3$ or $4$.  The optimization was done by 
using MATLAB's {\tt fminsearch}
function.  The fits are good in the middle of the distributions, but
off in the tails.  Moreover the Q-Q plots of the empirical and fitted
distributions bear this out.\\

\begin{figure}[ht]
\begin{center}
\includegraphics[height=8cm]{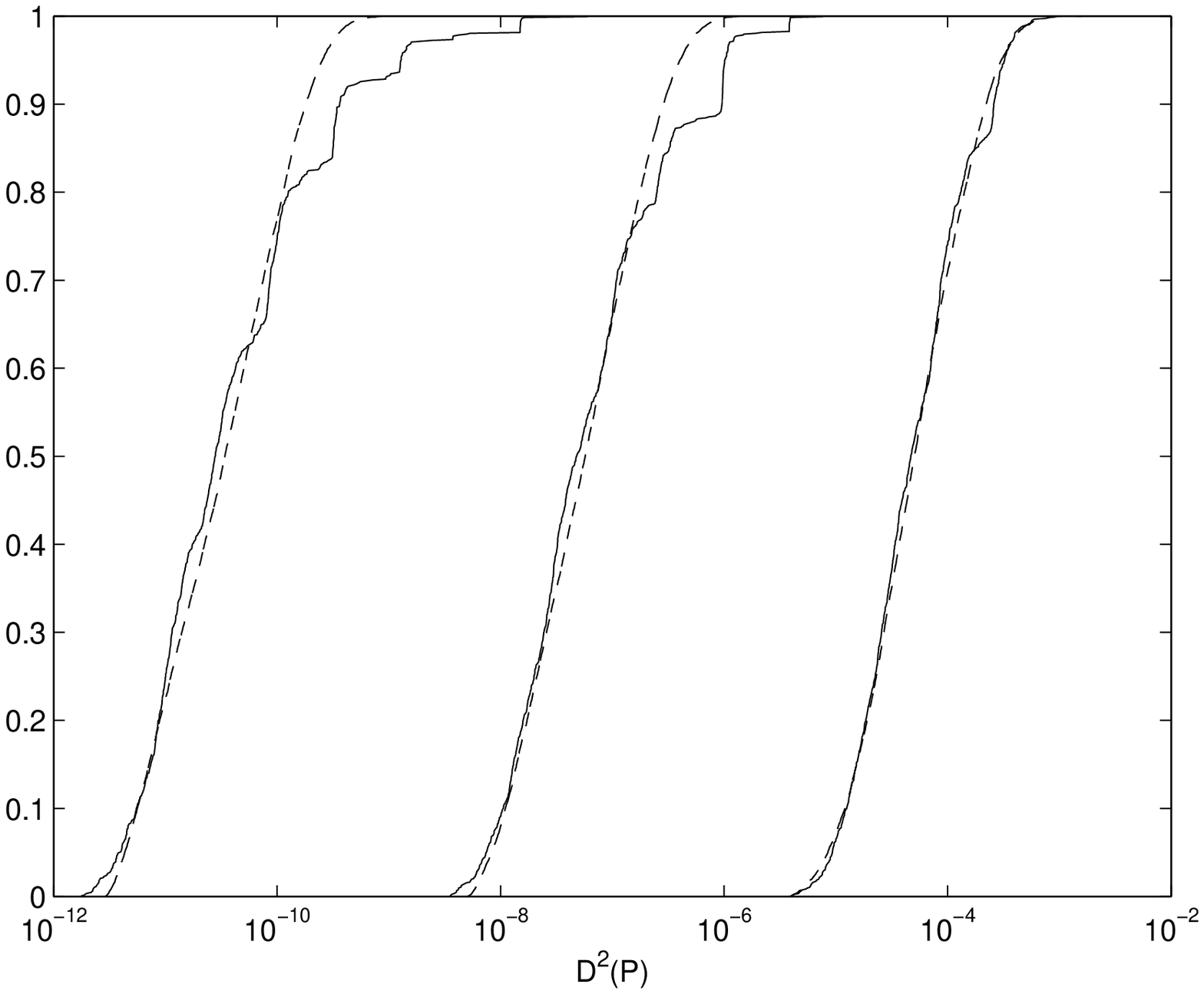}
\end{center}
\caption{\label{s2} The empirical distribution of the square
discrepancy  \eqref{discex} of a scrambled Sobol' net and its fitted
distribution for $s = 2$ and $N = 2^4, 2^7, 2^{10}$ from right to
left.}
\end{figure}

\begin{figure}[ht]
\begin{center}
\includegraphics[height=8cm]{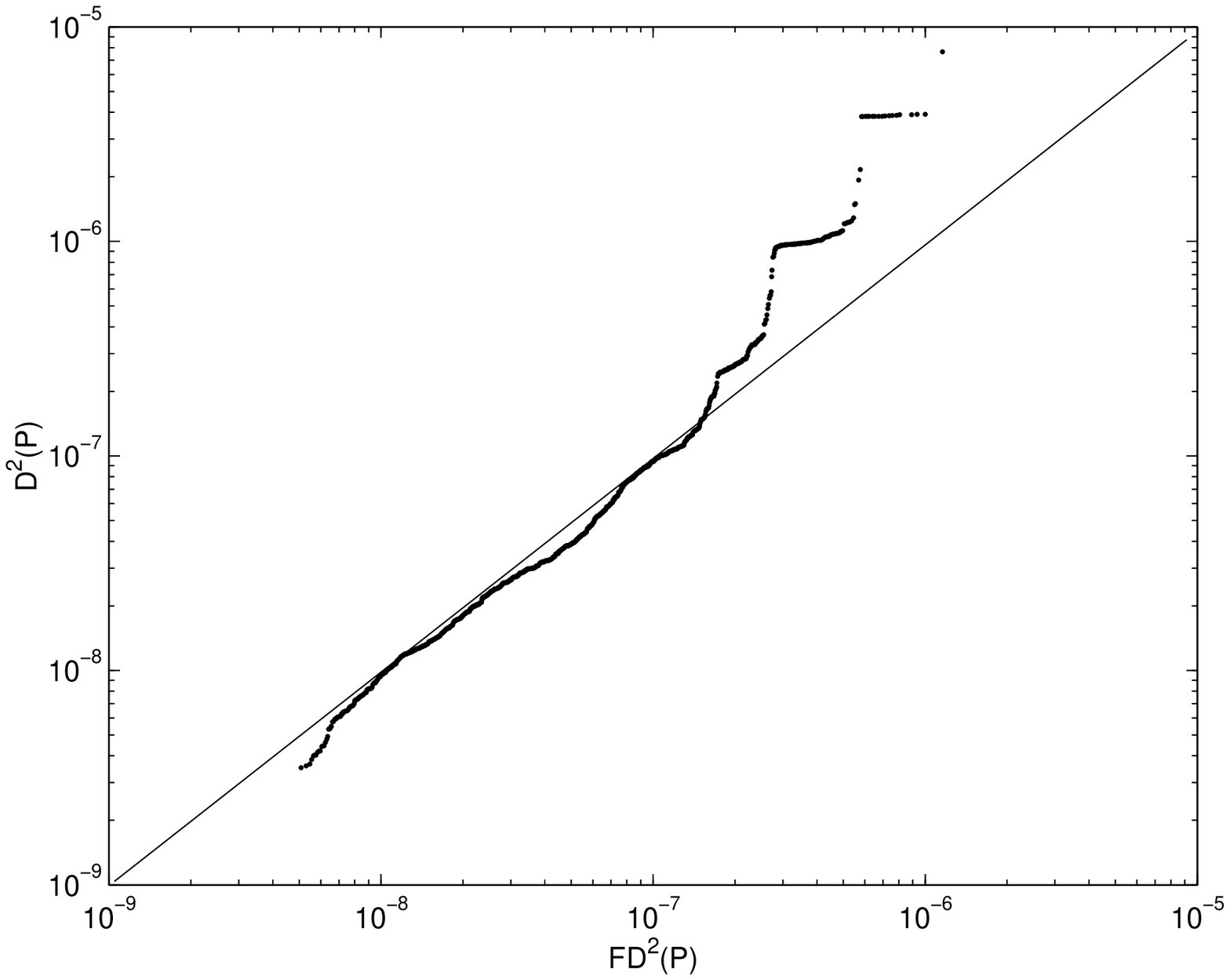}
\end{center}
\caption{\label{s2q7} Q-Q plot of the empirical distribution of the
square discrepancy versus the fitted distribution for a scrambled
Sobol' net with $s=2$ and $N = 2^7$.}
\end{figure}

\begin{figure}[ht]
\begin{center}
\includegraphics[height=8cm]{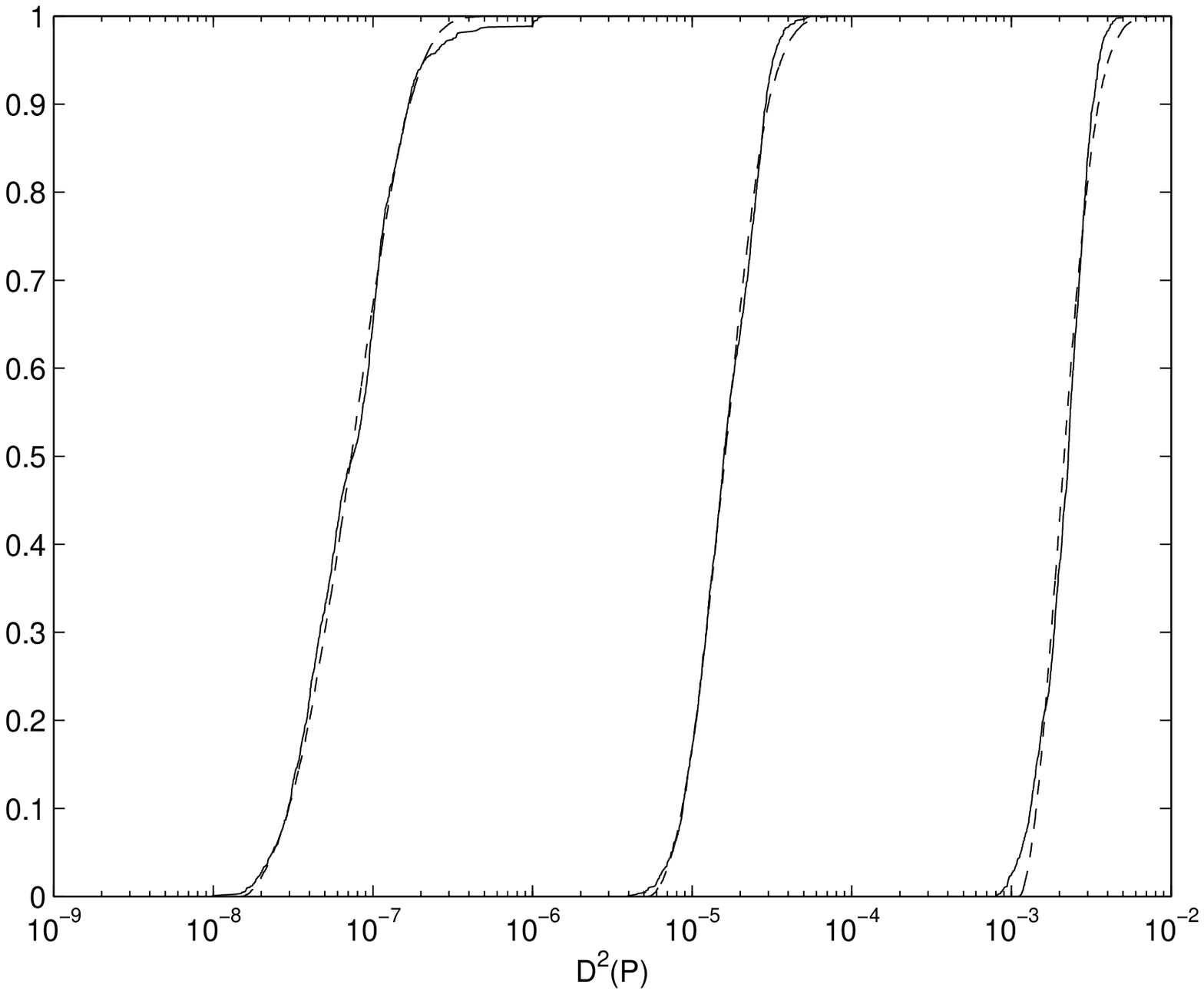}
\end{center}
\caption{\label{s5} The empirical distribution of the square
discrepancy  \eqref{discex} of a scrambled Sobol' net and its fitted
distribution for $s = 5$ and $N = 2^4, 2^7, 2^{10}$ from right to
left. }
\end{figure}
\begin{figure}[ht]
\begin{center}
\includegraphics[height=8cm]{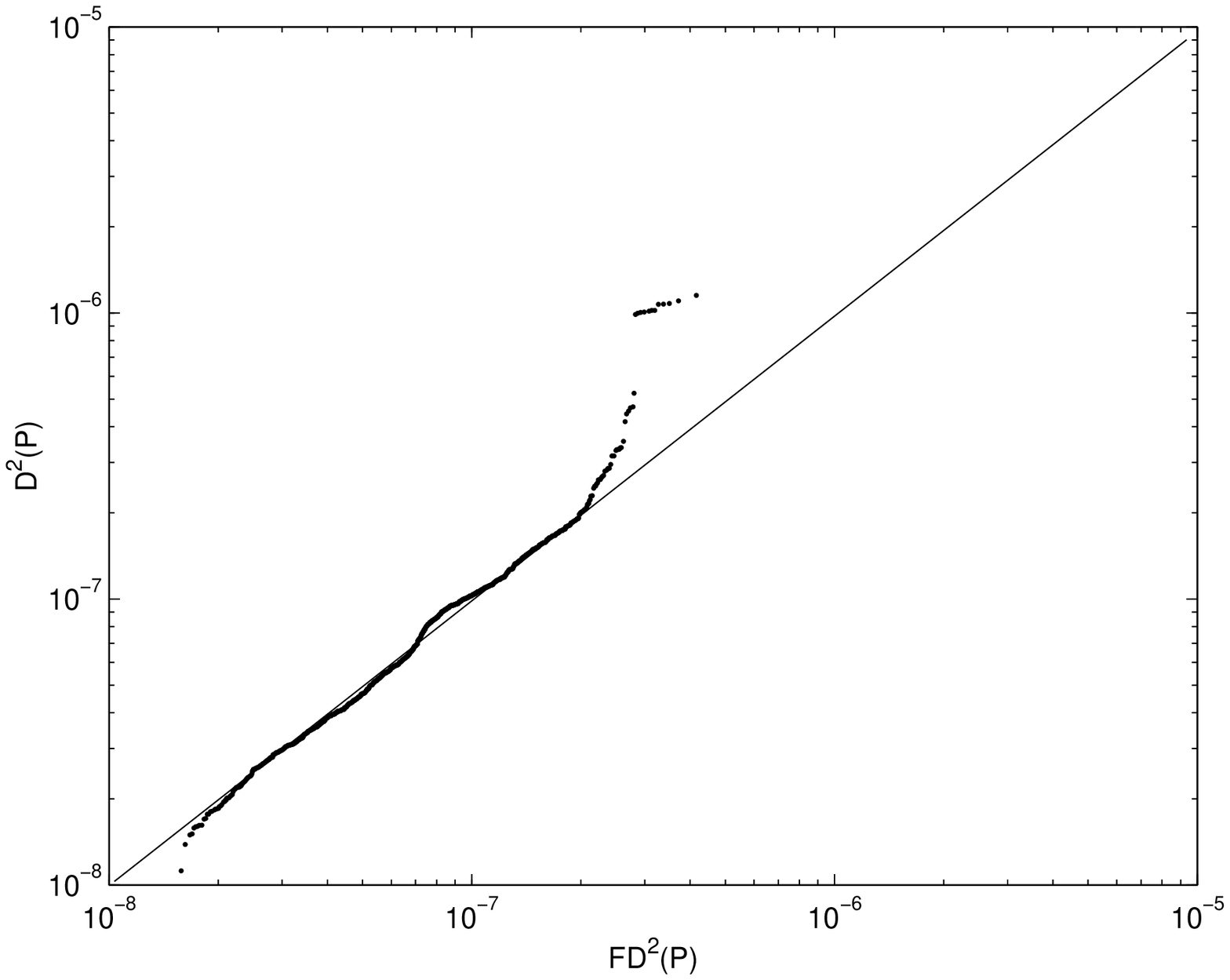}
\end{center}
\caption{\label{s5q10} Q-Q plot of the empirical distribution of the
square discrepancy versus the fitted distribution for a scrambled
Faure' net with $s=3$ and $N = 2^7$.}
\end{figure}

\begin{figure}[ht]
\begin{center}
\includegraphics[height=8cm]{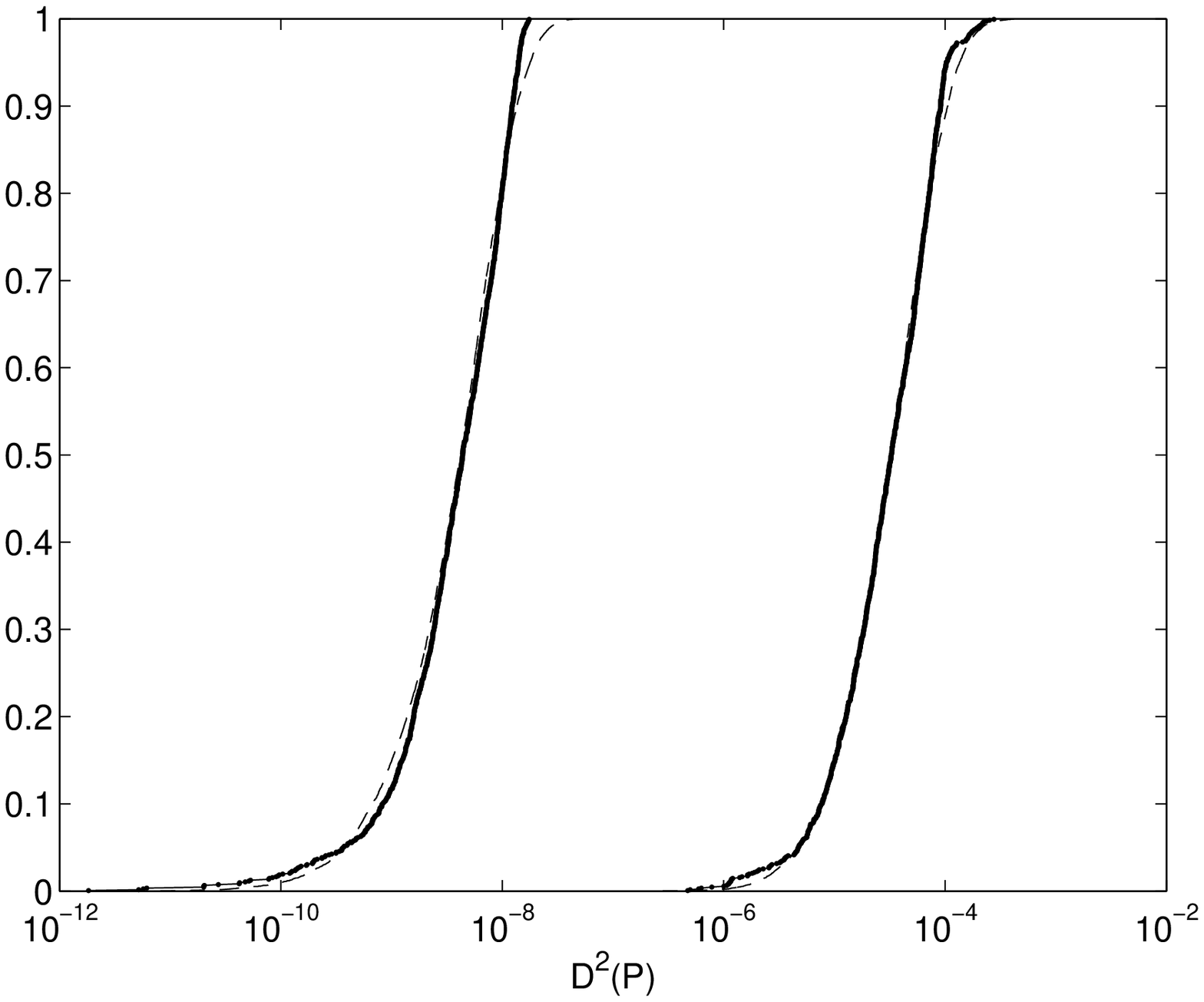}
\end{center}
\caption{\label{sf3} The empirical distribution of the square
discrepancy  \eqref{discex} of a scrambled Faure net and its fitted
distribution for $s = 3$ and $N = 3^3, 3^7$ from right to
left. }
\end{figure}

\begin{figure}[ht]
\begin{center}
\includegraphics[height=8cm]{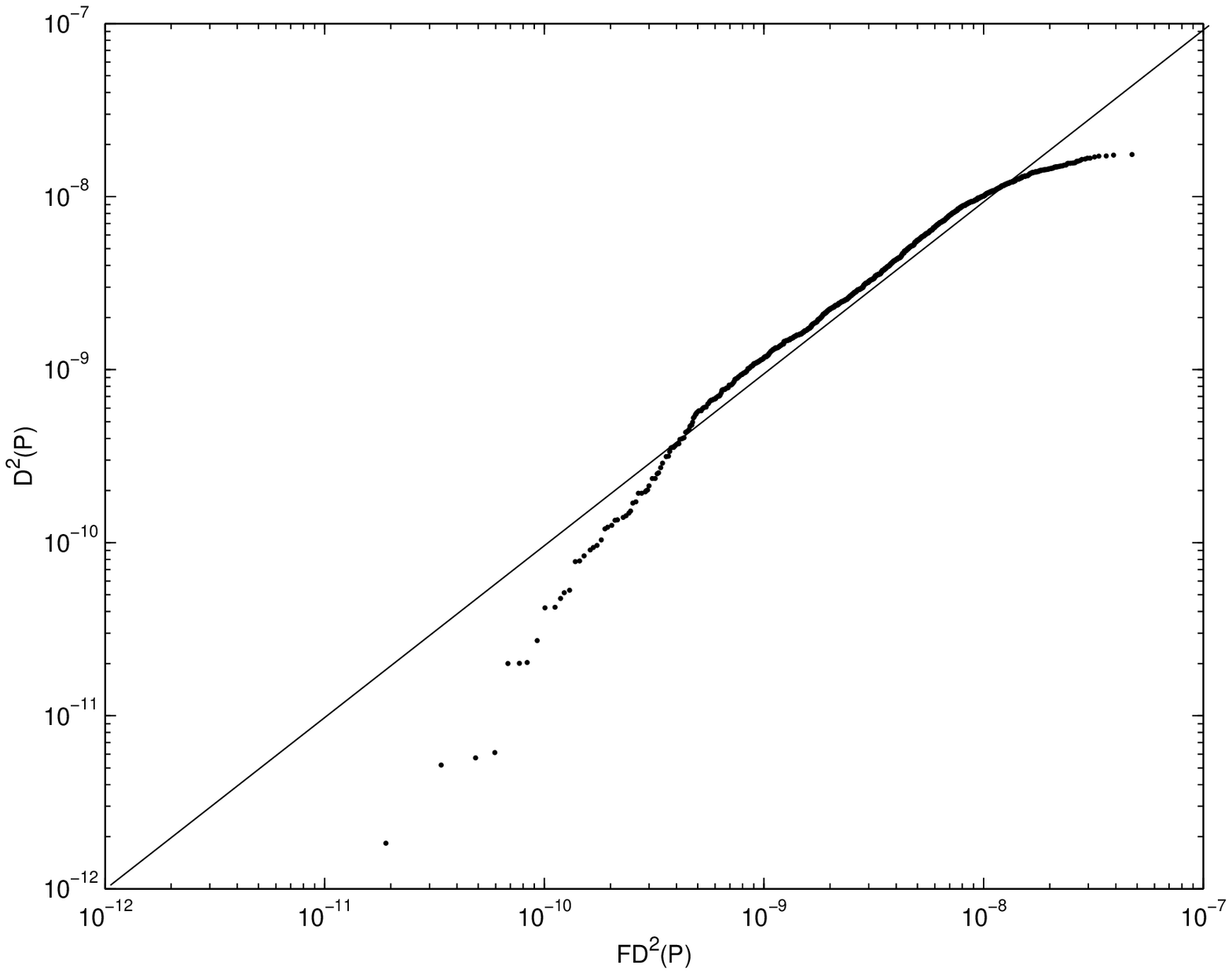}
\end{center}
\caption{\label{f3q7} Q-Q plot of the empirical distribution of the
square discrepancy versus the fitted distribution for a scrambled
Faure net with $s=3$ and $N = 3^7$.}
\end{figure}

\begin{figure}[ht]
\begin{center}
\includegraphics[height=8cm]{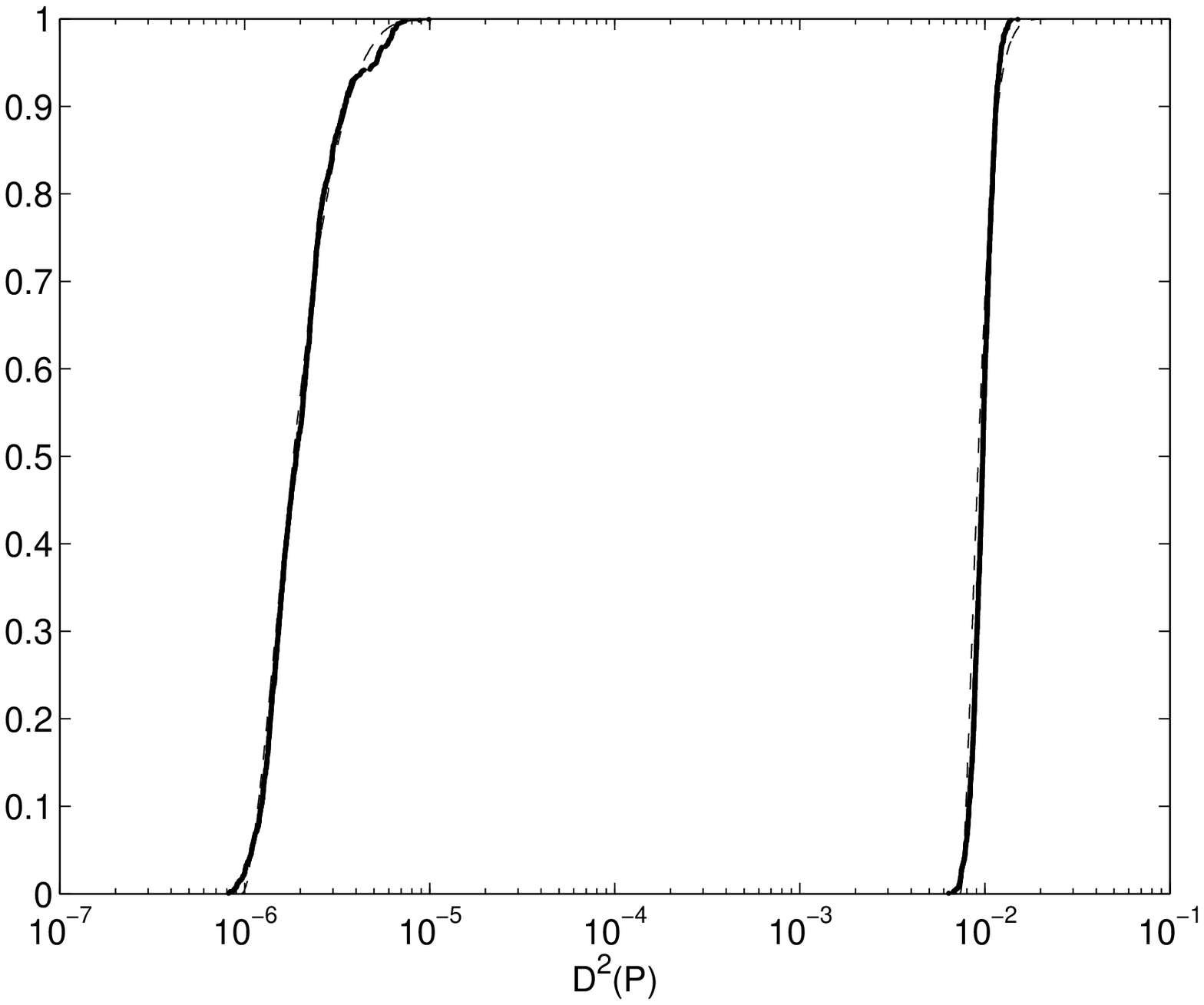}
\end{center}
\caption{\label{snx8} The empirical distribution of the square
discrepancy  \eqref{discex} of a scrambled Niederreiter-Xing net and its fitted
distribution for $s = 2$ and $N = 2^4, 2^{10}$ from right to
left. }
\end{figure}

\begin{figure}[ht]
\begin{center}
\includegraphics[height=8cm]{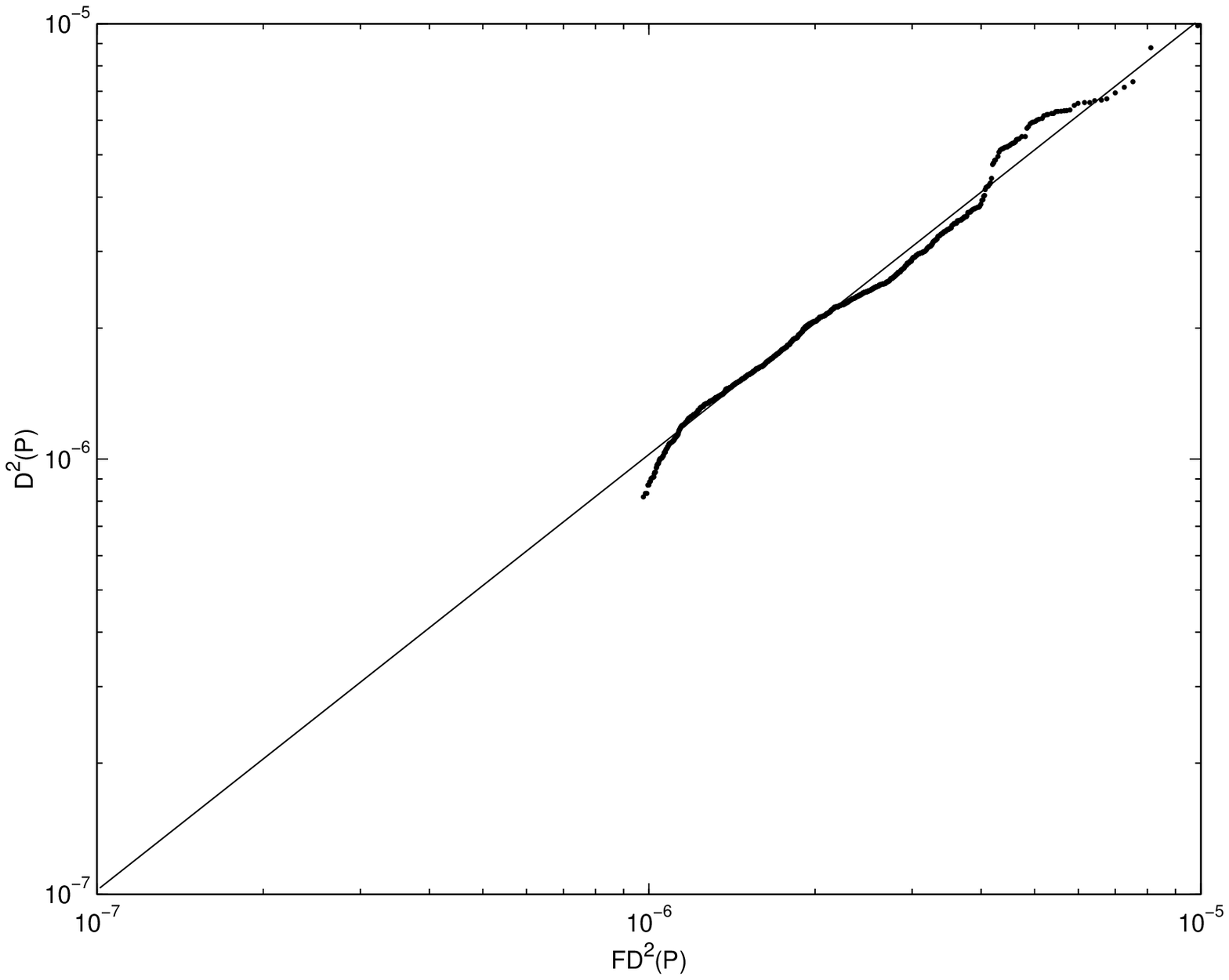}
\end{center}
\caption{\label{nx8q10} Q-Q plot of the empirical distribution of the
square discrepancy versus the fitted distribution for a scrambled
Niederreiter-Xing net with $s=8$ and $N = 2^{10}$.}
\end{figure}

There are several humps in the empirical distribution of the square
discrepancy for $s = 2$.  Moreover as $N$ increases, the number of
humps increases.  However, these humps are less noticeable as the
dimension increases.

The humps can be explained as follows:  The scrambling technique used
here, as described in Chapter 2, scrambles the
generator matrices of these digital nets and also gives them a digital
shift.  Since this technique preserves the digital nature of the nets,
it can be shown that for $l=0, 1, \ldots$ the sum of the digits of 
points numbered $l b^{2}, \ldots, (l+1)b^{2} -1 $, is zero modulo $b$.  
However, for
Owen's original scrambling this is not necessarily the case.  For example, in one
dimension Owen's scrambling is equivalent to Latin hypercube sampling.
Figures \ref{s1} and \ref{lhs1} show the empirical distribution of the
square discrepancy for a scrambled one-dimensional Sobol' net and for
Latin hypercube sampling.  Since the base of the Sobol' sequence is
$2$ the difference in the two graphs emerges for $N \ge b^{2}=4$.
Note that although the distributions of the square discrepancy for
Owen's original scrambling and the variant used here are somewhat
different, the means of the two distributions are the same.
\clearpage
\begin{figure}[t]
\begin{center}
\includegraphics[height=8cm]{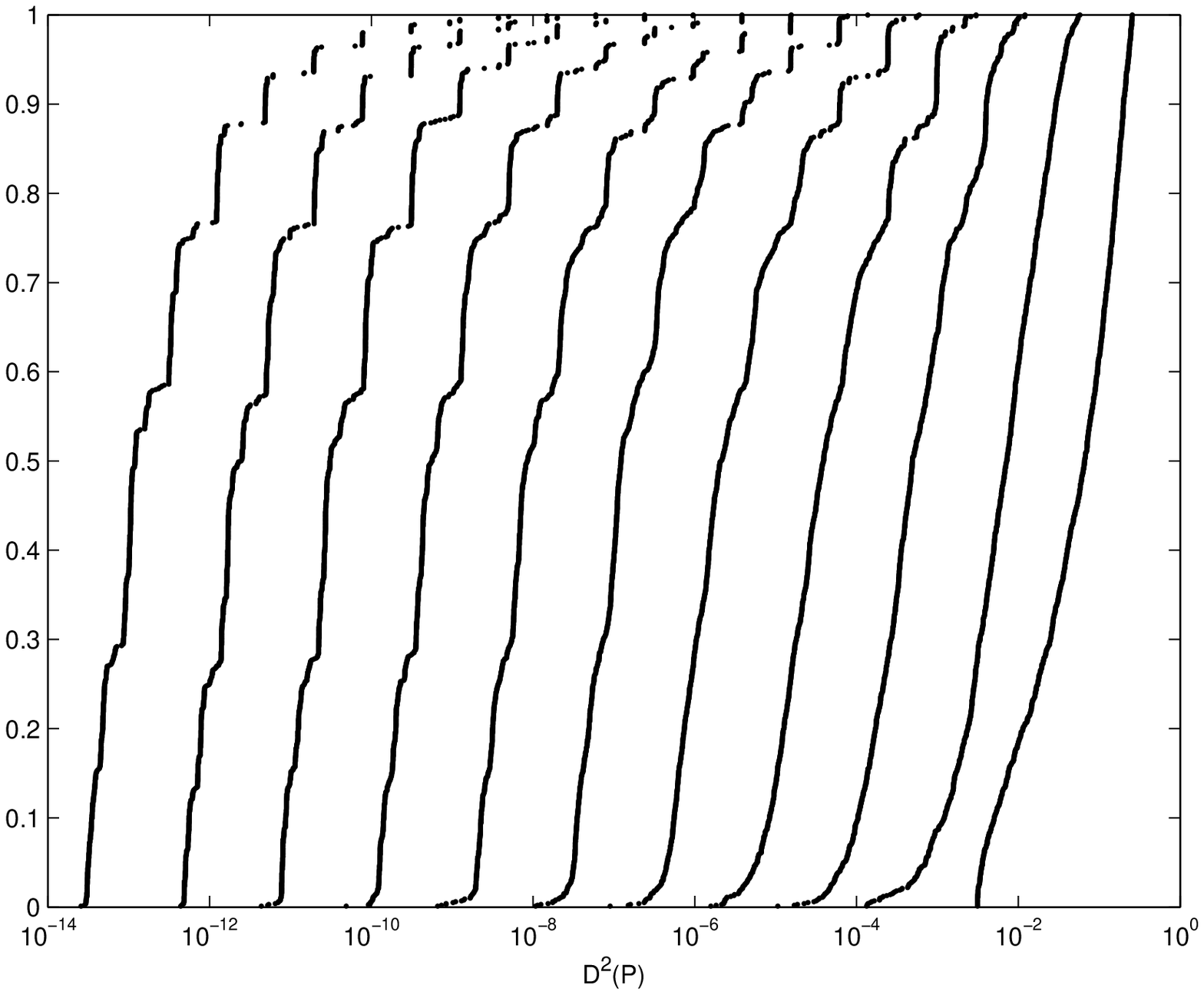}
\end{center}
\caption{\label{s1} The empirical distribution of the square discrepancy
of scrambled Sobol' net for $s = 1$ and $N = 2^0,\cdots,2^{10}$
from right to left.}
\end{figure}
\begin{figure}[b]
\begin{center}
\includegraphics[height=8cm]{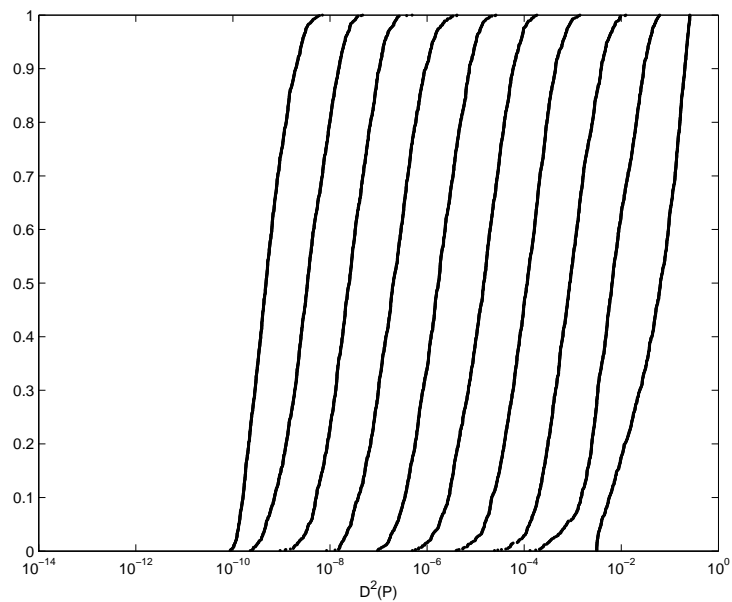}
\end{center}
\caption{\label{lhs1} The empirical distribution of the square discrepancy
of the points from Latin hypercube sampling for $s = 1$ and
$N = 2^0,\cdots,2^{10}$ from right to left.}
\end{figure}

The square discrepancy can be written as a sum of polynomials in
$\gamma$:
\begin{eqnarray*}
D^2(P,K) =& & (\gamma-\gamma_0) \tilde{D}^2_1(P,K)+
(\gamma-\gamma_0)(\gamma-\gamma_1) \tilde{D}^2_2(P,K) \\
& & + \cdots + 
(\gamma-\gamma_0) \cdots (\gamma-\gamma_{s-1}) \tilde{D}^2_s(P,K),
\end{eqnarray*}
Then the sum of polynomials can be rearranged as
$$
D^2(P,K) = \gamma D^2_1(P,K)+
\gamma^2 D^2_2(P,K)+ \cdots + \gamma^s D^2_s(P,K),
$$
where $D^2_j$ measures the uniformity of all $j$-dimensional
projections of $P$ \cite{Hic97a}.  A computationally efficient 
method for obtaining all of the $D^2_j$ is to compute $D^2(P,K)$ 
for $s$ different values of $\gamma_0,\cdots \gamma_{s-1}$ and 
then perform polynomial interpolation. Specifically Newton's divided 
difference formula has been applied for polynomial interpolation.

\begin{figure}[ht]
\begin{center}
\includegraphics[height=12cm]{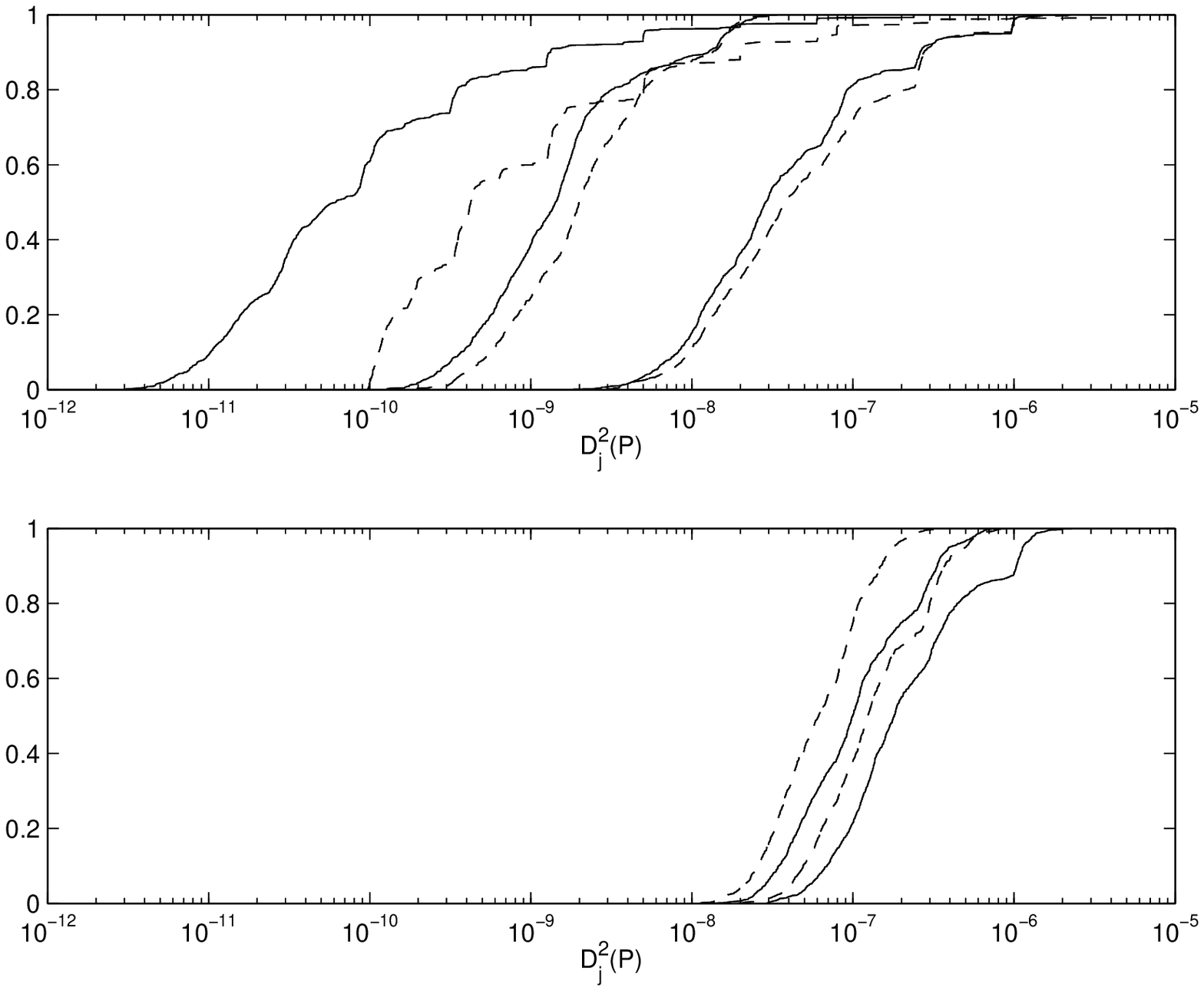}
\end{center}
\caption{\label{ssnx} The empirical distribution of $D^2_j(P)$
for the (3,9,5)-net of scrambled Sobol' (solid) and (2,9,5)-net of 
scrambled Niederreiter-Xing (dashed). Top figure shows $j = 2,5,1$ and
the bottom figure shows $j = 3,4$ from right to left.}
\end{figure}

Figure \ref{ssnx} plots the distribution of $D^2_j$ for
$5$-dimensional scrambled Sobol' and scrambled Niederreiter-Xing nets
\cite{NieXin96,NieXin98a} for $j=1, \ldots, 5$.  The $t$-values for
these nets are 3 and 2 respectively.  From this figure the scrambled Sobol' net has
smaller $D^2_j$ for $j = 1,2,5$ and the scrambled
Niederreiter-Xing sequence has smaller $D^2_j$ for $j = 3,4$.
This means that for integrands that can be well approximated by sums
of one and two-dimensional functions, the Sobol' net will tend to give
small quadrature error even though it has larger $t$ value.

Here the empirical distribution of randomly shifted
lattice points has been investigated as well (see Figures \ref{g5} and
\ref{g5q4}).  Here we use a Korobov type rank-1 lattice with generator
$\eta = 17797$  \cite{HicEtal00}.  Although there is no known theory on the
distribution of the randomly shifted lattice points, the fitted
distribution of the form used for nets seems to work well for large
enough $N$.

\begin{figure}[ht]
\begin{center}
\includegraphics[height=7cm]{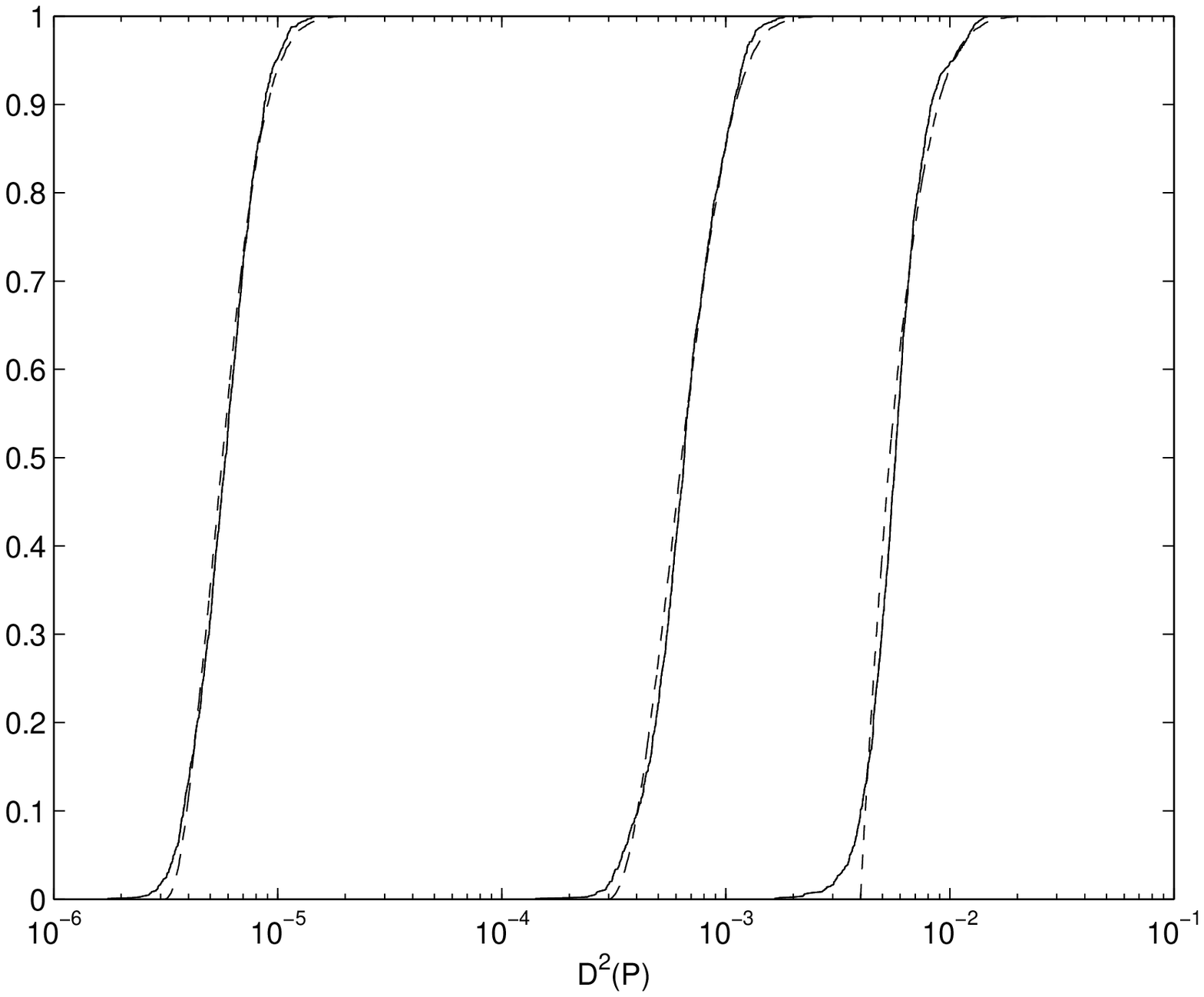}
\end{center}
\caption{\label{g5} The empirical distribution of the square
discrepancy \eqref{discex} of a randomly shifted rank-1 lattice and
its fitted distribution for $s = 5$ and $N = 2^4, 2^7, 2^{10}$ from
right to left.}
\end{figure}

\begin{figure}[ht]
\begin{center}
\includegraphics[height=7cm]{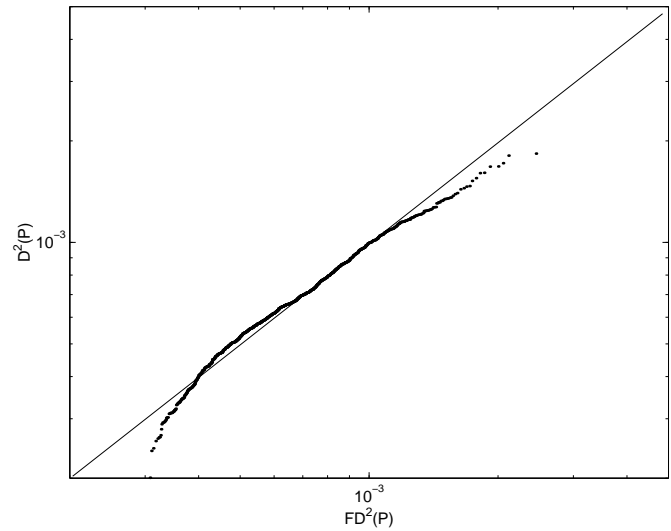}
\caption{\label{g5q4} A Q-Q plot of the empirical distribution of the
squared discrepancy versus the fitted distribution for a randomly shifted
rank-1 lattice with $s=5$ and $N = 2^7$.}
\end{center}
\end{figure}
\end{section}
\end{chapter}

\begin{chapter}{\texorpdfstring{$t$}{t}-parameter Optimization for Digital \texorpdfstring{$(t,m,s)$}{(t,m,s)}-nets by
Evolutionary Computation}
So far most $(t,m,s)$-nets are generated by number theoretic methods.
However, in this chapter 
we use optimization methods to generate a digital $(t,m,s)$-net.
This chapter explains a procedure for finding good generator matrices 
of digital $(t,m,s)$-nets. Here we consider imbedded generator matrices, 
that is one set of generator matrices is used for any $m \leq m_{max}$ and 
$s \leq s_{max}$. 
Many well known generator matrices of digital $(t,m,s)$-nets are imbedded 
matrices, namely Sobol' and Niederreiter. Such 
imbedded matrices have a certain advantage that the user needs only one 
set of matrices rather than several different sets of matrices that work for
different $s$ or $m$. Finding good imbedded nets is more difficult than 
finding nets with fixed $m$ and $s$. 

In Chapter 1 we introduced the definition of $t$ as the quality measure
of $(t,m,s)$-nets and $(t,s)$-sequences in terms of an elementary interval 
in base $b$. In this chapter the quality parameters $t$ are expressed
as the number of linearly independent rows of generator
matrices for digital $(t,m,s)$-nets, which can be computed from the 
generator matrices. Section 4.2 introduces the basic concept of evolutionary 
computation that is used for solving 
optimization problems for finding better nets. 
Section 4.3 explains the detail 
methodology. Finally we report on some numerical results of 
new generator matrices by making comparisons to other well known generator 
matrices of digital $(t,m,s)$-nets.

\begin{section}{Digital \texorpdfstring{$(t,m,s)$}{(t,m,s)}-Nets and \texorpdfstring{$(t,s)$}{(t,s)}-Sequences}

For $b$ a prime number the quality of a digital sequence is expressed 
by its $t$-value, which
can be determined from the generator matrices.  For any positive
integer $m$ let ${\bfc}_{jmk}^{T}$ be the row vector containing the
first $m$ columns of the $k^{\text{th}}$ row of ${\bfC}_j$.  Let $t$
be an integer, $0\leq t\leq m$, such that for all $s$-vectors
${\bfd}=(d_1,\ldots,d_s)$ of non-negative integers with $d_{1} +
\cdots + d_{s}=m-t$,
\begin{multline} \label{tcond}
    \text{the vectors } {\bfc}_{jmk},\ k=1,\ldots,d_j,\ j=1,\ldots,s
    \\
    \text{ are linearly independent over } \Zb.
\end{multline}
where $\Zb$ is a finite field in mod $b$.
Then for any non-negative integer $\nu$ and any $\lambda =0, \ldots,
b-1$ with $\lambda \le b - (\nu \mod b)$, the set $\{\bfy_{\nu b^{m}},
\ldots, \bfy_{(\nu+\lambda)b^{m}-1} \}$ defined by \eqref{digdef} is a
$(\lambda,t,m,s)$-net in base $b$.  (Note that a $(1,t,m,s)$-net is
the same as a $(t,m,s)$-net.)  If the same value of $t$ holds for all
non-negative integers $m$, then the digital sequence is a
$(t,s)$-sequence.

An efficient algorithm is available for the determination of the quality 
parameter of a digital net in the binary case. More detailed explanations
can be found in \cite{Sch99}. 
\end{section}

\begin{section}{Evolutionary Computation}

Evolutionary Computation \cite{Jor00} is an algorithm which is based on the 
principles of a natural evolution as a method to solve parameter optimization
problems.
An evolutionary algorithm (EA) constructs a population of candidate solutions for the
problem at hand. The solutions are evolving by iteratively applying a set 
of stochastic operators,
(or genetic operators) known as recombination, mutation, reproduction and 
selection. During the iterating process (or evolving process) each 
individual in the population receives a measure of its fitness in the 
environment. 

In EA each individual represents a potential solution to 
the problem at hand, and in any evolution program, is implemented as 
some (possibly complex structure) data structure.
Each solution is evaluated to give some measures of its ``fitness", then the new 
population is formed by selecting the more fit individuals. Some members 
of the new population undergo transformations by means of ``genetic" 
operators 
to form a new solution. There is a unary transformation (mutation type), which 
creates new individuals by a small change in a single individual.
Higher 
order transformation (crossover type) creates new individuals by 
combining parts from several individuals. After some number of generations
the program converges to a solution which might be a near-optimum solution. 
Traditionally the data structure of EA is binary representations and 
an operator set consisting only of binary crossover and binary mutation which
is not suitable for many applications. 

However a richer set of data structure (for chromosome representation) 
together with an expanded set of a genetic operators has been 
introduced. The chromosomes need not be represented by bit-strings.
Instead it could be floating number representation. The alternation 
process including other ``genetic" operators appropriates for the given 
structure and the given 
problem. These variations include variables length string, richer than binary 
string. These modified genetic operators meet the needs of particular 
applications. 

The following is a detailed description of stochastic operators in EA.

\begin{itemize}
\item Recombination perturbs the solution by decomposing distinct 
solutions and then randomly mixes their parts to form a novel solution. 
\item Mutation might play important role that introduces further diversity while
the algorithm is running by randomly (or stochastically) perturbing a 
candidate solution, since a large amount of diversity is usually introduced
only at the start of the algorithm by randomizing the genes in the population.
\item Selection evaluates the fitness value of the individuals and
purges poor solutions from population. It resembles the fitness of survivors
in nature.
\item Reproduction focuses attention on high fitness individuals and
replicates the most successful solutions found in a population.
\end{itemize}
The resulting process tends to find globally optimal solutions to the problem 
much in the same way as the natural populations of organisms adapt to their 
surrounding environment. 

A variety of evolutionary algorithms has been proposed. The major ones are
genetic algorithms \cite{Gol89}, evolutionary programming \cite{Fog95}, 
evolutionary strategies \cite{Kur92},
and genetic programming. They all share the common conceptual base of simulating 
the evolution of individual structures via processes of stochastic operators 
which are mentioned earlier. They have been applied to various problems 
where there was no other known problem solving strategy, and the problem domain is
NP-complete. That is the usual place where EAs solve the problem
by heuristically finding solutions where others fail.

\newpage

PSEUDO CODE

Algorithm EA is 

*********************************************************************

\hspace{1cm} // start with an initial time 

\hspace{2cm} Time := 0; 

\hspace{1cm} // initialize a usually random population of individuals 
 
\hspace{2cm} initpopulation Pop; 

\hspace{1cm} // evaluate fitness of all initial individuals in population 

\hspace{2cm}    evaluate Pop; 

\hspace{1cm} // test for termination criterion (time, fitness, etc) 

\hspace{2cm}    While not done do

\hspace{2cm} // select sub-population for offspring production 

\hspace{3cm}  SPop := selectparents Pop; 

\hspace{2cm}     // recombine the "genes" of selected parents

\hspace{3cm}        recombine SPop;

\hspace{2cm}     // perturb the mated population randomly or stochastically

\hspace{3cm}        mutate SPop;  

\hspace{2cm}     // evaluate its new fitness 

\hspace{3cm}        evaluate SPop;
 
\hspace{2cm}     // select the survivors from the actual fitness 

\hspace{3cm}        Pop := survive Pop,SPop;

\hspace{2cm} // increase the time counter 

\hspace{3cm}   Time := Time+1;

*********************************************************************
\newpage

\end{section}

\begin{section}{Searching Method}

Searching for finding good generator matrices can be considered as
solving a combinatorial optimization problem. As the size of $s$ and $m$
increases, the possible solutions of the problem become exponentially
large, $O(b^{m^2s})$. Evolutionary computation
is one of the methods for solving such problem.
In this section we describe the searching procedure for finding good 
digital $(t,m,s)$-nets by an evolutionary computation strategy.

\begin{subsection}{Environment Setting}

Considering generator matrix $\bfC_j$ where $j = 1,\cdots,s$.
\[ \left. \bfC_j   = 
\left( \begin{array}{cccccc} 
1 & C_{j,1,2} & C_{j,1,3} & \cdots & C_{j,1,m}\\
0 & 1 & C_{j,2,3} & \cdots & C_{j,2,m}\\
\vdots & \vdots & \ddots & \vdots &\vdots \\
 0 & \cdots & \cdots &1 &  C_{j,m-1,m} \\ 
 0 & \cdots & \cdots & \cdots& 1 \end{array} \right) \right.\]

First, we define an individual as $\bfA = (\bfA_1,\cdots,\bfA_s)$. 
Each sub population $\bfA_j$ can be 
expressed as $\bfA_j = (A_{j(1)},\cdots,A_{j(r)}) \in \Z_2$ where 
$r$ is the size of cells.
Here we consider the matrix $\bfC_j$ as a $m \times m$ upper triangular matrix
with all diagonal entries are one. Notice that we only consider the case 
$b = 2$. Then the size of $r$ becomes $(mm-m)/2$. The following explains 
the correspondence between $\bfA_j$ and $\bfC_j$.
\begin{eqnarray*}
  A_{j(1)}&=& C_{j,1,2},     \\
  A_{j(2)}&=& C_{j,1,3},     \\
   & \cdots & \\
  A_{j(m-1)}&=& C_{j,1,m},     \\
  A_{j(m)}&=& C_{j,2,3},     \\
   & \cdots & \\
  A_{j(r)} &=& C_{j,m-1,m}. 
\end{eqnarray*}
\end{subsection}

\begin{subsection}{Searching Strategy}
In previous section, we define an individual $\bfA_j$. 
However 
looking at optimization problem as a whole, 
the population becomes $\bfA^i_j$ where $i = 1,\cdots, u$ and $u$ is the 
number of individuals. In this problem we choose $u = 400$.
The followings explain the detail setting of stochastic operators.

(a) Recombination : Two points crossover has been adopted. The 
recombination only occurs for the same $j$ among $\bfA^i_j$s.  
The selection of the locations of crossover is chosen randomly within first
half for the first location and second half for two second location. 

(b) Mutation : 5 \% of randomly selected positions are assigned to perform 
mutation.
 
(c) Selection : Best 30 \% of individuals are selected based on a given
objective function.

(d) Reproduction : The best individuals are selected
from the current generation to remain as the part of the next generation. 
And the rest of the next generation is produced based on the replication of
the best individuals of the current generation. Since we are keeping the 
best previous individuals as themselves, aging has been introduced that each 
individual will be removed when the individuals are survived more than 
certain length of period.
 
\end{subsection}

\begin{subsection}{Objective Function}

There are various objective functions that can be used in this problem.
Here we select an objective function as the sum of $t$ values for $s$
dimensional projections. The details are following.\\
Let $T(m,s,{\bfC}_{1},\cdots,{\bfC}_{s_{max}})$ be the smallest 
$t$ for which ${\bfC}_{1},\cdots,{\bfC}_{s}$ generates a $(t,m,s)$-net
that the $t$ values of the digital $(t,m,s)$-net generated by 
$$
{\bfC}_{1m},\cdots,{\bfC}_{sm},
\hspace{0.7cm} s = 1,\cdots, s_{max} \mbox{ ; } m =1,\cdots,m_{max}
$$
where $\bfC_{jm}$ is the top $m\times m$ sub-matrix.

For given $\bfC_1,\cdots,\bfC_{s_{max}}$, we can use these matrices to
generate $(T(m,s),m,s)$-nets and extend the nets for
$m \leq m_{max}$ and $s \leq s_{max}$. 

Finally the objective function is written as
$$T({\bfC}_{1},\cdots,{\bfC}_{s_{max}}) = \sum_{s=1}^{s_{max}}
\sum_{m=1}^{m_{max}} T(m,s;{\bfC}_{1},\cdots,{\bfC}_{s_{max}}).$$
Our goal is finding the generator matrices that obtain
the smaller $T({\bfC}_{1},\cdots,{\bfC}_{s_{max}})$ values.

\end{subsection}
\begin{figure}[h]
\begin{center}
\includegraphics[height=18cm,width=7.4cm]{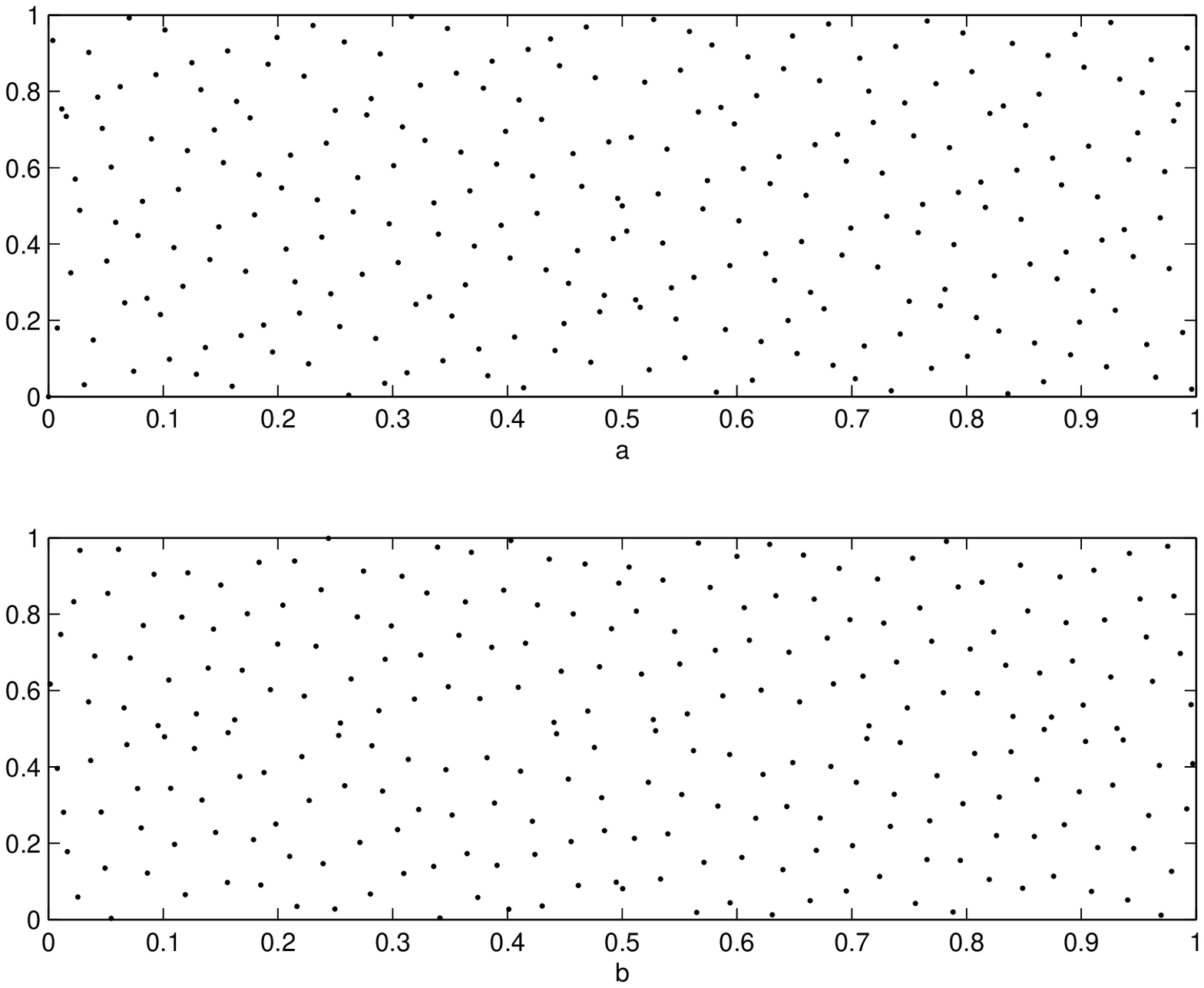}
\end{center}
\caption{\label{psh} Dimensions 16 and 18 of the a) EC and 
b) Scrambled EC' sequences for $N= 256$.}
\end{figure}
\end{section}

\begin{section}{Numerical Results}
The newly obtained digital $(t,m,s)$-net by an evolutionary computation
is named as the EC net.
The exact quality parameters $t$ of the generator matrices were calculated
by using the program provided by Schmid \cite{Sch99}.
Notice that initial populations are selected partially at random and 
partially by taking Sobol' generator matrices and randomly choosing
the $\bfC_j$ matrices.
Table \ref{hongt} shows the $T(m,s)$ values of generator matrices
of EC sequence. Here the $T(m,s)$ is defined as following:
$$ 
T(m,s) = \displaystyle \max_{\left. 
{m' \leq m \atop s'\leq s} \right. }
T(m',s',{\bfC}_{1},\cdots,{\bfC}_{s_{max}}).
$$ 
The same definition is applied for the remaining section.

{\renewcommand{\baselinestretch}{1.2}
\begin{table}
    \caption{\label{hongt} The $T(m,s)$-table for generator matrices of EC}
\begin{center}
\footnotesize
\begin{tabular}{*{22}{c}}
  s & m &1&2&3&4&5&6&7&8&9&10&11&12&13&14&15&16&17&18&19&20\\\hline
  1 &  &0&0&0&0&0&0&0&0&0&0&0&0&0&0&0&0&0&0&0&0\\
  2 &  &0&0&0&0&0&0&0&0&0&0&0&0&0&0&0&0&0&0&0&0\\
  3 &  &0&1&1&1&1&1&1&1&1&1&1&1&1&1&1&1&1&1&1&1\\ 
  4 &  &0&1&2&2&2&3&3&3&3&3&3&3&3&3&3&3&3&3&3&3\\ 
  5 &  &0&1&2&2&2&3&3&3&3&3&4&4&5&5&5&5&5&5&5&5\\ 
  6 &  &0&1&2&2&3&3&4&5&5&5&6&6&6&6&7&7&7&7&7&7\\ 
  7 &  &0&1&2&3&3&3&4&5&5&5&6&6&6&6&7&7&7&7&8&8\\ 
  8 &  &0&1&2&3&3&3&4&5&5&6&7&7&7&8&9&9&9&9&9&9\\ 
  9 &  &0&1&2&3&3&4&4&5&6&6&7&7&7&8&9&9&9&10&11&11\\ 
  10&  &0&1&2&3&3&4&4&5&6&7&7&8&9&9&9&9&9&10&11&11\\ 
  11&  &0&1&2&3&4&5&5&5&6&7&7&8&9&10&10&10&10&11&11&11\\ 
  12&  &0&1&2&3&4&5&5&5&6&7&8&8&9&10&10&10&10&11&12&13\\ 
  13&  &0&1&2&3&4&5&5&6&6&7&8&8&9&10&10&10&10&11&12&13\\ 
  14&  &0&1&2&3&4&5&5&6&7&7&8&8&9&10&11&12&13&13&13&13\\ 
  15&  &0&1&2&3&4&5&5&6&7&7&8&8&9&10&11&12&13&13&13&13\\ 
  16&  &0&1&2&3&4&5&5&6&7&7&8&9&10&10&11&12&13&13&13&13\\ 
  17&  &0&1&2&3&4&5&5&6&7&7&8&9&10&10&11&12&13&13&13&13\\ 
  18&  &0&1&2&3&4&5&5&6&7&7&8&9&10&10&11&12&13&13&14&15\\ 
  19&  &0&1&2&3&4&5&5&6&7&7&8&9&10&10&11&12&13&13&14&15\\ 
  20&  &0&1&2&3&4&5&5&6&7&7&8&9&10&10&11&12&13&13&14&15\\ 
  21&  &0&1&2&3&4&5&6&6&7&7&8&9&10&10&11&12&13&13&14&15\\ 
  22&  &0&1&2&3&4&5&6&6&7&7&8&9&10&10&11&12&13&13&14&15\\ 
  23&  &0&1&2&3&4&5&6&6&7&7&8&9&10&11&12&12&13&13&14&15\\ 
  24&  &0&1&2&3&4&5&6&6&7&7&8&9&10&11&12&12&13&13&14&15\\ 
  25&  &0&1&2&3&4&5&6&6&7&7&8&9&10&11&12&12&13&13&14&15\\ 
  26&  &0&1&2&3&4&5&6&6&7&7&8&9&10&11&12&12&13&13&14&15\\ 
  27&  &0&1&2&3&4&5&6&6&7&7&8&9&10&11&12&12&13&13&14&15\\ 
  28&  &0&1&2&3&4&5&6&6&7&7&8&9&10&11&12&12&13&13&14&15\\ 
  29&  &0&1&2&3&4&5&6&6&7&7&8&9&10&11&12&12&13&13&14&15\\ 
  30&  &0&1&2&3&4&5&6&6&7&7&8&9&10&11&12&12&13&13&14&15\\ 
  31&  &0&1&2&3&4&5&6&6&7&7&8&9&10&11&12&12&13&13&14&15\\ 
  32&  &0&1&2&3&4&5&6&6&7&7&8&9&10&11&12&12&13&13&14&15\\ 
  33&  &0&1&2&3&4&5&6&6&7&7&8&9&10&11&12&12&13&13&14&15\\ 
  34&  &0&1&2&3&4&5&6&6&7&7&8&9&10&11&12&13&14&15&15&15\\ 
  35&  &0&1&2&3&4&5&6&6&7&7&8&9&10&11&12&13&14&15&15&15\\ 
  36&  &0&1&2&3&4&5&6&6&7&7&8&9&10&11&12&13&14&15&15&15\\ 
  37&  &0&1&2&3&4&5&6&6&7&7&8&9&10&11&12&13&14&15&15&15\\ 
  38&  &0&1&2&3&4&5&6&6&7&7&8&9&10&11&12&13&14&15&15&15\\ 
  39&  &0&1&2&3&4&5&6&6&7&8&9&10&10&11&12&13&14&15&15&15\\ 
  40&  &0&1&2&3&4&5&6&7&7&8&9&10&10&11&12&13&14&15&15&15\\\hline 
\end{tabular}
\end{center}
\end{table}
}

Table \ref{hst} shows the comparison of $T(m,s)$
values between generator matrices of Sobol' and EC nets for $s=40$ and $m=20$.
The $T({\bfC}_{1},\ldots,{\bfC}_{s_{max}})$ value for Sobol' sequence is $5421$ 
and EC sequence is $5381$. For the comparison of overall
$T(m,s)$ values,  
if the two methods have the same $T(m,s)$
values then the entry is marked as *.
But if one net is better, which means having a smaller 
$T(m,s)$ value,
then it is marked as the initial of the better net.
Also if the difference of $T(m,s)$ value is greater than $1$ then the difference
is indicated in front of the initial. 
For example, $2E$ means the generator matrices of EC sequence
has a smaller $T(m,s)$ value than Sobol' sequence by $2$. 
From the table \ref{hst}, there are some entries that Sobol'
is better than EC and vice versa. However EC sequence
has smaller $T({\bfC}_{1},\ldots,{\bfC}_{s_{max}})$ values than the 
Sobol' sequence. The EC sequence tends to have smaller $T(m,s)$ values that the
Sobol' sequence for larger $m$ and 
$s$. 

{\renewcommand{\baselinestretch}{1.4}
\begin{table}
 \caption{\label{hst} Comparison of $T(m,s)$ values between generator matrices
  of Sobol' and EC for $s= 40$ and $m = 20$.}    
\begin{center}
\scriptsize
\begin{tabular}{*{22}{@{\extracolsep{3mm}}c}}
  s & m &1&2&3&4&5&6&7&8&9&10&11&12&13&14&15&16&17&18&19&20\\\hline
  10&   &*&*&*&*&*&*&*&*&*&*&S&*&S&S&S&*&*&*&*&*\\
  11&   &*&*&*&*&*&S&S&*&*&S&*&*&*&*&E&2E&2E&E&E&E\\
  12&   &*&*&*&*&*&S&*&*&*&*&*&*&*&*&E&2E&2E&E&*&S\\
  13&   &*&*&*&*&*&S&*&S&*&*&*&*&*&*&E&2E&2E&E&*&*\\ 
  14&   &*&*&*&*&*&S&*&S&S&*&*&*&*&*&*&*&S&*&E&E\\
  15&   &*&*&*&*&*&S&*&S&S&*&*&E&E&*&*&*&S&*&E&E\\
  16&   &*&*&*&*&*&S&*&S&S&*&*&*&*&*&*&*&S&*&E&E\\
  17&   &*&*&*&*&*&S&*&S&S&*&*&*&*&*&*&*&S&*&E&E\\
  18&   &*&*&*&*&*&S&*&*&*&E&*&*&*&E&E&E&E&E&*&S\\
  19&   &*&*&*&*&*&S&*&*&*&E&*&*&*&E&E&E&E&E&*&S\\
  20&   &*&*&*&*&*&S&*&*&*&E&*&*&*&E&E&E&E&E&*&S\\
  21&   &*&*&*&*&*&S&S&*&*&E&*&*&*&E&E&E&E&E&*&S\\
  22&   &*&*&*&*&*&S&S&*&*&E&*&*&*&E&E&E&E&E&*&S\\
  23&   &*&*&*&*&*&S&S&*&*&E&*&*&*&*&*&E&E&E&*&S\\
  24&   &*&*&*&*&*&S&S&*&*&E&*&*&*&*&*&E&E&E&*&S\\
  25&   &*&*&*&*&*&S&S&*&*&E&*&*&*&*&*&E&E&E&*&S\\
  26&   &*&*&*&*&*&S&S&*&*&E&*&*&*&*&*&E&E&E&*&S\\
  27&   &*&*&*&*&*&S&S&*&*&E&*&*&*&*&*&E&E&E&*&S\\
  28&   &*&*&*&*&*&S&S&*&*&E&*&*&*&*&*&E&E&E&*&S\\
  29&   &*&*&*&*&*&S&S&*&*&E&*&*&*&*&*&E&E&E&*&S\\
  30&   &*&*&*&*&*&S&S&*&*&E&*&*&*&*&*&E&E&E&*&S\\
  31&   &*&*&*&*&*&S&S&*&*&E&E&E&E&*&*&E&E&E&*&*\\
  32&   &*&*&*&*&*&S&S&*&*&E&E&E&E&*&*&E&E&E&*&*\\
  33&   &*&*&*&*&*&*&S&*&*&E&E&E&E&*&*&E&E&E&*&*\\
  34&   &*&*&*&*&*&*&*&*&*&E&E&E&E&*&*&*&*&S&S&*\\
  35&   &*&*&*&*&*&*&*&E&*&E&E&E&E&*&*&*&*&S&*&*\\
  36&   &*&*&*&*&*&*&*&E&E&2E&E&E&E&*&*&*&*&S&*&*\\
  37&   &*&*&*&*&*&*&*&E&E&2E&E&E&E&*&*&*&*&*&E&2E\\
  38&   &*&*&*&*&*&*&*&E&E&2E&E&E&E&*&*&*&*&*&E&2E\\
  39&   &*&*&*&*&*&*&*&E&E&E&*&*&E&*&*&*&*&*&E&2E\\
  40&   &*&*&*&*&*&*&*&*&E&E&*&*&E&*&*&*&*&*&E&2E\\\hline
\end{tabular}
\end{center}
\end{table}
}

We also compare the generator matrices of Niederreiter-Xing sequence
with EC sequence. However, notice that the generator 
matrices of Niederreiter-Xing sequence are 
originally not constructed as imbedded matrices like EC and Sobol' sequences do. 
Therefore there are different generator matrices available for different
$s$. Here we compute the $T(m,s)$ values of the generator matrices of 
Niederreiter-Xing like other imbedded matrices. 
Two matrices are chosen for comparison, ``nxs18m30" for $s = 18$ and 
``nxs32m40" for $s = 32$. The generator matrices are obtained from 
\cite{Pir02}.

From tables \ref{hnt18} and \ref{hnt32}, the EC sequence obtains smaller 
$T(m,s)$ values than Niederreiter-Xing sequence, and the differences are 
larger for smaller $s$. It is understandable that the generator matrices 
are constructed best for fixed dimension like $s= 18$ and $32$ in this 
example. However the 
tables show EC sequence even obtains smaller $T(m,s)$ values for $s_{max}$.

{\renewcommand{\baselinestretch}{1.4}
\begin{table}
 \caption{\label{hnt18}Comparison of $T(m,s)$ values between generator 
 matrices of Niederreiter-Xing and EC for $s= 18$ and $m = 20$.}
\begin{center}
\scriptsize
\begin{tabular}{*{22}{@{\extracolsep{3mm}}c}}
  s & m &1&2&3&4&5&6&7&8&9&10&11&12&13&14&15&16&17&18&19&20\\\hline
  1 &   &E&2E&3E&4E&5E&6E&6E&6E&6E&6E&6E&6E&6E&6E&6E&6E&6E&6E&6E&6E\\
  2 &   &E&2E&3E&4E&5E&6E&6E&6E&6E&6E&7E&8E&8E&8E&8E&8E&8E&8E&8E&8E\\
  3 &   &E&E&2E&3E&4E&5E&5E&5E&5E&5E&6E&7E&7E&7E&7E&7E&7E&7E&7E&7E\\
  4 &   &E&E&E&2E&3E&3E&3E&3E&3E&3E&4E&5E&5E&6E&6E&6E&6E&6E&6E&6E\\
  5 &   &E&E&E&2E&3E&3E&3E&3E&3E&3E&3E&4E&3E&4E&4E&4E&4E&5E&5E&5E\\
  6 &   &E&E&E&2E&2E&3E&2E&E&E&E&E&2E&2E&3E&3E&4E&4E&4E&4E&4E\\
  7 &   &E&E&E&E&2E&3E&2E&E&E&E&E&2E&2E&3E&3E&4E&4E&4E&4E&4E\\
  8 &   &E&E&E&E&2E&3E&2E&E&E&E&E&E&E&E&E&2E&2E&3E&3E&3E\\
  9 &   &E&E&E&E&2E&2E&2E&E&*&E&E&E&2E&E&E&2E&2E&2E&E&E\\
  10&   &E&E&E&E&2E&2E&2E&E&*&*&E&*&*&*&E&2E&3E&3E&3E&4E\\
  11&   &E&E&E&E&E&E&E&E&*&*&E&*&*&N&*&E&2E&2E&3E&4E\\
  12&   &E&E&E&E&E&E&E&E&E&E&E&2E&2E&2E&3E&4E&4E&3E&2E&2E\\
  13&   &E&E&E&E&E&E&E&*&E&E&E&2E&2E&2E&3E&4E&4E&3E&2E&2E\\
  14&   &E&E&E&E&E&E&E&E&*&E&E&2E&2E&2E&2E&2E&E&E&E&2E\\
  15&   &E&E&E&E&E&E&E&E&*&E&E&2E&2E&2E&2E&2E&E&E&E&2E\\
  16&   &E&E&E&E&E&E&E&E&*&E&E&E&E&2E&2E&2E&E&E&E&2E\\
  17&   &E&E&E&E&E&E&E&E&*&E&E&E&E&2E&2E&2E&E&E&E&2E\\
  18&   &E&E&E&E&E&E&E&E&*&E&E&E&E&2E&2E&2E&E&E&*&*\\\hline

\end{tabular}
\end{center}
\end{table}
}

{\renewcommand{\baselinestretch}{1.4}
\begin{table}
    \caption{\label{hnt32} Comparison of $T(m,s)$ values between generator matrices
  of Niederreiter-Xing and EC for $s= 32$ and $m = 20$.}
\begin{center}
\scriptsize
\begin{tabular}{*{22}{@{\extracolsep{3mm}}c}}
  s & m &1&2&3&4&5&6&7&8&9&10&11&12&13&14&15&16&17&18&19&20\\\hline
  1 &   &E&2E&3E&4E&4E&4E&4E&4E&4E&4E&4E&4E&4E&4E&4E&4E&4E&4E&4E&4E\\
  2 &   &E&2E&3E&4E&4E&4E&4E&4E&4E&4E&4E&4E&4E&4E&4E&4E&4E&4E&4E&4E\\
  3 &   &E&E&2E&3E&3E&3E&3E&3E&4E&4E&4E&4E&4E&5E&5E&5E&5E&5E&5E&5E\\
  4 &   &E&E&E&2E&2E&2E&3E&3E&3E&3E&4E&5E&6E&6E&6E&6E&6E&6E&6E&6E\\
  5 &   &E&E&E&2E&2E&2E&3E&3E&3E&3E&3E&4E&4E&4E&4E&4E&4E&4E&4E&4E\\
  6 &   &E&E&E&2E&E&2E&2E&2E&3E&3E&2E&2E&3E&3E&2E&2E&3E&4E&5E&6E\\
  7 &   &E&E&E&E&E&2E&2E&2E&3E&3E&2E&3E&3E&3E&2E&2E&3E&4E&4E&5E\\
  8 &   &E&E&E&E&E&2E&2E&2E&3E&2E&E&2E&2E&E&*&*&E&2E&3E&4E\\
  9 &   &E&E&E&E&E&E&2E&2E&2E&2E&E&2E&2E&E&*&*&E&E&E&2E\\
  10&   &E&E&E&E&E&E&2E&2E&2E&E&E&E&*&*&*&E&E&E&E&2E\\
  11&   &E&E&E&E&*&*&E&2E&2E&E&E&E&*&N&*&*&*&*&E&2E\\
  12&   &E&E&E&E&*&*&E&2E&2E&E&*&E&*&N&*&*&*&*&*&*\\
  13&   &E&E&E&E&*&*&E&E&2E&E&*&E&*&N&*&*&E&E&*&*\\
  14&   &E&E&E&E&*&*&E&E&E&E&*&E&*&*&N&2E&2E&N&N&*\\
  15&   &E&E&E&E&*&*&E&E&E&E&*&E&E&E&*&N&N&N&*&E\\
  16&   &E&E&E&E&*&*&E&E&E&E&*&*&*&E&*&N&N&N&*&E\\
  17&   &E&E&E&E&*&*&E&E&E&E&*&*&*&E&*&N&N&*&E&E\\
  18&   &E&E&E&E&*&*&E&E&E&E&*&*&*&E&*&N&N&*&*&N\\
  19&   &E&E&E&E&*&*&E&E&E&E&*&*&*&E&*&N&N&*&*&N\\
  20&   &E&E&E&E&*&*&E&E&E&E&*&*&*&E&*&N&N&*&*&N\\
  21&   &E&E&E&E&*&*&*&E&E&E&*&*&*&E&*&N&N&*&*&N\\
  22&   &E&E&E&E&*&*&*&E&E&2E&E&*&*&E&*&N&N&*&*&N\\
  23&   &E&E&E&E&*&*&*&E&E&2E&E&*&*&*&N&*&N&*&*&*\\
  24&   &E&E&E&E&*&*&*&E&E&2E&E&*&*&*&N&*&N&*&*&*\\
  25&   &E&E&E&E&*&*&*&E&E&2E&E&*&*&*&*&E&*&*&*&*\\
  26&   &E&E&E&E&*&*&*&E&E&2E&E&*&*&*&*&E&*&*&*&*\\
  27&   &E&E&E&E&*&*&*&E&E&2E&E&*&*&*&*&E&*&*&*&*\\
  28&   &E&E&E&E&*&*&*&E&E&2E&E&*&*&*&*&E&*&*&*&*\\
  29&   &E&E&E&E&*&*&*&E&E&2E&E&*&*&*&*&E&*&E&*&*\\
  30&   &E&E&E&E&*&*&*&E&E&2E&E&*&*&*&*&E&*&E&*&*\\
  31&   &E&E&E&E&*&*&*&E&E&2E&E&*&*&*&*&E&*&E&*&*\\
  32&   &E&E&E&E&*&*&*&E&E&2E&E&*&*&*&*&E&*&E&*&*\\\hline
\end{tabular}
\end{center}
\end{table}
}

We also investigate the $T(m,s)$ values of the lower dimensional 
projections for EC and Niederreiter-Xing sequences. Tables 
\ref{hnlt18} and \ref{hnlt32} show the three different dimensional 
projections. 
\begin{eqnarray*}
T^{avg}_r(\bfC_1,\ldots,\bfC_{s_{max}}) &=& 
\frac{1}{\left( {s_{max}\atop r} \right)} 
\sum_{\left. {u \subseteq 1:s_{max}}\atop {|u| = r} \right.}
T(m,r;\bfC_j, j \in u) \\
T^{max}_r(\bfC_1,\ldots,\bfC_{s_{max}})  &=& 
\displaystyle 
\max_{\left. {u \subseteq 1:s_{max}}\atop {|u| = r} \right.}
T(m,r;\bfC_j, j \in u)
\end{eqnarray*}
where $r = 1,2$, and $18$ for $s_{max} = 18$ and $r = 1,2$, and $32$ for 
$s_{max} = 32$.

From tables \ref{hnlt18} and \ref{hnlt32} the generator matrices of EC has smaller $t$ values for
most entries. It implies that EC sequences has better equidistribution 
property for lower as well as higher dimensional projections.
The $t$ values of $T(m,18,NX_{Prop})$ and $T(m,32,NX_{Prop})$
are taken from \cite{Pir01,Pir02}. Notice that $T(m,18,NX_{Prop})$
values are usually smaller than $T(m,18,NX)$, because $T(m,18,NX_{Prop})$
is obtained by applying propagation rule \cite{Nie99} to the generator 
matrices.

\begin{table}
    \caption{\label{hnlt18} The $T^{avg}_r$ and $T^{max}_r$
values for $r = 1,2$ and $18$ of generator matrices of 
Niederreiter-Xing and EC for $s_{max}= 18$, and 
$m_{max} = 20$.}
\begin{center}
\footnotesize
\begin{tabular}{rrrrrrrrcccccccccccccccccccccccc} 
 m &1&2&3&4&5&6&7&8&9&10\\\hline 
 $T_1^{avg}(EC)$&0&0&0&0&0&0&0&0&0&0 \\
 $T_1^{avg}(NX)$&0.4&0.6&0.9&1&1.3&1.6&1.2&1.4&1.6&1.3 \\
 $T_1^{max}(EC)$&0&0&0&0&0&0&0&0&0&0 \\
 $T_1^{max}(NX)$&1&2&3&4&5&6&5&6&5&3\\
 $T_2^{avg}(EC)$&0&0.5&0.9&1.2&1.6&1.8&1.9&2.2&2.4&2.5\\  
 $T_2^{avg}(NX)$&0.6&1&1.6&2&2.6&3.1&3&3.4&3.4&3.2\\
 $T_2^{max}(EC)$&0&1&2&3&4&5&5&6&7&8\\
 $T_2^{max}(NX)$&1&2&3&4&5&6&6&7&7&8\\
 $T(m,18,EC)$&0&1&2&3&4&5&5&6&7&7\\ 
 $T(m,18,NX)$&1&2&3&4&5&6&6&7&7&8\\
 $T(m,18,NX_{Prop})$&1&2&3&4&4&5&6&7&7&8\\\hline
 m&11&12&13&14&15&16&17&18&19&20\\
 $T_1^{avg}(EC)$&0&0&0&0&0&0&0&0&0&0 \\
 $T_1^{avg}(NX)$&1&0.9&0.9&0.9&0.7&0.7&1&1.2&1.2&1.1\\
 $T_1^{max}(EC)$&0&0&0&0&0&0&0&0&0&0 \\
 $T_1^{max}(NX)$&4&4&4&5&4&5&6&7&8&7\\
 $T_2^{avg}(EC)$&2.8&2.9&3.2&3.2&3.3&3.4&3.5&3.6&3.6&3.7\\
 $T_2^{avg}(NX)$&3.5&3.7&3.9&4.2&3.9&3.9&4.2&4.7&5.0&4.9\\
 $T_2^{max}(EC)$&7&7&8&7&8&9&9&10&11&11\\
 $T_2^{max}(NX)$&9&10&11&12&13&14&11&12&10&11\\
 $T(m,18,EC)$&8&9&10&10&11&12&13&13&14&15\\
 $T(m,18,NX)$&9&10&11&12&13&14&14&14&14&15\\
 $T(m,18,NX_{Prop})$&9&9&10&11&12&12&12&13&14&15\\\hline
\end{tabular}
\end{center}
\end{table}

\begin{table}
    \caption{\label{hnlt32} The $T^{avg}_r$ and
    $T^{max}_r$ values for $r = 1,2$ and $32$ of 
generator matrices of Niederreiter-Xing and EC for $s_{max}= 32$, 
and $m_{max} = 20$.}
\begin{center}
\footnotesize
\begin{tabular}{rrrrrrrrcccccccccccccccccccccccc} 
 m &1&2&3&4&5&6&7&8&9&10\\\hline 
 $T_1^{avg}(EC)$&0&0&0&0&0&0&0&0&0&0 \\
 $T_1^{avg}(NX)$&0.4&0.6&0.9&0.8&0.6&0.8&1&1.1&1.1&1.1\\
 $T_1^{max}(EC)$&0&0&0&0&0&0&0&0&0&0 \\
 $T_1^{max}(NX)$&1&2&3&4&3&4&5&4&4&5\\
 $T_2^{avg}(EC)$&0&0.5&0.8&1.2&1.6&1.9&2.1&2.4&2.6&2.8\\  
 $T_2^{avg}(NX)$&0.9&1.5&1.9&2.3&2.3&2.5&3&3.2&3.5&3.7\\
 $T_2^{max}(EC)$&0&1&2&3&4&5&6&6&7&7\\
 $T_2^{max}(NX)$&1&2&3&4&4&5&6&7&8&9\\
 $T(m,32,EC)$&0&1&2&3&4&5&6&6&7&7\\ 
 $T(m,32,NX)$&1&2&3&4&4&5&6&7&8&9\\
 $T(m,32,NX_{Prop})$&1&2&3&4&4&5&6&7&8&9\\\hline
 m&11&12&13&14&15&16&17&18&19&20\\
 $T_1^{avg}(EC)$&0&0&0&0&0&0&0&0&0&0 \\
 $T_1^{avg}(NX)$&0.8&1&1.1&1&0.9&0.8&0.3&0.6&0.7&1\\
 $T_1^{max}(EC)$&0&0&0&0&0&0&0&0&0&0 \\
 $T_1^{max}(NX)$&4&5&5&6&4&4&3&4&3&4\\
 $T_2^{avg}(EC)$&3&3.2&3.3&3.3&3.5&3.6&3.8&3.8&3.9&4\\
 $T_2^{avg}(NX)$&3.5&3.7&3.9&4.2&4.0&4.3&4.3&4.5&4.8 4.9\\
 $T_2^{max}(EC)$&8&8&9&10&11&12&11&11&12&12\\
 $T_2^{max}(NX)$&9&9&9&10&11&12&11&11&12&13\\
 $T(m,32,EC)$&8&9&10&11&12&12&13&13&14&15\\
 $T(m,32,NX)$&9&9&10&11&12&13&13&14&14&15\\
 $T(m,32,NX_{Prop})$&9&9&10&11&12&13&13&14&14&15\\\hline
\end{tabular}
\end{center}
\end{table}

Figures \ref{hnxs} and \ref{hnxns} plot the root mean square
discrepancy of randomly scrambled and unscrambled Niederreiter-Xing and 
EC sequences. The discrepancy is the same discrepancy used in 
Chapter 2. Also 100 random replications are performed for the scrambled one.
The choices of dimension are $s=18$ and $s = 32$. The number of scrambled 
digits is chosen to be 31. From the figures for the unscrambled case the EC
sequence has a smaller discrepancy than the Niederreiter-Xing sequence. 
But the convergence rate of the discrepancy seems to be nearly the same. 
However for the scrambled case the EC sequence shows better convergence 
rate than the Niederreiter-Xing sequence.

\begin{figure}[ht]
\begin{center}
\includegraphics[height=8cm]{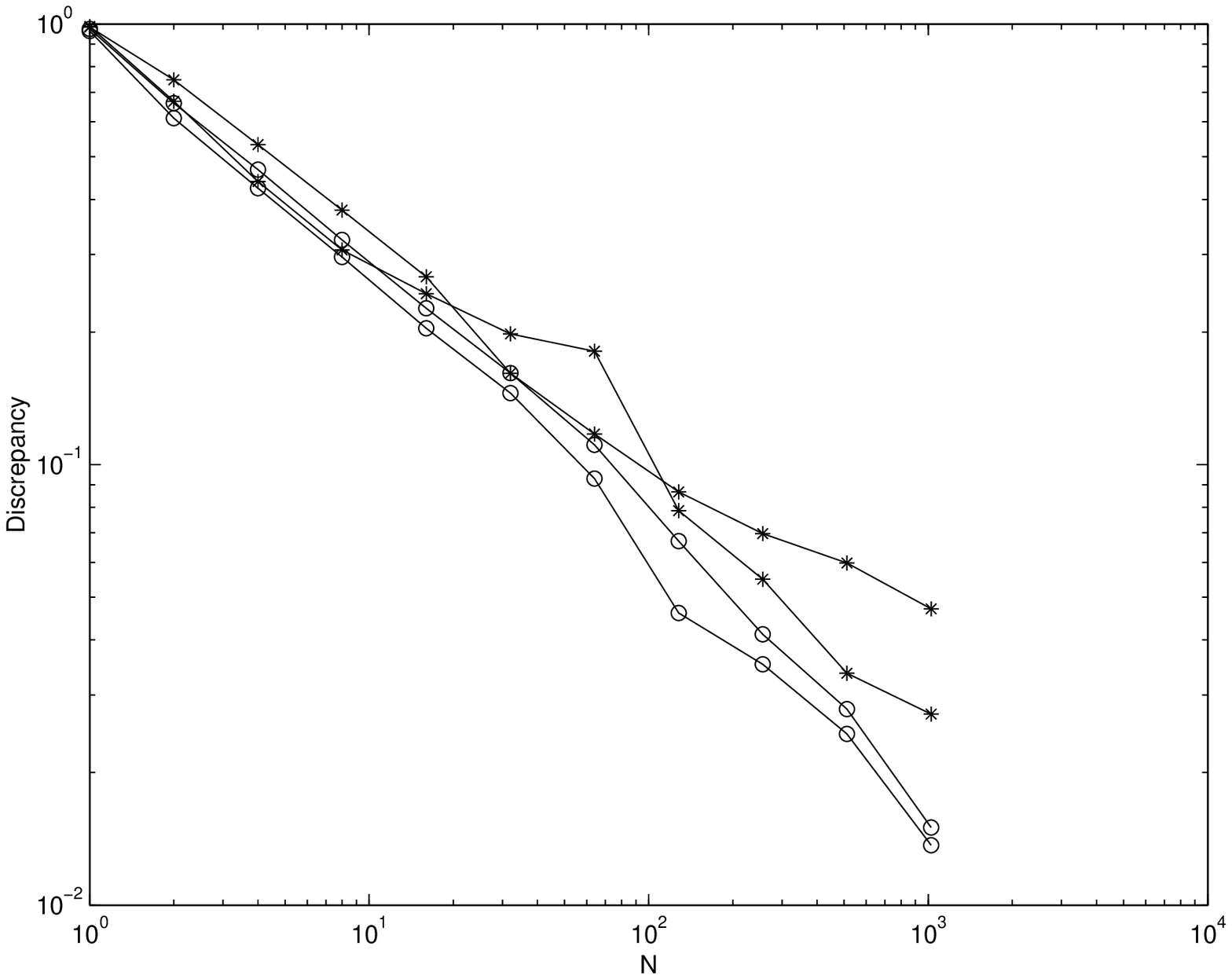}
\end{center}
\caption{\label{hnxs} The scaled root mean square discrepancy with
$\alpha = 2$ as a function of $N$ for Owen-scrambled EC ($o$) and
Niederreiter-Xing ($*$) nets for $s=18$ and $32$.}
\end{figure}
\begin{figure}[ht]
\begin{center}
\includegraphics[height=8cm]{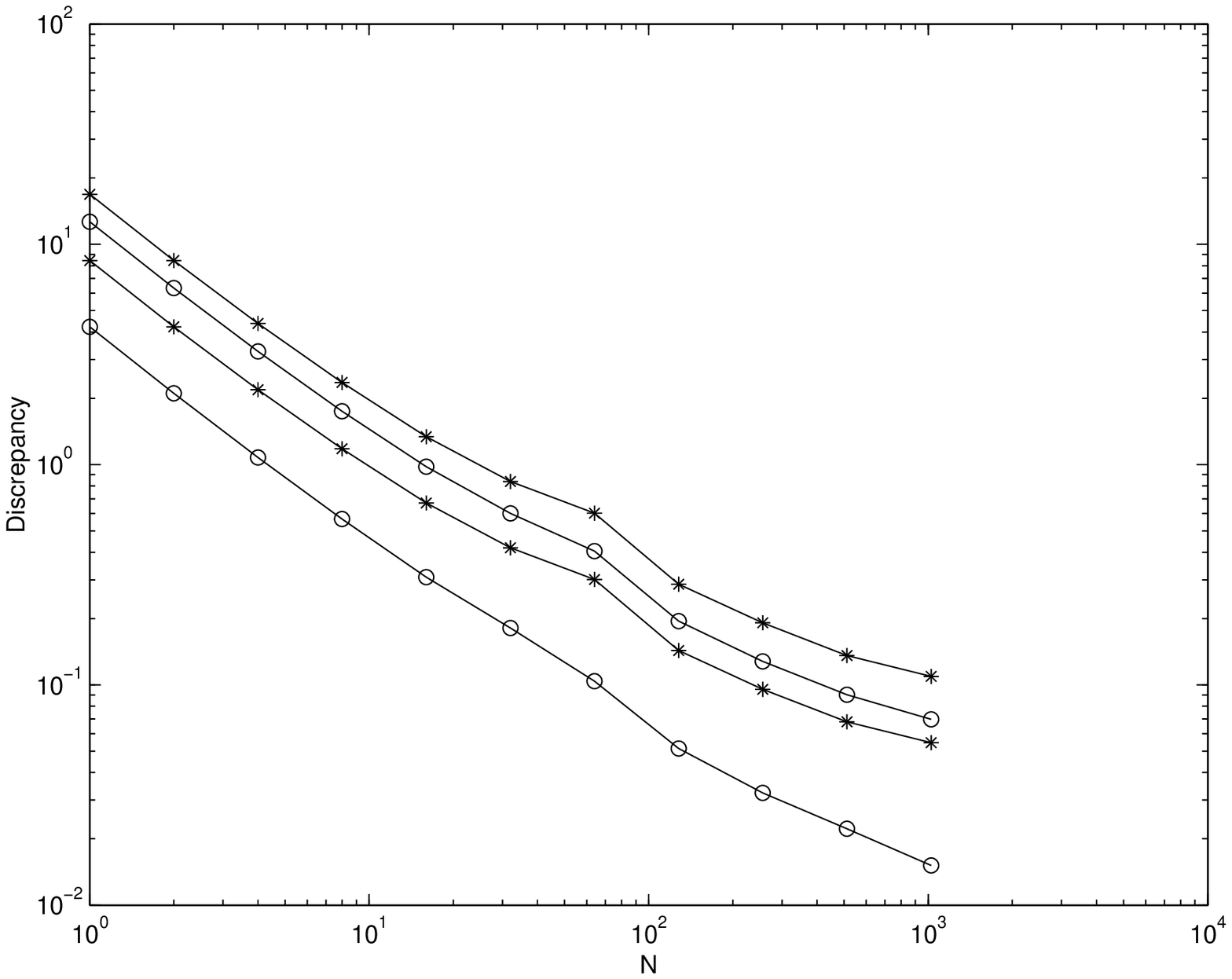}
\end{center}
\caption{\label{hnxns} The scaled root mean square discrepancy with
$\alpha = 2$ as a function of $N$ for non-scrambled EC' ($o$) and
Niederreiter-Xing ($*$) nets for $s=18$ and $32$.}
\end{figure}

The table \ref{wgt} is the comparison of $T(m,s,EC)$ 
and the updated $t$ values from the tables \cite{CLM99,Sch99b}, denoted $A$.
Here, we made comparisons between a $T(m,s,EC)$ value and the $t$ value
for $s+1$, because the EC net can be extended to $s+1$ by simply 
adding a skewed diagonal matrix to generator matrices of EC. 
It is not surprising to see that $t$ values from the references
exhibit smaller $t$ values than $T(m,s)$ values, 
since $t$ values which from \cite{CLM99,Sch99b} are the best attainable 
$t$ values for fixed $(t,m,s)$-net.
For smaller $s$ and larger $m$ the difference are large.
However propagation rules \cite{Nie99} can be adopted for further 
improvement for the EC net. Notice that as $s$ increases the 
difference becomes smaller and there are some places that $T(m,s)$ 
values and the $t$ values are same.

{\renewcommand{\baselinestretch}{1.2}
\begin{table}
    \caption{\label{wgt} Comparison of $t$ values 
from \cite{CLM99,Sch99b} and $T(m,s)$ values of
generator matrices from EC.}
\begin{center}
\scriptsize
\begin{tabular}{*{22}{c}}
  s & m &1&2&3&4&5&6&7&8&9&10&11&12&13&14&15&16&17&18&19&20\\\hline
  1 &   &*&*&*&*&*&*&*&*&*&*&*&*&*&*&*&*&*&*&*&*\\
  2 &   &*&*&*&*&*&*&*&*&*&*&*&*&*&*&*&*&*&*&*&*\\
  3 &   &*&*&*&*&*&*&*&*&*&*&*&*&*&*&*&*&*&*&*&*\\ 
  4 &   &*&*&A&A&A&2A&2A&2A&2A&2A&2A&2A&2A&2A&2A&2A&2A&2A&2A&2A\\ 
  5 &   &*&*&A&A&*&A&A&A&A&A&2A&2A&3A&3A&3A&3A&3A&3A&3A&3A\\ 
  6 &   &*&*&A&A&A&A&2A&2A&2A&2A&3A&3A&3A&3A&4A&4A&4A&4A&4A&4A\\ 
  7 &   &*&*&*&A&A&A&A&2A&2A&A&2A&2A&2A&2A&3A&3A&3A&3A&4A&4A\\ 
  8 &   &*&*&*&A&A&*&A&2A&2A&2A&3A&3A&2A&3A&4A&4A&4A&4A&4A&4A\\ 
  9 &   &*&*&*&A&A&A&A&2A&2A&2A&3A&2A&2A&3A&3A&3A&3A&4A&5A&5A\\ 
  10&   &*&*&*&A&A&A&A&A&2A&3A&3A&3A&4A&3A&3A&3A&2A&2A&3A&3A\\ 
  11&   &*&*&*&A&2A&2A&A&A&2A&2A&2A&3A&3A&4A&4A&3A&2A&2A&2A&2A\\ 
  12&   &*&*&*&A&2A&2A&A&A&2A&2A&3A&3A&3A&4A&3A&3A&2A&2A&2A&3A\\ 
  13&   &*&*&*&A&2A&2A&A&2A&2A&2A&3A&2A&3A&4A&3A&3A&2A&2A&2A&3A\\ 
  14&   &*&*&*&A&2A&2A&A&2A&2A&2A&2A&2A&3A&3A&4A&4A&5A&4A&3A&3A\\ 
  15&   &*&*&*&*&A&2A&A&2A&2A&2A&2A&2A&3A&3A&4A&4A&5A&4A&3A&3A\\ 
  16&   &*&*&*&*&A&2A&A&2A&2A&2A&2A&3A&3A&3A&3A&4A&5A&4A&3A&3A\\ 
  17&   &*&*&*&*&A&2A&A&A&2A&2A&2A&2A&3A&3A&3A&3A&4A&4A&3A&3A\\ 
  18&   &*&*&*&*&A&2A&A&A&2A&A&2A&2A&3A&3A&3A&3A&3A&3A&4A&5A\\ 
  19&   &*&*&*&*&A&2A&A&A&2A&A&2A&2A&3A&3A&3A&3A&3A&3A&4A&5A\\ 
  20&   &*&*&*&*&A&2A&A&A&2A&A&2A&2A&3A&2A&3A&3A&3A&2A&3A&4A\\ 
  21&   &*&*&*&*&A&2A&2A&A&2A&A&2A&2A&3A&2A&3A&3A&3A&2A&2A&3A\\ 
  22&   &*&*&*&*&A&2A&2A&A&2A&A&2A&2A&2A&2A&2A&3A&3A&2A&2A&2A\\ 
  23&   &*&*&*&*&A&2A&2A&A&A&A&2A&2A&2A&3A&3A&3A&3A&2A&2A&2A\\ 
  24&   &*&*&*&*&A&2A&2A&A&A&A&A&2A&2A&3A&3A&2A&2A&A&A&2A\\ 
  25&   &*&*&*&*&A&2A&2A&A&A&A&A&2A&2A&3A&3A&2A&2A&A&A&2A\\ 
  26&   &*&*&*&*&A&2A&2A&A&A&A&A&2A&2A&2A&3A&2A&2A&A&A&2A\\ 
  27&   &*&*&*&*&A&2A&2A&A&A&A&A&2A&2A&2A&3A&2A&2A&A&A&2A\\ 
  28&   &*&*&*&*&A&2A&2A&A&A&A&A&2A&2A&2A&3A&2A&2A&A&A&2A\\ 
  29&   &*&*&*&*&A&2A&2A&A&A&A&A&2A&2A&2A&3A&2A&2A&A&A&A\\ 
  30&   &*&*&*&*&A&2A&2A&A&A&A&A&2A&2A&2A&3A&2A&2A&A&A&A\\ 
  31&   &*&*&*&*&*&A&2A&A&A&A&A&2A&2A&2A&3A&2A&2A&A&A&A\\ 
  32&   &*&*&*&*&*&A&2A&A&A&*&A&2A&2A&2A&2A&A&A&*&A&A\\ 
  33&   &*&*&*&*&*&A&2A&A&A&*&A&2A&2A&2A&2A&A&A&*&A&A\\ 
  34&   &*&*&*&*&*&A&2A&A&A&*&A&2A&2A&2A&2A&2A&2A&2A&A&A\\ 
  35&   &*&*&*&*&*&A&2A&A&A&*&A&2A&2A&2A&2A&2A&2A&2A&A&A\\ 
  36&   &*&*&*&*&*&A&2A&A&A&*&A&A&2A&2A&2A&2A&2A&2A&A&A\\ 
  37&   &*&*&*&*&*&A&2A&A&A&*&A&A&2A&2A&2A&2A&2A&2A&A&A\\ 
  38&   &*&*&*&*&*&A&2A&A&A&*&A&A&2A&2A&2A&2A&2A&2A&A&A\\ 
  39&   &*&*&*&*&*&A&2A&A&A&A&2A&2A&2A&2A&2A&2A&2A&2A&A&A\\ 
  40&   &*&*&*&*&*&A&2A&2A&A&A&2A&2A&2A&2A&2A&2A&2A&2A&A&A\\\hline 
\end{tabular}
\end{center}
\end{table}
}
\end{section}

\begin{section}{Discussion}
The new EC net has a smaller $T(m,s)$ value than the Sobol' net for higher 
dimension and larger $m$. The new sequence shows 
smaller $T(m,s)$ values in overall comparing to Niederreiter-Xing sequence
which has been considered as a imbedded net.
For the chosen two generator matrices of the Niederreiter-Xing 
net, the new net
still shows smaller $T(m,s)$ values in some $m$ for $s_{max}$. Moreover 
the new net shows smaller $t$ values for lower dimensional projections 
as well. 

However there is a disadvantage of finding the generator matrices by 
optimization method, since we only can find the generator matrices
of a $(t,m,s)$-net. In contrast the finite size of 
generator matrices for the Sobol' and Niederreiter-Xing sequences can 
be extended indefinitely.

Our computation was carried out on a Unix station in C. The new generator 
matrices were obtained after less than 100 generations. 
The actual time was about 2 weeks. 
There is a problem of increasing $m$ due to the computation time,
since the computation time is more sensitive to $m$ than $s$.
Partially it is due to increasing numbers of possible solutions by a 
factor for 2. Also the computing time for $t$ sharply increases as 
$m$ increases.

Also notice that the new generator matrices are from
the initial populations which are partially taken from the
generator matrices of the Sobol' sequence. Therefore
the new matrices are more likely to be better than the 
Sobol' net which is constructed by the number theoretic method. 

Choosing an objective function is another difficult issue, 
In this thesis we consider the total sum 
of $t$ values as an objective function. This leads the new net to have 
smaller $t$ values in larger $s$ and $m$.
There are many different objective 
functions to be considered such as looking at the sum of  
$t$ values for different dimensional projections or/and looking at 
the maximum $t$ values of their projections. 

Finding the generator matrices based on the Niederreiter-Xing sequence is 
another interesting thing to try, because the generator matrices of the 
Niederreiter-Xing sequence are full matrices and it is the most recently 
developed one. Nevertheless it is easier to find 
better generator matrices based on the existing matrices rather than 
the complete new matrices.

\end{section}
\end{chapter}

\begin{chapter}{Application}

\begin{section}{Multivariate Normal Distribution}
Consider the following multivariate normal probability
\begin{equation} \label{multprob}																																																																																									
\frac{1}{\sqrt{|\bfSigma|(2\pi)^s}}\int_{[\bfa,\bfc]} 
e^{-\frac{1}{2}\bftheta^{T}\bfSigma^{-1}\bftheta}\ d\bftheta,
\end{equation}
where $\bfa$ and $\bfc$ are known s-dimensional vectors that define
the interval of integration, and $\Sigma$ is a given $s\times s$
positive definite covariance matrix.  One or more of the components of
$\bfa$ and $\bfc$ may be infinite.  Since the original form is not
well-suited for numerical quadrature, Alan Genz \cite{Gen92,Gen93}
proposed a transformation of variables that result in an integral over
the $s-1$ dimensional unit cube.  This transformation is used here.
\begin{figure}[h]
\begin{center}
\includegraphics[height=10cm]{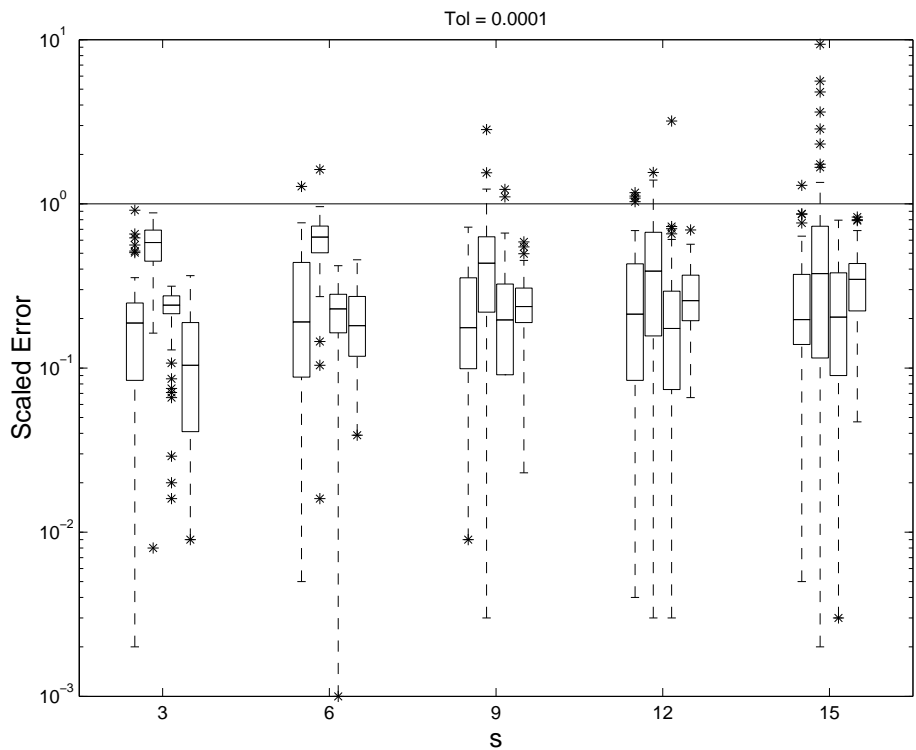}
\caption{\label{mvpe} Box and whisker plots of scaled
errors, $E/\epsilon$, for 50 randomly chosen test problems
\eqref{multprob}.  For each dimension $s$ results from left to the
four algorithms as listed in the text.}
\end{center}
\end{figure}
\begin{figure}[h]
\begin{center}
\includegraphics[height=10cm]{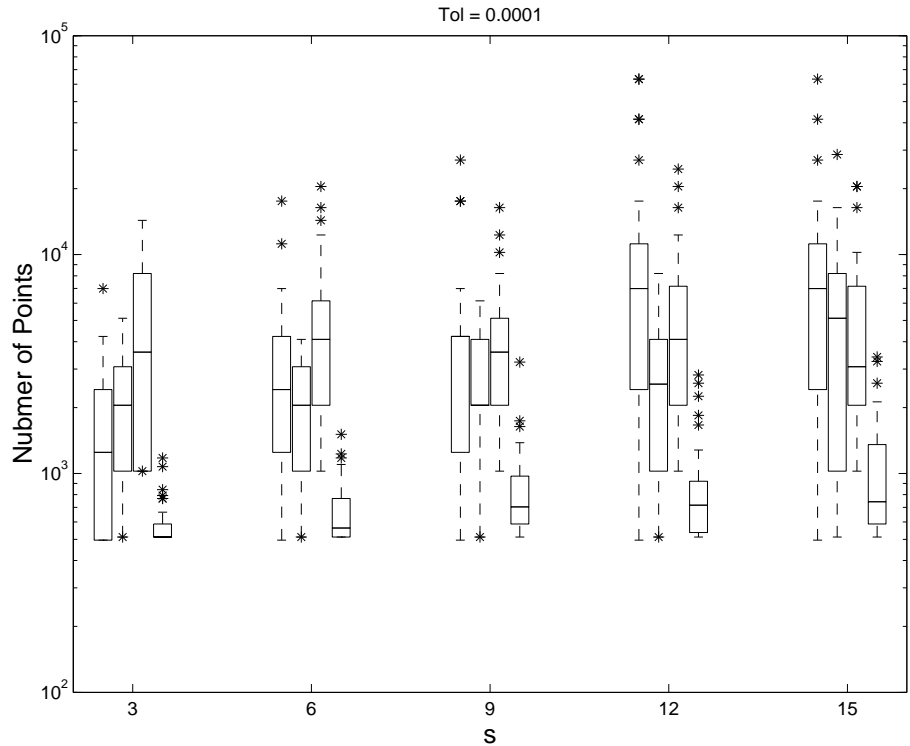}
\caption{\label{mvpp} Box and whisker plots of the computation times
in seconds for 50 randomly chosen test problems \eqref{multprob}.  For
each dimension $s$ results from left to the four algorithms as listed
in the text.}
\end{center}
\end{figure}
\begin{figure}[h]
\begin{center}
\includegraphics[height=10cm]{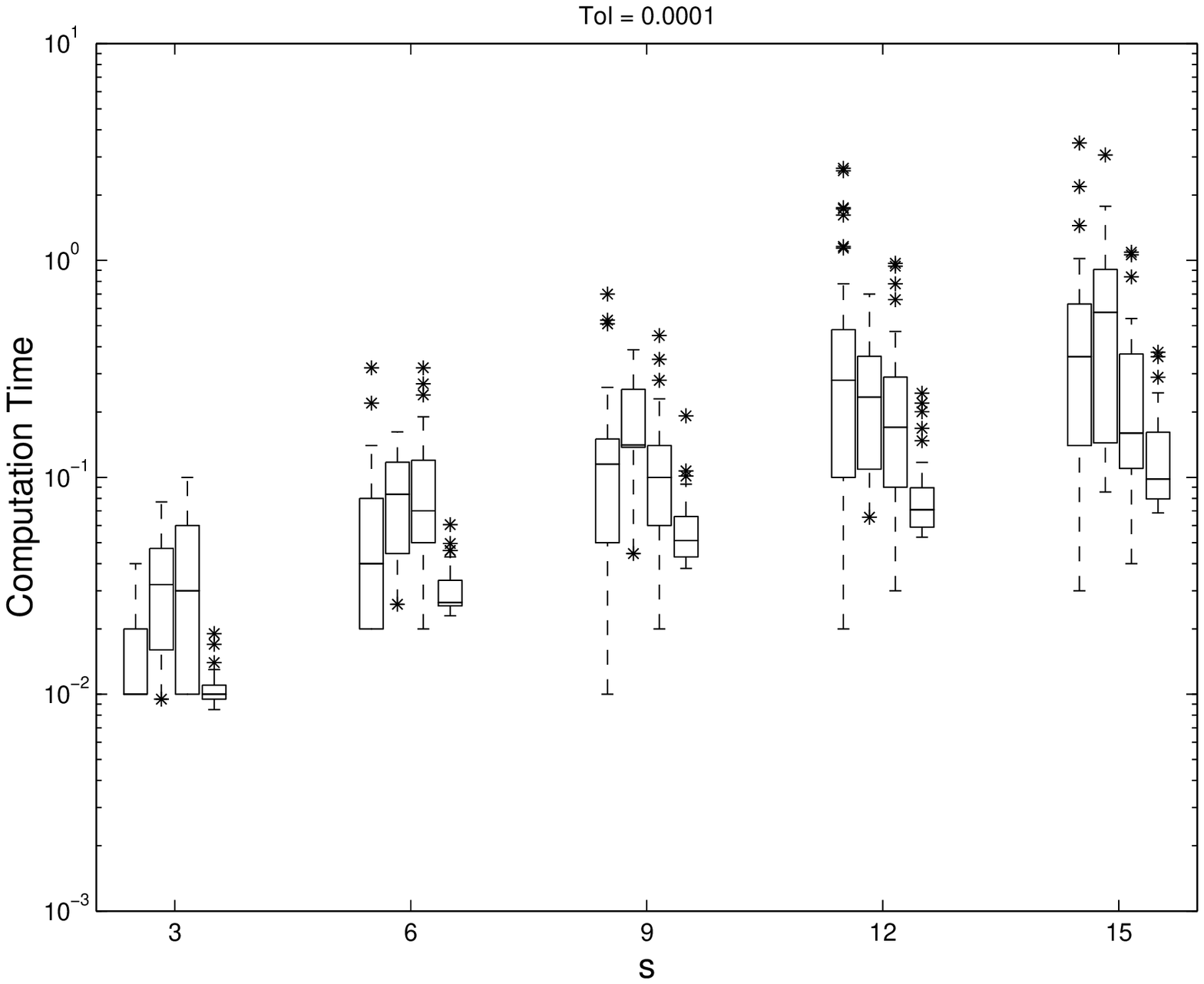}
\caption{\label{mvpt} Box and whisker plots of the number of points
used for solving 50 randomly chosen test problems \eqref{multprob}. 
For each dimension $s$ results from left to right, the four algorithms as
listed in the text.}
\end{center}
\end{figure}

The choices of the parameters in \eqref{multprob} are made as in 
\cite{Gen93,HicHon97a}:
\begin{gather*}
a_1 = \cdots = a_s = -\infty, \qquad c_j \text{ i.i.d. uniformly on }
[0,\sqrt{s}], \\
\Sigma \text{ generated randomly as in \cite{Gen93}}.
\end{gather*}

Four different types of algorithms are compared for this problem:
\begin{enumerate}
    
    \renewcommand{\labelenumi}{\roman{enumi}.}
    \item the adaptive algorithm MVNDNT of Alan Genz, which can be
    found at \\ {\tt http://www.sci.wsu.edu/math/faculty/genz/software/mvndstpack.f}

    \item a Korobov rank-1 lattice rule implemented as part of the NAG
    library,

    \item a randomly-shifted extensible rank-1 lattice sequence of
    Korobov type \cite{HicHon97a,HicEtal00} with generator vector $(1,
    17797, \ldots, 17797^{s-1})$, and

    \item the scrambled Sobol' sequence described here.
\end{enumerate}
The periodizing transformation $x_{j}' = |2x_j-1|$, $j=1, \ldots, s$
has been applied to the integrand for the second and third algorithms ii.\ and
iii., as it appears to increase the accuracy of these two
algorithms.  The computations were carried out in Fortran on an Unix
workstation in double precision.  An absolute error tolerance of
$\varepsilon = 10^{-4}$ was chosen and compared to the absolute error
$E$.  Error estimation for each algorithm is described in 
\cite{HicEtal00}.
Because the true value of the integral is unknown for this
test problem, the Korobov algorithm with $\varepsilon = 10^{-8}$ was
used to be ``exact" value for computing error.  For each $s$, 50 test
problems were generated randomly.  Figure \ref{mvpe}, \ref{mvpp}, and
\ref{mvpt} show the box and whisker plots of the scaled absolute error
$E/\epsilon$, the number of points used, and the computation time in
seconds for the four methods with various dimensions.  The boxes are
divided by the median and contain the middle half of the values.  The
whiskers show the full range of values and the outliers are plotted as
$*$.

The scaled error indicates how conservative the error estimate is. 
Ideally the scaled error should be close to but not more than one.  If
the scaled error is greater than one, then the error estimation is not
careful enough.  However, if the scaled error is much smaller than one,
then the algorithm is wasting time by being too conservative.  From
the Figures \ref{mvpe}, \ref{mvpp} and \ref{mvpt} all four methods perform
reasonably well in error estimation.  The scrambled Sobol' 
points, however, use fewer function evaluations and less computation 
time.  Also, their performance is less dependent on the dimension.

\end{section}

\begin{section}{Physics Problem}
The following multidimensional integral arising in physics problems is
considered by Keister \cite{Kei96}:
\begin{equation} \label{physprob}
\int_{\vec R^s} \cos(\|\bfx\|)e^{-\|x\|^2}\ d\bfx = \pi^{s/2}\int_{[0,1]^s} \cos
\left( \sqrt{\sum_{j = 1}^s\frac{[\Phi^{-1}(\bfy_j)^2]} {2}} \right)\ 
d\bfy,
\end{equation} 
where $\|\cdot\|$ denotes the Euclidian norm in $\vec R^s$, and $\Phi$
denotes the standard multivariate Gaussian distribution function. 
Keister provides an exact formula for the answer and compared the
quadrature methods of McNamee and Stenger \cite{McNSte67} and Genz and
Patterson \cite{Gen82,Pat68} for evaluating this integral. 
Papageorgiou and Traub \cite{PapTra97} applied the generalized Faure
sequence from their FINDER library to this problem.  The exact value of
the integral is reported in \cite{Kei96,PapTra97}.

\begin{figure}[h]
\begin{center}
\includegraphics[height=10cm]{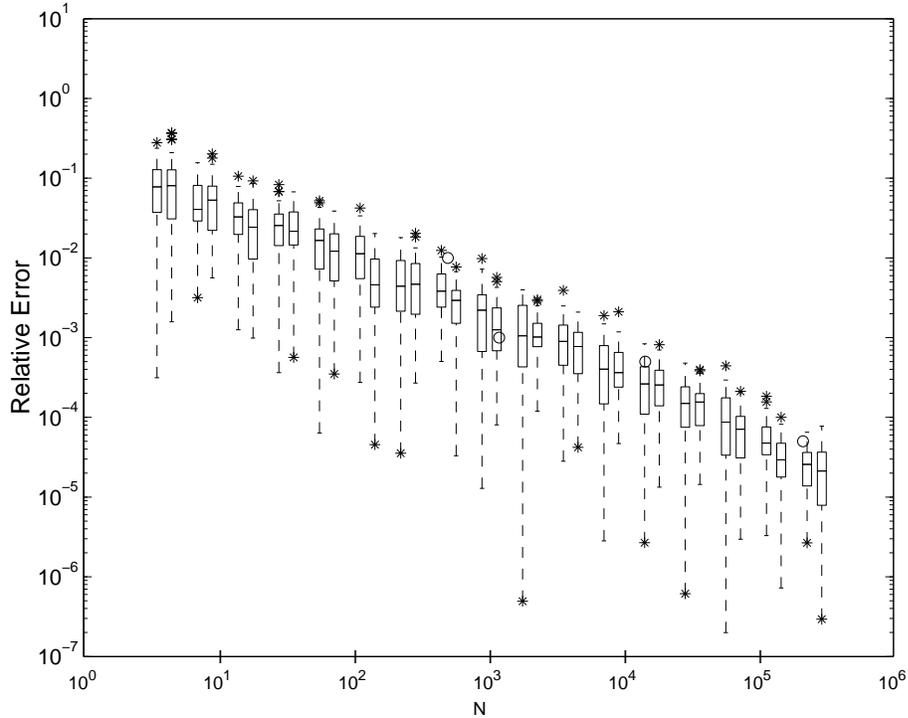}
\caption{\label{phyp} The relative error obtained in approximating
\eqref{physprob} for $s= 25$ for FINDER's generalized Faure 
algorithm $(\circ)$.  The box and whisker
plots show the performance of the randomly scrambled Sobol' points and 
randomly shifted extensible rank-1 lattice,
for the values of $N = 4,8, \ldots ,2^{18}$.}
\end{center}
\end{figure}

Since the methods of McNamee and Stenger and Genz and Patterson
performed much worse than FINDER's generalized Faure sequence, they
are ignored here.  Instead we compare the scrambled Sobol' sequences 
implemented here to the performance of FINDER
and the extensible lattice rules described before.  To be consistent
with the results reported in \cite{PapTra97} the relative error is
computed as a function of $N$.  Figure \ref{phyp} shows the results of
the numerical experiments for dimension $s=25$.  Box and Whisker plots
show how well 50 randomized nets and lattice rules perform.  The three
algorithms appear to be quite competitive to each other.  In some
cases randomized Sobol'  performs better than the other two sequences.

\end{section}

\begin{section}{Multidimensional Integration}
Consider the following multidimensional integration problem:
\begin{equation} \label{mip}																																																																																									
\int_{\Cs} \prod_{j = 1}^s [1+a_j(y_j-\frac{1}{2})] d\bfy,
\end{equation}
where the exact value of the integration is 1.
We numerically compute this problem for the case $a_j = 0.4+\frac{j}{10}$ where 
$j = 1,\cdots,s$ with scrambled Sobol' and EC nets. 
Figure \ref{sec} show the root mean square relative error of
(5.3) for the scrambled Sobol' and scrambled EC nets.
Figure \ref{sec} plots the root mean square relative error computed as a 
function of $N$ for $s = 14$ and $24$. 
For scrambling 100 different replications are made.
Figure \ref{sec} shows that the scrambled 
EC sequence has smaller relative error than the scrambled
Sobol' sequence for all choices of $N$.

\begin{figure}[ht]
\begin{center}
\includegraphics[height=8cm]{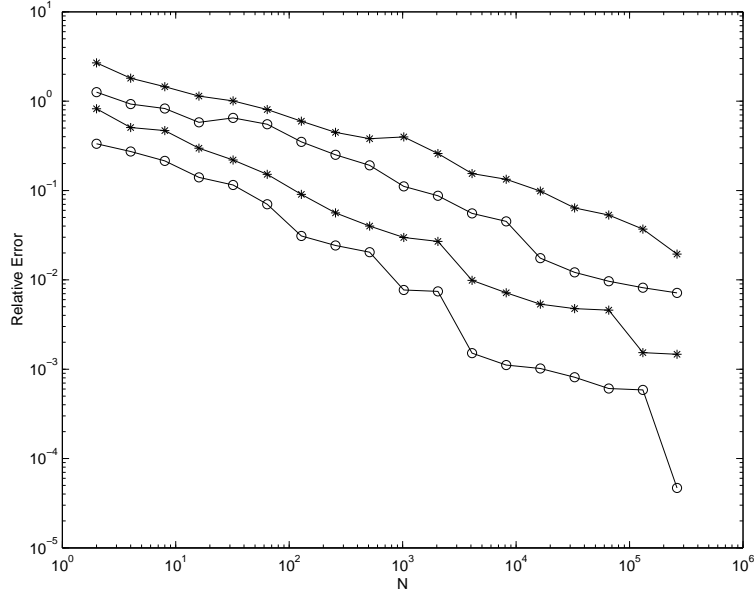}
\end{center}
\caption{\label{sec} The relative error of the numerical 
integration of (5.3) for scrambled Sobol($*$) and scrambled
EC($o$) nets for $s = 14$ and $24$.}
\end{figure}

Numerical experiments on multidimensional integrations are performed 
by using Genz' test function package \cite{Gen87}. 
We test his six integral families, 
namely Oscillatory, Product Peak, Corner Peak, Gaussian, $C^0$ function,
and Discontinuous for five different sequences.
The compared sequences are Niederreiter-Xing, scrambled Niederreiter-Xing, 
EC, scrambled EC, and imbedded Niederreiter-Xing sequences. For imbedded
Niederreiter-Xing sequence which we choose the
generator matrices ``nxs20m30". The choices of parameters are made 
as in \cite{Gen87}. For scrambled sequences, 20 different replications 
are performed. The digit accuracy is obtained by computing
{\tt{max(log(relative error),0))}}. Therefore if the digit accuracy is 
zero then the relative error is greater than 1. The larger value of 
the digit accuracy implies the better performance.
From Figures \ref{gz1}-\ref{gz6} it is shown that scrambled nets perform 
almost always the best. Considering non-scrambled sequences
the EC sequence shows the better performance as the dimension of 
the problem increases. The Niederreiter-Xing sequence tends to perform 
well in lower dimension, and the imbedded Niederreiter-Xing sequence always performs 
worse than the Niederreiter-Xing sequence. This is expected because the 
Niederreiter-Xing sequence are not designed as imbedded sequences.
We also test the problems with Sobol' and scrambled Sobol' sequences.
The performances of Sobol' and EC sequences are nearly the same.

\begin{figure}[ht]
\begin{center}
\includegraphics[height=8cm]{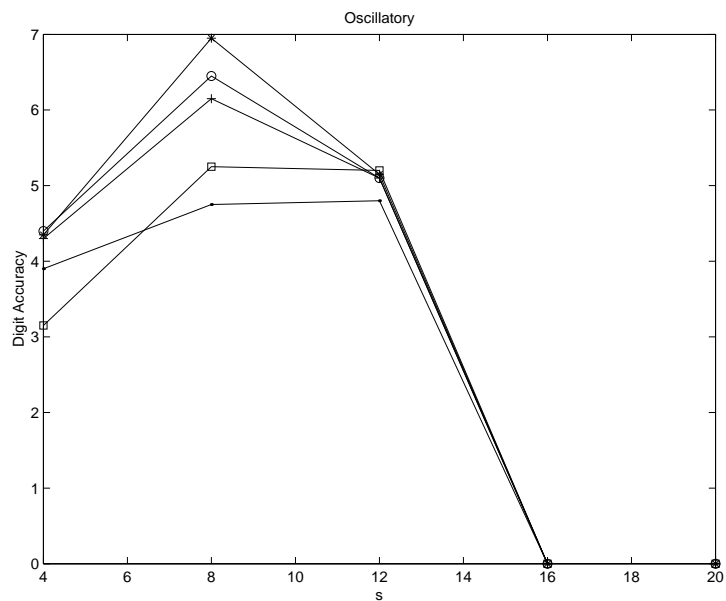}
\end{center}
\caption{\label{gz1} The digit accuracy of the numerical 
integration of Oscillatory functions for EC ($.$), 
scrambled EC ($o$),
Niederreiter-Xing ($+$), scrambled Niederreiter-Xing ($*$), 
and imbedded Niederreiter-Xing ($\square$) nets.}
\end{figure}
\begin{figure}[ht]
\begin{center}
\includegraphics[height=8cm]{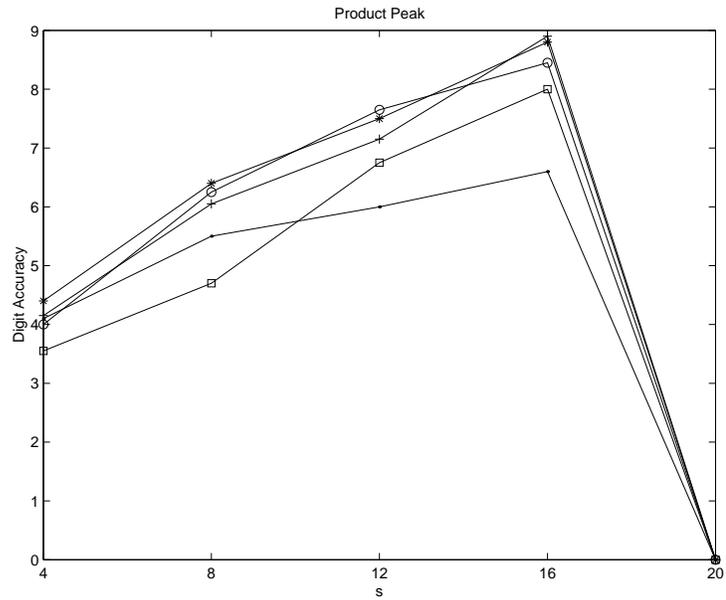}
\end{center}
\caption{\label{gz2} Same as \ref{gz1} for the numerical 
integration of Product Peak functions.}
\end{figure}
\begin{figure}[ht]
\begin{center}
\includegraphics[height=8cm]{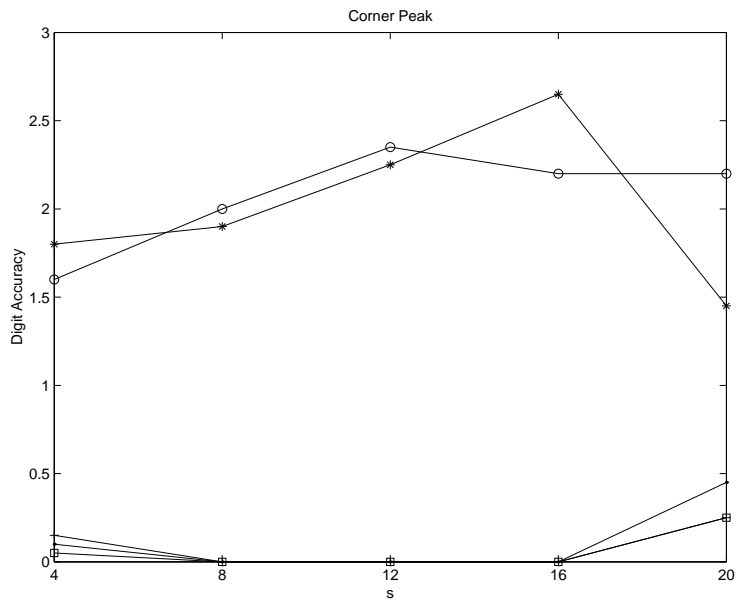}
\end{center}
\caption{\label{gz3} Same as \ref{gz1} for the numerical 
integration of Corner Peak functions. }
\end{figure}
\begin{figure}[ht]
\begin{center}
\includegraphics[height=8cm]{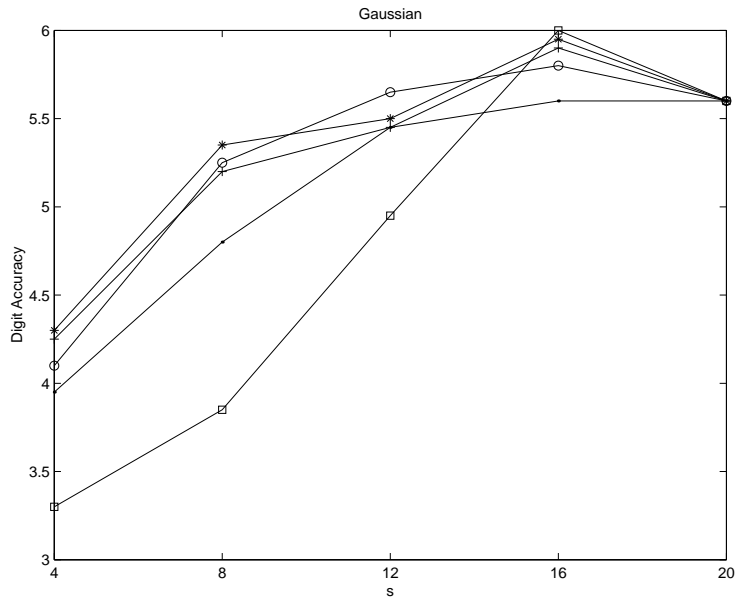}
\end{center}
\caption{\label{gz4} Same as \ref{gz1} for the numerical 
integration of Gaussian functions. }
\end{figure}
\begin{figure}[ht]
\begin{center}
\includegraphics[height=8cm]{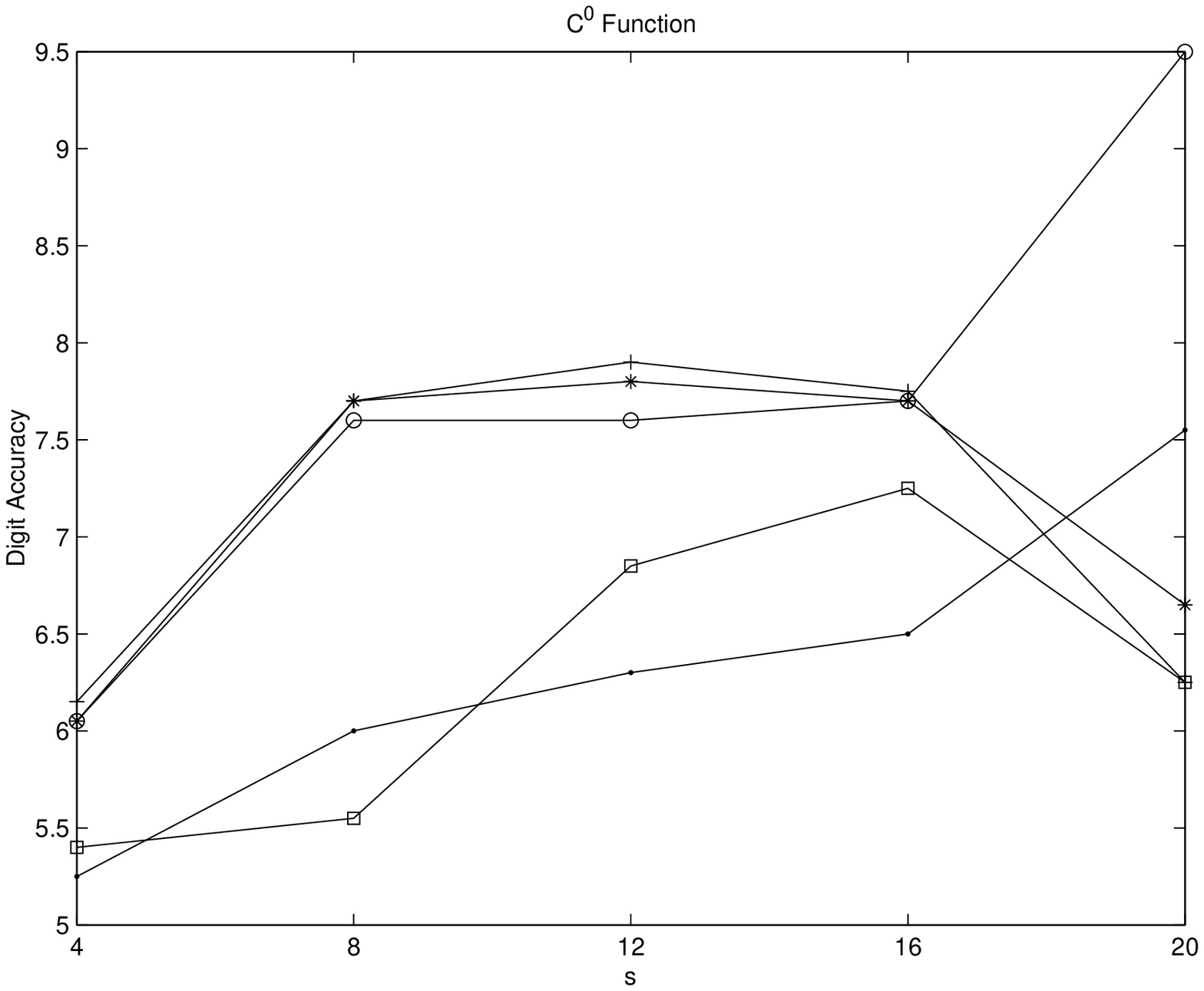}
\end{center}
\caption{\label{gz5}  Same as \ref{gz1} for the numerical 
integration of $C^o$ functions. }
\end{figure}
\begin{figure}[ht]
\begin{center}
\includegraphics[height=8cm]{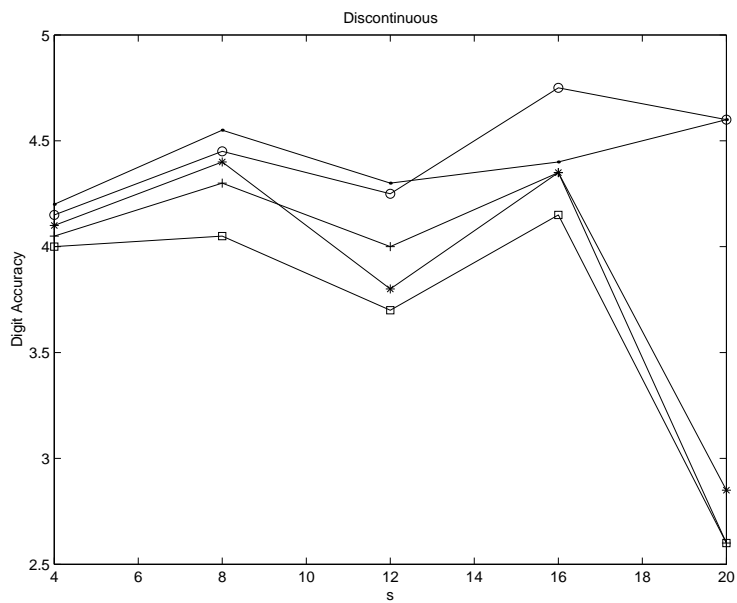}
\end{center}
\caption{\label{gz6} Same as \ref{gz1} for the numerical 
integration of Discontinuous functions. }
\end{figure}
\end{section}
\end{chapter}

\newpage
\begin{chapter}{Conclusion}
In Chapter 2 the randomization of Owen and of Faure and Tezuka have been
implemented for the best known digital $(t,s)$-sequences.  Both types
of randomization involve multiplying the original generator matrices by
randomly chosen matrices and adding random digital shifts.  Base 2
sequences can be generated faster by taking the advantage of the
computer's binary arithmetic system.  Using a gray code also speeds up
the generation process.  The cost of producing a scrambled sequence is
about twice the cost of producing a non-scrambled one.

Scrambled sequences often have smaller discrepancies than their
non-scrambled counterparts.  Moreover, random scrambling facilitates
error estimation.  On test problems taken from the literature 
scrambled digital sequences performed as well as or better 
than other quadrature methods (see Chapter 5).

In Chapter 3 the empirical distribution of the square discrepancy of scrambled
digital nets has been fit to a mixture of chi-square random variables
as suggested by the central limit theorem by Loh \cite{Loh01}.  Apart
from some technical difficulties that have been discussed there is
a good fit between the empirical and theoretical distributions.
Although the digital scrambling schemes used in
\cite{Mat98,HonHic00a,YueHic02a} gives the same mean square
discrepancy as Owen's original scheme, the distribution of the mean
square discrepancies varies somewhat.  It seems that the square
discrepancy has a smaller variance for Owen's scheme.  It is shown
how one may study the uniformity of lower dimensional projections of
sets of points by decomposing the square discrepancy into pieces,
$D^2_j$.  In particular, Sobol' nets seem to have better
uniformity for low dimensional projections than the more recent
Niederreiter-Xing nets.

In Chapter 4 new digital $(t,m,s)$-nets are 
generated by using the evolutionary computation method.  
The new net shows better equidistribution
properties for large $s$ and $m$ compared to the the Sobol' sequence.
We also compare the exact quality parameters $t$ of the new nets
with the Niederreiter-Xing nets, where we consider the Niederreiter-Xing
nets as an imbedded net. It is shown that the new net has overall
smaller $t$ values for different dimensional projections.
For the $s = 18$ and $32$, the new net has better $t$ values than 
the Niederreiter-Xing net in some $m$ for $s_{max}$. 
By testing the new nets with multidimensional integration
problems, we find that the new nets performs better than Sobol' and 
Niederreiter-Xing nets for certain problems.
Chapter 5 shows some promising results that
it is possible to improve existing nets by the optimization method.
 
\end{chapter}

\cleardoublepage
\phantomsection
\addcontentsline{toc}{chapter}{Bibliography}

\newcommand{\etalchar}[1]{$^{#1}$}
\def\Ignore#1{}\def\Ignore#1{}\def\notesupp#1{}

\end{document}